\documentclass[12pt]{article}
\usepackage{graphicx}
\usepackage[percent]{overpic}
\usepackage{amsmath,amscd}
\usepackage{amsthm,amssymb}
\newtheorem{theorem}{Theorem}[section]
\newtheorem{prop}{Proposition}[section]
\newtheorem{cor}{Corollary}[section]
\newtheorem{lemma}{Lemma}[section]

\newtheorem{definition}{Definition}[section]
\title{Moduli spaces for Lam\'e functions and Abelian differentials
of the second kind}

\def\C{\mathbf{C}}
\def\bC{\mathbf{\overline{C}}}
\def\CP{\mathbf{CP}}
\def\Z{\mathbf{Z}}

\def\A{\mathbf{A}}
\def\V{\mathbf{V}}
\def\id{\mathrm{id}}
\def\Sph{\mathrm{Sph}}
\def\R{\mathbf{R}}
\def\mod{\mathrm{mod\,}}
\def\LW{\mathrm{LW}}
\def\Rea{\mathrm{Re}\,}
\def\tr{\mathrm{tr}\, }
\def\Ima{\mathrm{Im}\, }
\def\Lame{\mathrm{Lame}}
\def\T{\mathbf{T}}
\def\L{\mathbf{L}}
\def\H{\mathbf{H}}
\def\[{\lfloor}
\def\]{\rfloor}
\author{Alexandre Eremenko\thanks{Supported by NSF grant
DMS-1665115.},
Andrei Gabrielov${}^*$,\\ Gabriele Mondello\thanks{Partially supported by
INdAM research group GNSAGA and by grant PRIN 2017 ``Moduli and Lie theory''.}\; and
Dmitri Panov\thanks{Supported by EPSRC grant EP/S035788/1.}}
\begin{document}
\maketitle
\begin{abstract}
The topology of the moduli space for Lam\'e functions
of degree $m$ is determined: this is a Riemann surface which consists
of two connected components when $m\geq 2$; we find the Euler
characteristics and genera of these components. 
As a corollary we prove a conjecture of R. Maier on degrees of Cohn's
polynomials.
These results are obtained
with the help of a geometric description of these Riemann surfaces,
as quotients
of the moduli spaces for certain singular flat triangles.

An application is given to the study of
metrics of constant positive curvature
with one conic singularity with the angle $2\pi(2m+1)$ on a torus. We show
that the degeneration locus of such metrics is 
a union of smooth analytic curves and we enumerate these curves. 


\end{abstract}

\section{Introduction }
\subsection{Statement of main results}
\label{intro}

Consider an elliptic curve in the form of Weierstrass
\begin{equation}\label{W}
u^2=4x^3-g_2x-g_3,\quad g_2^3-27g_3^2\neq 0.
\end{equation}
{\em Lam\'e equation} (algebraic version) is the second order differential
equation
\begin{equation}\label{lame1}
\left(\left( u\frac{ d}{dx}\right)^2-m(m+1)x-\lambda\right)w=0,
\end{equation}
with respect to a function $w(x)$, where $m\geq 0$ is an integer. The more familiar form of the Lam\'e equation is
\begin{equation}\label{lame2}
\frac{d^2W}{dz^2}-\left(m(m+1)\wp(z)+\lambda\right) W=0,
\end{equation}
which is obtained from (\ref{lame1}) by the change of the
independent variable $x=\wp(z),\; u=\wp'(z)$, $W(z)=w(\wp(z))$.
Here $\wp$ is the Weierstrass
function of the lattice $\Lambda$ with invariants
$$g_2=60\sum_{\omega\in\Lambda\backslash\{0\}}\omega^{-4},
\quad g_3=140\sum_{\omega\in\Lambda\backslash\{0\}}\omega^{-6}.$$

Changing the variables in (\ref{lame1}) to $x_1=kx,\; w_1(x_1)=w(x_1/k)$,
we obtain a new equation (\ref{lame1}) with parameters
\begin{equation}\label{scaling}
(k\lambda,k^2g_2,k^3g_3),\quad k\in\C^*.
\end{equation}
Two equations obtained by such a change of the variables
are called {\em equivalent.}
The set of equivalence
classes is the {\em moduli space for Lam\'e equations}, $\Lame_m$.
The quotient of $\C^3\backslash\{0\}$ by the $\C^*$ action
(\ref{scaling}) is the weighted projective space $\CP(1,2,3)$, see
for example %\cite[Section 5.1]{Reid},
\cite{H}, and
$\Lame_m$ is obtained from it by removing the curve $g_2^3-27g_3^2=0$.
It parametrizes projective structures on tori with one conic singularity
with angle $2\pi(2m+1)$. For a general
discussion of projective structures
with conic singularities we refer to \cite{Loray}.
$\Lame_m$ has two singularities corresponding to two points with non-trivial
stabilizers: $(0,1,0)$ and $(0,0,1)$, which are associated with
square and hexagonal
tori in (\ref{W}).

The homogeneous function 
$$J=\frac{g_2^3}{g_2^3-27g_3^2}$$
is called the {\em absolute invariant} of the elliptic curve (\ref{W});
it defines a map $\pi_m:\Lame_m\to\C_J$ which is called the
{\em forgetful map}.

A {\em Lam\'e function} is a non-trivial solution $w$ of (\ref{lame1}) with
the property that $w^2$ is a polynomial in $x$. It is easy to see
that the degree of this polynomial must be $m$.
For given
$g_2,g_3$, such solutions exist if and only if the {\em accessory parameter}
$\lambda$ satisfies
an algebraic equation 
\begin{equation}
\label{F}
F_m(\lambda,g_2,g_3)=0,
\end{equation}
where $F_m$ is a polynomial in three variables, monic with
respect to $\lambda$,
which is invariant under the $\C^*$ action (\ref{scaling}).
For each $\lambda$ satisfying (\ref{F}), the Lam\'e function
is unique up to a constant factor, \cite[23.41]{WW}. 

Two Lam\'e functions
are {\em equivalent} if $w(x)=cw_1(kx),\; c,k\in\C^*$, and the set of equivalence
classes is the {\em moduli space for Lam\'e functions} $\L_m\subset\Lame_m$.
We consider this space $\L_m$ as an abstract Riemann surface which is obtained
by taking the quotient of the surface
$$\{(\lambda,g_2,g_3)\in\C^3:F_m(\lambda,g_2,g_3)=0,\; g_2^3-27g_3^2\neq 0\}$$
by the $\C^*$ action~(\ref{scaling}).
It is the normalization of the curve defined by equation (\ref{F})
in $\CP(1,2,3)\backslash\{ g_2^3-27g_3^2=0\}$.
Since $F_m$ in (\ref{F}) is a monic polynomial, the forgetful map $\L_m\to\C_J$
is proper.

Equations (\ref{F}) for $m\leq 8$ are
explicitly written
in Table 3 in \cite{Maier}, see also Appendix to the present paper. 
Maier calls $F_m$ the
{\em Lam\'e spectral polynomials}.

Lam\'e equation and Lam\'e functions
have long history going back to the work of Gabriel Lam\'e
\cite{Lame1,Lame2}, and they have been intensively studied ever since,
because of their
importance for mathematical physics. Good reference
for the classical work is \cite{WW},
%$(see also \cite{Heine} and \cite{Halphen}), 
and a modern survey is contained in
the first three sections of \cite{Maier}.

Most of the classical work on Lam\'e functions was
concentrated on the real case with positive discriminant 
($g_2,g_3,\lambda$ are real and $g_2^3-27g_3^2>0$), and we are not aware of any
systematic study of 
general properties of $F_m$ and~$\L_m$.

In this paper we determine the topology of the Riemann surfaces $\L_m$.
\def\N{\mathbf{N}}
To state our main result,
we recall the notion of a $2$-dimensional orbifold 
\cite{Milnor},\cite{Thurston}.
It is a compact Riemann surface $S$ with a function
$n:S\to\N_{\geq 1}\cup\{\infty\}$,
where $\N_{\geq 1}$ is
the set of positive integers, such that $n(x)=1$ for all points, except
finitely many {\em orbifold points} of order $n(x)>1$.
Points with $n(x)=\infty$ are interpreted as
punctures.
For example, the moduli space for tori $\C_J$ has a natural orbifold
structure with $S=\bC$ and  
three orbifold points:
$n(0)=3$, $n(1)=2$, and $n(\infty)=\infty$. $J=0$ corresponds to
the hexagonal torus, and $J=1$ to the square torus.

The orbifold Euler characteristic is defined as
$$\chi^O=\chi(S)-\sum_{x\in S} \left(1-1/n(x)\right),$$
where $\chi(S)=2-2g(S)$ is the ordinary Euler characteristic of
the underlying compact surface,
and $g$ is its genus. Here we follow \cite{Thurston}.
Notice that in \cite{EMP} a different definition
is used: the Euler characteristic in \cite{EMP} is smaller
by a factor of $2$.
A ramified covering $f:S_1\to S_2$ is called an {\em orbifold map}
if $n(f(x))$ divides $n(x)\deg_x(f)$ for all $x\in S_1$, and an
{\em orbifold covering} if
\begin{equation}\label{ocov}
 n(f(x))=n(x)\deg_x f\quad\mbox{for all}\quad x\in S_1. 
\end{equation}
If $f$ is an orbifold covering, the Riemann-Hurwitz formula
gives
\begin{equation}\label{rh}
\chi^O(S_1)=\deg(f)\chi^O(S_2).
\end{equation}
We introduce the following functions of non-negative integers $m$.
\begin{eqnarray}\label{degFI}
d_m^I&:=&\left\{\begin{array}{ll}m/2+1,& m\equiv0\;
(\mod 2)\\
(m-1)/2,& m\equiv 1\; (\mod 2),\end{array}\right.\\
\label{degFII}
d_m^{II}&:=&3\lceil m/2\rceil.
\end{eqnarray}
It is easy to see that $d_m^I$ and $d_m^{II}$ always have
opposite parity: when $m\in\{0,3\}\;(\mod 4),$ $d_m^I$ is odd
and $d_m^{II}$ is even; when $m\in\{1,2\}\;(\mod 4),$ $d_m^I$ is even
and $d_m^{II}$ is odd.
\begin{eqnarray}\label{epsilon0}
\epsilon_0&=&\left\{\begin{array}{ll}0,&\mbox{if}\; m\equiv1\;(\mod 3),\\
1&\mbox{otherwise,}\end{array}\right.\\ \label{epsilon2}
\epsilon_1&=&\left\{\begin{array}{ll}0,&\mbox{if}\; m\in\{1,2\}\;(\mod 4),\\
1&\mbox{otherwise.}\end{array}\right.
\end{eqnarray}
One can restate the definitions of $\epsilon_j$ as follows:
$\epsilon_0=0$ if and only if $d_m^I$ is divisible by $3$,
and $\epsilon_1=0$ if and only if $d_m^I$ is even.

Our main result is the following.

\begin{theorem}\label{theorem1}
When $m\geq 2$, $\L_m$ is a Riemann surface consisting
of two connected components which we call $\L_m^I$
and $\L_m^{II}$, while
$\L_0$ and $\L_1$ are connected: $\L_0=\L_0^I$ and
$\L_1=\L_1^{II}$.

The Riemann surface $\L_m$ has a natural orbifold structure with $\epsilon_0$
orbifold points of order 3 on $\L_m^I$, and one orbifold point
of order 2 which belongs to $\L_m^I$ when $\epsilon_1=1$ and
to $\L_m^{II}$ otherwise.
The component $\L_m^I$ has $d_m^I$ punctures, and the component $\L_m^{II}$ has
$2d_m^{II}/3=2\lceil m/2\rceil$
punctures.

The restrictions of the forgetful map to these components are orbifold maps
and their degrees are 
$d_m^I$ and $d_m^{II}$. 
The orbifold Euler characteristics are
\begin{equation}\label{euler}
\chi^O(\L_m^I)=-(d_m^I)^2/6\quad\mbox{and}\quad\chi^O(\L_m^{II})=
-(d_m^{II})^2/18.
\end{equation}
\end{theorem}

{\bf Remark.} 
Ordinary Euler characteristics $\chi$ are
obtained from the $\chi^O$ by adding the {\em orbifold corrections}
which in our case are
\begin{equation}\label{ocorr}
E^I=(4\epsilon_0+3\epsilon_1)/6\quad\mbox{and}\quad
E^{II}=(1-\epsilon_1)/2.
\end{equation}
Euler characteristics can be expressed as functions of $m$
rather than $d$, see Appendix.

It is well known that for $m\geq 2$ each polynomial $F_m$ factors
into four factors in
$\C[\lambda,e_1,e_2,e_3]$, where $e_i$ are related to $g_2,g_3$ by the equation
$$4x^3-g_2x-g_3=4(x-e_1)(x-e_2)(x-e_3),$$
see for example \cite[Thm. 2]{T} for an explicit statement of this.
However there is no discussion of irreducibility of these four factors
in the literature. Our theorem says that for $m\geq 2$, $F_m$ has exactly
two irreducible factors in $\C[\lambda,g_2,g_3]$ and implies that
the four factors of $F_m$ in $\C[\lambda,e_1,e_2,e_3]$ are irreducible.

Theorem \ref{theorem1} implies
the formulas for the genera of $\L_m^K,\; K\in\{ I,II\}$ in terms of $m$
or $d_m^K$ which are given in the Appendix.

We give several applications of Theorem \ref{theorem1}.
As a first application, we prove that the two irreducible
components of the surface $F_m(\lambda,g_2,g_3)=0$ in (\ref{F}) 
have no singularities in $\C^3$ except the lines $(0,t,0)$ and $(0,0,t)$.

To obtain a non-singular curve in $\CP^2$ parametrizing $\L_m$, we use
Legendre's family of elliptic curves
\begin{equation}\label{le}
v^2=z(z-1)(z-a),\quad a\in\C_a:=\C\backslash\{0,1\}.
\end{equation}
For the $J$-invariant of this curve we have
\begin{equation}\label{Ja}
J=\psi(a)=\frac{4}{27}\frac{(a^2-a+1)^3}{a^2(1-a)^2}.
\end{equation}
If $\bC_J$ is considered as an orbifold with $n(0)=3$, $n(1)=2$,
$n(\infty)=\infty$,
and in $\bC_a$ we set 
$n(a)=\infty$ for $a\in\{0,1,\infty\}$, and $n(a)=1$ otherwise,
then $a\mapsto \psi(a)$ is an orbifold
covering.

The form of the Lam\'e equation corresponding to (\ref{le}) is 
\begin{equation}\label{legendre}
Py''+\frac{1}{2}P'y'-((m(m+1)z+B)y=0,\quad P(z)=4z(z-1)(z-a),
%4z(z-1)(z-a)y''+(6z^2-4(\lambda+1)z+2a)y'-(m(m+1)z+B)y=0,
\end{equation}
where the accessory parameter
$B$ is an affine function of $\lambda$ (see Section \ref{nons} for the
details of this transformation). Lam\'e functions correspond to non-trivial
solutions $y$ of (\ref{legendre}) such that
$y^2$ is a polynomial in $z$ of degree $m$. For such a solution to exist, a
polynomial equation
\begin{equation}\label{H}
H_m(B,a)=0
\end{equation}
must be satisfied. The Riemann surface defined by this equation
will be denoted by $\H_m$. It is the normalization of
the algebraic curve $H_m(B,a)=0$ in $\C_a\times\C_B$.
For $m\geq 2$, we will show that
it consists of four irreducible components,
$\H_m^j,\; j\in\{0,1,2,3\}.$
These components are defined as follows: $\H_m^0$
corresponds to polynomial solutions $y$ when $m$ is even and
to solutions of the form $y=Q\sqrt{P}$ when $m$ is odd. The other
three components $j\in\{1,2,3\}$ correspond to solutions
$$Q(z)\sqrt{z},\quad Q(z)\sqrt{z-1},\quad Q(z)\sqrt{z-a},\quad
\mbox{when}\quad m\quad\mbox{is odd},$$
and 
$$Q(z)\sqrt{(z-1)(z-a)},\quad Q(z)\sqrt{z(z-a)},\quad Q(z)\sqrt{z(z-1)},$$
when $m$ is even,
where $Q$ is a polynomial. The forgetful maps $\sigma_m^j:\H_m^j\to\C_a$
are defined by $(B,a)\mapsto a$.

The polynomial $H_m$ is a product of four 
factors $H_m^j$, and we have ramified coverings
$\Psi_m^0:\H_m^0\to \L_m^I,$
and $\Psi_m^j:\H_m^j\to\L_m^{II}$, $j\in\{1,2,3\}$
such that
$$\pi_m^K\circ\Psi_m^j=\psi\circ\sigma_m^j,$$
where $\pi_m^K:\L_m^K\to\C_J$ and $\sigma_m^j:\H_m^j\to\C_a$
are the forgetful maps, and $\psi$ is the
function (\ref{Ja}).
We will show that these maps $\Psi_m^j$, are
orbifold coverings with respect to the appropriate orbifold structures
defined on some compactifications of $\L_m^K$ and $\H_m^j$.

This will permit us to compute the genera of
components of $\H_m$ via the
Riemann--Hurwitz formula. Once the genera and degrees are known
one can conclude that these curves are non-singular by the 
``genus--degree formula'':
\begin{equation}\label{gd}
g\leq(d-1)(d-2)/2,
\end{equation}
where we have equality only for non-singular curves.
We consider compactifications $\overline{\H}_m^j$ obtained
from $\H_m^j$ by filling in the punctures. Equivalently,
$\overline{\H}_m^j$ is the normalization of the algebraic curve
obtained as the projective closure of the zero set of equation
$H_m^j=0$.
\begin{theorem}\label{theorem2}
The maps $\overline{\H}_m^j\to\CP^2$ for $j\in\{0,\ldots,3\}$ 
are non-singular embeddings, in particular $\overline{\H}_m^j$
are irreducible. The
degrees of the ramified coverings $\Psi_m^j$
are $6$ for $j=0$ and $2$ for $j\in\{1,2,3\}$. 
\end{theorem}
\def\Q{\mathbf{Q}}

So we have interesting sequences of non-singular planar curves
$\overline{\H}^j_m$
defined over $\Q$ for which 
degrees and genera have been explicitly determined. Only a few
such examples of high genus are known to the authors,
see~\cite{Bries}, \cite{MO}.

One can deduce Theorem~\ref{theorem1} from Theorem~\ref{theorem2}:
if we know that $\overline{\H}_m^j$ are non-singular,
we can find their genera from the equality in (\ref{gd}) and obtain all
topological characteristics of $\L_m$ using the orbifold coverings
$\Psi_m^j$.
\begin{cor}\label{corollary1}
All singularities of
irreducible components of
surfaces (\ref{F}) are contained in the
lines $(0,t,0)$ and $(0,0,t)$.
\end{cor}
Computation for $m\leq 6$ shows that the only singularities
of surfaces (\ref{F}) are the zeros of the discriminant, but we do
not prove this in this paper.

These results allow us to prove a conjecture of Maier about degrees
of Cohn's polynomials \cite[Conj. 3.1(ii)]{Maier}. We recall the definition.
 Let $F_m=F_m^IF_m^{II}$ be the irreducible
factorization. Let $D_m^K$ be the discriminant of $F_m^K$ with respect to
$\lambda$. Then $D_m^K(g_2,g_3)$ is quasi-homogeneous, that is 
the curves $D_m^K(g_2,g_3)=0$ are
invariant under the scaling transformations
(\ref{scaling}). Therefore, the equations $D_m^K(g_2,g_3)=0$ can be
rewritten as $C_m^K(J)=0$, and these $C_m^K$ are called
{\em Cohn's polynomials}.
\begin{cor}\label{corollary2} (Maier's conjecture)
 $\deg C_m^I=\lfloor ((d_m^I)^2-d_m^I+4)/6\rfloor$ and
$\deg C_m^{II}=d_m^{II}(d_m^{II}-1)/2$, where $d_m^K=\deg_\lambda F_m^K,$
as in (\ref{degFI}), (\ref{degFII}).
\end{cor}

Our second application, and the original motivation of this work
is 
the problem of describing degeneration of
metrics of constant positive curvature with conic
singularities which recently attracted substantial attention,
\cite{L1,L2,L4,EMP,L3,LW,MP1,MP2}. Let $S$ be a compact surface,
and $\alpha_1,\ldots,\alpha_n$ positive
numbers. Consider Riemannian metrics on $S$
with $n$ conic singularities with the angles $2\pi\alpha_j$.
Each such metric defines a conformal structure on S with $n$ marked points,
so we have the forgetful map assigning this conformal structure to the metric.
The goal is understanding the space of such metrics and the properties
of the forgetful map.

In this paper
we restrict ourselves to the case when $S$ is a torus with one singularity
with the angle $2\pi\alpha$, where $\alpha=2m+1$ is an odd integer.
Following \cite{MP2} we denote the set of all such metrics by
$\Sph_{1,1}(\alpha)$. One can embed $\Sph_{1,1}(\alpha)$ into $\Lame_m$:
the image of this embedding consists of those Lam\'e equations
whose projective monodromy is {\em unitarizable}
(that is conjugate to a subgroup of $PSU(2).$)

We have the forgetful map
$\Sph_{1,1}(\alpha)\to \C_J$
which assigns to each metric its conformal class. It is known that
when $\alpha>1$ is not an odd integer, this
forgetful map is proper and surjective \cite{BT}, \cite{MP2},
and $\Sph_{1,1}(\alpha)$
is properly embedded in $\Lame_m$.

This is not the case for odd integers $\alpha$, the fact discovered
in \cite{LW} (see also \cite{BE} for a shorter proof of the
main result of \cite{LW}).
As the conformal class varies, a spherical  metric can degenerate.

Let us define the set $\LW_m\in\Lame_m$ consisting
of all Lam\'e equations whose projective monodromy consists
of collinear translations (by the periods of the integral (\ref{Ab})
below). In Section~\ref{mono} we show that
\begin{theorem}\label{LWm}
\begin{equation}\label{lw}
\partial\Sph_{1,1}(2m+1)=\LW_m,
\end{equation}
where the boundary is with respect to $\Lame_m$.
\end{theorem}

Then we describe the set $\LW_m$. 
\begin{theorem}\label{theorem4} The set $\LW_m$ consists of $m(m+1)/2$
curves. These curves and their projections on $\C_J$ are smooth
real analytic curves
(images of intervals under analytic functions with non-vanishing derivatives).
\end{theorem}
We propose to call projections of the curves $\LW_m$ to $\C_J$
Lin--Wang curves.
They can be seen in 
the pictures in \cite{LW,BE} and our Figs.~12, 13 for $m=1$; 
in \cite{L3,L4} and our Figs, 14, 15 for $m=2$, and our Figs.~16--18 for $m=3$.
Large part of the papers \cite{L1}, \cite{L2} and \cite{L3} is
dedicated to analytic study of these curves for small $m$;
here we propose a different,
geometric description of them. 

\subsection{Description of the method}\label{description}

Our main tool in this paper is
a new geometric interpretation of Lam\'e functions and their moduli
space $\L_m$.
Lam\'e functions correspond to what we call 
{\em translation structures} on tori with one conic singularity 
with the angle $2\pi(2m+1)$ and $m$ simple poles.
(``Simple pole'' refers to the developing map; its differential
has double poles).

Let $(S,O)$ be an elliptic curve, that is $S$ is a torus with a marked
point
$O\in S$, and
$m\geq 0$ an integer. Translation structures we are talking about
can be identified with Abelian differentials of
the second kind $g(z)dz$ on $S$
with
single zero of multiplicity $2m$ at $O$, and $m$ double poles, subject
to the condition that {\em all residues vanish}. Two translation structures
are {\em equivalent} if the differentials differ by a non-zero
constant factor.

We refer to a survey \cite{Z} of translation structures. Structures
considered in this survey have no poles and correspond to Abelian
differentials of the first kind.
To explain the name ``translation structure'',
consider the Abelian integral
\begin{equation}\label{Ab}
f(z)=\int_{z_0}^z g(\zeta)d\zeta.
\end{equation}
This is a multi-valued function on $S$ with a single critical
point at $O$, and the monodromy of $f$ consists of
translations by the periods of the integral coming from
the fundamental group of $S$. This function $f$
is a developing map of a singular
flat structure on $S$: it has one conic singularity
with the angle $2\pi(2m+1)$ at $O$ and $m$ simple poles; the monodromy of
this structure consists of translations, and the local monodromy at
all points is trivial.
\begin{prop}\label{trans-structures} The correspondence
$w\mapsto\Omega= dx/(uw^2)$ is a bijection
between the space of Lam\'e functions and the
space of triples $(S,O,\Omega)$,
where $S$ is a torus, $O\in S$ a point, and $\Omega$ is
an Abelian differential on $S$, which switches sign under
the conformal involution, and has a single zero of
multiplicity $2m$ at the point $O$,
and $m$ double poles with vanishing residues.
This bijection defines a biholomorphic map between $\L_m$ and
the moduli space of Abelian differentials of considered type, up to scaling.
\end{prop}
\vspace{.1in}

One can pull back the {\em flat} metric on $\C$ via $f$ and obtain a flat
metric on the torus with one conic singularity at $O$ with the angle $2\pi(2m+1)$
and $m$ simple ``poles''. A {\em pole} of a flat metric is a point which
has a punctured neighborhood
isometric to $\{ z\in\C:|z|>R\}$ with Euclidean metric,
for some $R>0$.
We call our torus equipped with this metric
a {\em flat singular torus}. Two flat singular tori are
considered {\em equivalent} if there is an orientation-preserving
diffeomorphism between them
multiplying the metric by a non-zero constant.

To study flat singular tori
we cut each of them into two congruent flat singular triangles.
{\em Congruent} means ``related by an {\em orientation-preserving} isometry''. 
\begin{definition}
{\em A} flat singular triangle (FT) {\em is a closed disk
with three marked boundary points which are called} corners,
{\em equipped with a flat
metric with conic singularities at the corners and possibly
simple poles inside the disk or on the open boundary arcs (}sides{\em)
between adjacent corners,
and such that the sides 
are geodesic. A side passing through a pole
must be ``unbroken'' at this pole: in the chart $\{ z\in\C:|z|>R\}$
it corresponds to two rays of the} the same {\em line.}
\end{definition}

Alternatively an FT can be described as a triple
$(D,\{ a_j\},f)$, where $D$ is a closed disk in the complex plane,
$(a_1,a_2,a_3)$
three distinct boundary points, and $f$ a locally univalent meromorphic
function on $D\backslash\{ a_j\}$ with conic singularities at
$a_j$, which means
$$f(z)=f(a_j)+(z-a_j)^{\alpha_j}h_j(z),$$
where $\alpha_j>0$, and $h_j$ is analytic near $a_j$, $h_j(a_j)\neq 0$,
and such that the images of the sides
$f([a_j,a_{j+1}])\subset \ell_j\cup\{\infty\}$ where $\ell_j$
are three straight lines (not necessarily distinct)
in the complex plane.

Two such triples $(D,\{ a_j\},f)$ and $(D',\{ a_j^\prime\}, g)$
are {\em equivalent} if there is
a conformal\footnote{We assume that the complex plane has
the standard orientation, and conformal maps preserve it
otherwise we call them anti-conformal.}
homeomorphism $\phi:D\to D'$, $\phi(a_j)=a_j^\prime$,
and complex constants $c_1,c_2$ such that
$f=c_1g\circ\phi+c_2$.

These two definitions of FT are equivalent: for a given triple, we pull back
the Euclidean metric from $\C$ via $f$ and obtain a metric on $D$
satisfying the first definition. In the opposite direction, given such
a metric we obtain $f$ as its developing map.

It is easy to see that the sum of the interior angles $\pi\alpha_j$ of a flat
singular
triangle is an odd integer multiple of $\pi$. If one
angle is an integer multiple of $\pi$,
then all three of them
are 
integer multiples of $\pi$.

An FT can be visualized by making a picture of its image in the plane under 
the developing map $f$. Such a picture consists of three lines
(not necessarily distinct), three (pairwise distinct) points
of their intersections
$f(a_j)$, 
and a marking of the angles. See Figs. 1a, 1b, 1d, 1e, 2a-f.
We mark the angles by little arcs near the images of the corners,
and write $a_j$ instead of $f(a_j)$ in the pictures.
Two such pictures define
the same FT if they are related by a complex affine map.

An FT is called {\em balanced} if its
angles $\pi\alpha_1,\pi\alpha_2, \pi\alpha_3$ satisfy the
triangle inequalities
\begin{equation}\label{balanced}
\alpha_i\leq\alpha_j+\alpha_k,
\end{equation}
for all permutations of $(i,j,k)$.
For example, a triangle with angle sum $\pi$
in Fig.~1a is balanced if and only if all its angles are $\leq\pi/2$.
A triangle with angle sum $3\pi$
in Fig.~1b is balanced if and only if the largest angle
is $\leq 3\pi/2$.

All flat singular triangles
whose angles are integer multiples of $\pi$ are balanced.
A balanced triangle is called
{\em marginal} if we have equality in (\ref{balanced})
for at least one permutation $(i,j,k)$.
Otherwise it is called {\em strictly balanced}.
\vspace{.1in}

{\em We abbreviate the expression ``balanced flat singular triangle''
as BFT.}
\vspace{.1in}

Our main technical result is the following
\vspace{.1in}

\noindent
\begin{theorem}\label{theorem3}
Every flat singular torus has a 
decomposition into two congruent BFT.
When the triangles are strictly balanced,
this decomposition is unique up to a cyclic permutation of
the corner labels. 
If they are marginal, there are at most two such decompositions:
a marginal triangle and its reflection define the same torus.
\end{theorem}

This is similar to Theorem 1.3 in \cite{EMP} for spherical tori
with one singularity. Theorem \ref{theorem3} gives a parametrization
of our moduli space $\L_m$ by a simpler moduli space $\T_m$ for
BFT's with the sum of the angles $\pi(2m+1)$. This last space admits a nice
partition into open cells which is used to prove Theorem~\ref{theorem1}.

The plan of the paper is the following. In section \ref{slam} we recall basic
facts about Lam\'e equations and Lam\'e functions, and explain
the connection between Lam\'e functions and translation structures.

In section \ref{sBFT} we discuss BFT
and define
a map $\Phi=\Phi_m:\T_m\to \L_m$.
Explicit local coordinates on $\T_m$
are described in Section~\ref{complexanalytic}, and we show that
$\Phi$ is complex analytic and proper. The proof of properness
is based on the study of geodesics on flat singular tori.

In section \ref{sth3} we prove the first part of Theorem \ref{theorem3},
surjectivity of $\Phi$.

In section \ref{spaces} we develop
a classification of BFT and explicitly describe a partition
of $\T_m$ into open cells.
We show that $\T_m$ consists of two connected components,
and that $\Phi$ is in fact $3$-to-$1$ on
the subset of strictly balanced triangles, and $6$-to-$1$ on
the subset of marginal triangles. Factoring $\T_m$ by an appropriate equivalence
relation we obtain a space $\T^*_m$ and show that the induced map
$\Phi^*_m:\T^*_m\to\L_m$ is injective, thus completing the proof
of Theorem~\ref{theorem3}.

In section \ref{ident} we analyze the natural
partition of $\T_m$
and prove Theorem~\ref{theorem1}.

In section \ref{nons} we prove Theorem \ref{theorem2}
and its two corollaries. This is based
on Lemma~\ref{lemma2} which is proved in
Section~\ref{tarasov}.

The last two sections are devoted to spherical metrics:
in Section \ref{mono} we discuss the monodromy of Lam\'e equations, and
in Section \ref{linwang} we 
produce equations and pictures of Lin--Wang curves, enumerate them,
and show that they are smooth and real analytic.
\vspace{.1in}

{\em Since we refer to some figures many times, in different places,
all figures are collected at the end of the paper, after the
reference list.}
\vspace{.1in}

We thank Walter Bergweiler, Robert Maier and Vitaly
Tarasov for useful discussions. Walter Bergweiler produced Figs.~12--18
and a Maple program generating polynomials $F_m$.
Vitaly Tarasov proved Lemma~\ref{lemma2}. We also thank
Eduardo Chavez Heredia  for bringing \cite{ST} to our attention.

\section{Lam\'e equations, Lam\'e functions and
\newline
translation structures}\label{slam}

We use the form (\ref{lame2}) of the Lam\'e equation.
Every solution $W$ of (\ref{lame2}) is meromorphic in the $z$-plane.
Indeed, by the existence theorem for linear ODE, the singularities of
solutions belong to the lattice $\Lambda$,
and plugging a power series for $W(z)$ at $0$ shows that
there are two linearly independent meromorphic solutions.
\vspace{.1in}

{\em Proof of Proposition~\ref{trans-structures}.}
\vspace{.1in}

We start with a Lam\'e function and assign to it a translation structure.
Let $W$ be a meromorphic
solution of (\ref{lame2}) whose square is even and $\Lambda$-periodic. Then 
\begin{equation}\label{period}
W(z+\omega)=\pm W(z),\quad \omega\in\Lambda.
\end{equation}
Now 
\begin{equation}\label{w2}
W_1(z):=W(z)\int^z W^{-2}(\zeta)d\zeta
\end{equation}
is another solution of (\ref{lame1}), linearly independent of $W$,
which can be seen by direct computation, so it must be also
meromorphic, since every solution of (\ref{lame2}) is meromorphic.
It follows that all residues of $W^{-2}(\zeta)d\zeta$ vanish.
Now the ratio of two solutions $f=W_1/W$ is an Abelian integral
of the second kind,
and as a ratio of two solutions of a second order linear equation,
it also satisfies the Schwarz equation
\begin{equation}\label{schwa}
\frac{f'''}{f'}-\frac{3}{2}\left(\frac{f''}{f'}\right)^2=
-2\left( m(m+1)\wp+\lambda\right),
\end{equation}
thus all critical points of $f$ are in $\Lambda$
and $f$ is $(2m+1)$-to-$1$ at these points.

So every Lam\'e function defines a translation structure of the desired type.
The differential
$W^{-2}(z)dz$ descends via the map 
$$z\mapsto (\wp(z),\wp^\prime(z))=(x,u),$$
to the differential $dx/(w^2u)$ on the curve (\ref{W})
as stated in
Proposition~\ref{trans-structures}. 

Conversely, suppose that a translation structure of the described type is
given and its developing map is
\begin{equation}\label{abel}
f(z)=\int^z g(\zeta)d\zeta
\end{equation}
is given. Then $g$ is an even elliptic function with $m$ double poles
per period parallelogram, vanishing residues,
and zeros of order $2m$ at the points of $\Lambda$.
So $g=W^{-2}$ for some meromorphic function $W$ satisfying 
(\ref{period}). Now we define $W_1$ by (\ref{w2}),
and a direct computation shows that
$$WW_1^\prime-W'W_1=1.$$
This means that $W$ and $W_1$ are two linearly independent solutions
of some equation $W''=PW$, where $P$ is an elliptic function
with periods $\Lambda$. As the only pole
of $P$ can occur at a critical point of $f$, we conclude that poles 
of $P$ must be at the points of $\Lambda$, and a simple calculation
with power series at $0$ shows that
$$P(z)=m(m+1)z^{-2}+O(1),\quad z\to 0,$$
so $P(z)=m(m+1)\wp(z)+\lambda$ with some $\lambda\in\C$,
and thus $W$ is a Lam\'e function.

We recall that two Lam\'e functions $W_1$ and $W_2$ are  equivalent if
$W_1(z)=cW_2(kz)$ for some $c$ and $k$ in $\C^*$.
Translation structures are equivalent
if their developing maps $f_1,f_2$ are related by post-composition
with an affine map: $f_1(z)=af_2(kz)+b$, $a,k\in\C^*,\; b\in\C$.

We proved that {\em equivalence classes of degree $m$ Lam\'e functions
are in one-to-one correspondence with classes of translation
structures on the torus with one conic
singularity with the angle $2\pi(2m+1)$.}

It is clear that this correspondence is continuous and holomorphic,
therefore it is biholomorphic.\hfill$\Box$
\vspace{.1in}

Let us recall how the polynomial $F_m$ 
is computed (another method is described in Section 6).
It is convenient to use the algebraic form of the Lam\'e equation
(\ref{lame1}) which can be also written as
\begin{equation}\label{newlame}
w''+\frac{1}{2}\left(\sum_{j=1}^3\frac{1}{\zeta-e_j}\right)w'=
\frac{m(m+1)\zeta+\lambda}{4(\zeta-e_1)(\zeta-e_2)(\zeta-e_3)}w,
\end{equation}
where $w=y\circ\wp$. 
Since the only singularities in $\C$ of
this equation are $e_j$, and the local exponents
at these singularities are $\{0,1/2\}$, each Lam\'e function can be written
as
\begin{equation}\label{lfunction}
w(\zeta)=c\prod_{j=1}^3(\zeta-e_j)^{k_j/2}\prod_{j=1}^n(\zeta-\zeta_j),
\end{equation}
where $k_j\in\{0,1\}$ and
\begin{equation}\label{num}
\sum_{j=1}^3 k_j+2n=m.
\end{equation}
Plugging $\zeta=\zeta_k$ into (\ref{newlame}) we obtain after simple
calculations (see, for example, \cite[23.21]{WW})
the following system of equations for $\zeta_j$
\begin{equation}\label{bete}
2\sum_{j:j\neq k}^n\frac{1}{\zeta_k-\zeta_j}+
\sum_{j=1}^3\frac{k_j+1/2}{\zeta_k-e_j}=0,\quad 1\leq k\leq n.
\end{equation}
This system of equations determines Lam\'e functions. According to a
theorem of Heine and Stieltjes \cite[23.46]{WW}, system (\ref{bete}) has
at most $n+1$ solutions and exactly $n+1$ for generic $e_j$. Moreover,
it has exactly $n+1$ solutions when all $e_j$ are real.
System (\ref{bete}) is a very special case of the Bethe ansatz equations
which frequently occur in mathematical physics and in the study
of metrics with conic singularities. 

Lam\'e functions are classified according to the values of $k_j$ in (\ref{lfunction}):
traditionally the number
$$1+\sum_{j=1}^3k_j\in\{1,2,3,4\}$$
is called the {\em kind} of a Lam\'e function \cite{WW}.

Using (\ref{num}) and Stieltjes theorem
we conclude that the total number of Lam\'e functions
for a given generic lattice is $2m+1$, and this is the degree
of the Lam\'e spectral polynomial $F_m$ with respect to $\lambda$.
The number of Lam\'e functions is exactly $2m+1=d_m^I+d_m^{II}$
for lattices with
real $e_j$, or equivalently with real $g_j$ with $g_2^3-27g_3^2>0$.

It is easy to see that for even $m\geq 2$ there exist Lam\'e functions
of the first and third kind; we denote the corresponding subsets by $\L_m^I$
and $\L_m^{II}$. Similarly, when $m\geq 3$ is odd, we define $\L_m^I$ as the
set of Lam\'e functions of the fourth kind, and $\L_m^{II}$
as the set of Lam\'e functions of the second kind. 
This explains why
$\L_m$ has {\em at least two components} when $m\geq 2$.
The more difficult result, which is a part of Theorem~\ref{theorem1},
is that these components are in fact irreducible. 
\vspace{.1in}

{\em Connection with spin structures.}
\vspace{.1in}

Connected components of moduli spaces of Riemann surfaces endowed
with nonzero holomorphic differentials were classified by
Kontsevich and Zorich \cite{KZ}; the case
of meromorphic differentials was treated by Boissy \cite{Bo}.
A consequence of Theorem 4.1 in \cite{Bo} is that the moduli space of triples
$(S,O,\Omega)$, where $(S,O)$ is an elliptic curve
and $\Omega$ is
a meromorphic differential on $S$ with a zero at $O$ of order $2m$
and $m$ poles $q_1,\dots,q_m$ of order $2$ (and arbitrary residues),
has exactly two connected components.
Moreover, such components are distinguished by the
{\it{spin invariant}} (already defined in \cite[Sect. 2.2]{KZ}).
In our particular case, the spin invariant of $(S,O,\Omega)$ is {\it{odd}}
if there exists a function on $S$ with simple
poles at $q_1,\dots,q_m$ and a zero of order $m$
at $O$, and it is {\it{even}} if such function does not exist.

\begin{prop}
A Lam\'e function $w$ is in $\L_m^I$ if the spin
invariant of the corresponding translation surface is odd,
and $w$ is in $\L_m^{II}$ if the spin invariant is even.
\end{prop}

{\em Proof.}
Let $(S,O,\Omega)$ be the translation surface associated
to the Lam\'e function $w$.
In particular, $\Omega=\varphi/w^2$,
where $\varphi$ is a nonzero holomorphic
differential on $S$, see Proposition~\ref{trans-structures}
Then
the spin invariant of $(S,O,\Omega)$ is odd
if and only if $w$ is a well-defined function on $S$,
which happens if and only
if $w=Q(x)$ or $w=Q(x)u$ for a suitable $Q\in\mathbb{C}[x]$.
Hence, the spin invariant
is odd if and only if $w$ 
is of type $I$.\hfill$\Box$

As a consequence of  \cite[Thm. 4.1]{Bo} we obtain 

\begin{cor}
The space $\mathbf{L}_m$ of Lam\'e functions
is the disjoint union of the subset
$\mathbf{L}_m^I$
of Lam\'e functions of type $I$ and of the subset
$\mathbf{L}_m^{II}$
of Lam\'e functions of type $II$.
\end{cor}
The techniques of \cite{Bo}
do not allow to study connected components
of moduli of meromorphic differentials with vanishing residues.
So it does not follow from \cite{Bo}
that $\mathbf{L}_m$ has exactly two connected components, and so that
$\mathbf{L}_m^I$ and $\mathbf{L}_m^{II}$ are connected.

\section{Balanced flat singular triangles}\label{sBFT}

In this section we consider balanced flat singular triangles (BFT)
with marked corners $a_1,a_2,a_3$,
enumerated according to the positive orientation of the boundary
(so that the region is on the left when we trace the boundary).

We denote the interior angles at these corners by $\pi\alpha_1,\pi\alpha_2,
\pi\alpha_3$. As we already noticed, the sum of the angles is
an odd multiple of $\pi$, more precisely
\begin{equation}\label{suma}
\alpha_1+\alpha_2+\alpha_3=4n+2k+1=2m+1,
\end{equation}
where $n$ is the number of interior poles, and $k$ is the number 
of boundary poles. 
To prove (\ref{suma}) we recall the argument used in the
proof of the Schwarz--Christoffel formula. Consider the developing map $f$ defined
in the upper half-plane with corners at $(a_1,a_2,a_3)=(0,1,\infty)$.
Since the monodromy of $f$ is affine, $f''/f'$ must
extend to the complex plane as
a rational function whose poles are symmetric with respect to
the real line, and
we have
$$\frac{f''}{f'}(z)\sim\frac{\alpha_j-1}{z-a_j},\quad z\to a_j,
\quad j\in\{1,2\},\quad
\frac{f''}{f'}(z)\sim-\frac{\alpha_3+1}{z},\quad z\to\infty,$$
and 
$$\frac{f''}{f'}(z)\sim -\frac{2}{z-z_j},\quad z\to z_j$$
at the poles. So (\ref{suma}) follows by the Residue Theorem.

Fig.~1a shows a triangle with $n=k=0$. In Fig.~1b $n=0, k=1$.
Fig.~1d shows a triangle with $n=1, k=0$ (left) and two triangles with
$n=0, k=2$ (right).
Fig.~1e shows three triangles with $n=1, k=0$.  Fig.~2
shows all
types of triangles with sum of the angles $5\pi$ ($m=2$).

We also recall that {\em either none of the $\alpha_j$ or all of them
are integers.}
\begin{prop}\label{prop1}
An FT with non-integer $\alpha_j$ is completely
determined by the angles, and any positive angles $\pi\alpha_j$ where
$\alpha_j$ are not integers and their sum is odd can occur.

For integer angles, the necessary and sufficient condition of
existence of FT is that the sum of $\alpha_j$ is odd
and (\ref{balanced}) is satisfied. For any such angles,
there are three 
one-parametric families of FT's.
\end{prop}
A similar result for spherical triangles was proved in \cite{E}.
\vspace{.1in}

{\em Proof.} The first statement is essentially due to Klein \cite{Klein}.
If our triangle is modeled
on the upper half-plane $D=\overline{H}$ with vertices $(a_1,a_2,a_3)=(0,1,\infty)$
then the developing map $f:\overline{H}\to\bC$ satisfies
the Schwarz equation with three singularities $(0,1,\infty)$
and $\alpha_j$
determine this equation completely.
When the $\alpha_j$ are not integers,
the monodromy representation
corresponding to this equation is non-trivial, and there is only
one choice, up to an affine transformation, of a solution
with affine monodromy.

If all $\alpha_j$ are integers, then the monodromy representation
is trivial, and the developing map $f$ extends
from its domain $D$ to the whole Riemann sphere and we obtain
a rational
function with three critical points $a_1$, $a_2$ and $a_3$.
The images of all three sides under the developing map $f$ belong the same
line $\ell$. Preimage $f^{-1}(\ell)\cap \overline{D}$ is called the {\em net}
of $f$ (see \cite{EG}). The net is a cell decomposition of $D$
with three vertices $a_1,a_2,a_3$. The $1$-cells of this decomposition
are disjoint chords of the disk $D$ and three arcs of $\partial D$.
The number of chords in
the interior of $D$ which are adjacent to $a_j$ is $\alpha_j-1\geq 0$.
If $m_j$ is the number of chords from $a_i$ to $a_k$, then
$\alpha_j-1=m_i+m_k$, and these three equations have unique non-negative
solution $(m_1,m_2,m_3)$ if and only if the integers $\alpha_j\geq 1$ satisfy
(\ref{balanced}) and their sum is odd.
Thus the angles determine the net. Once the net
is given, the developing map can be recovered from the values $f(a_1),f(a_2),
f(a_3)$ which in the considered case belong to a line.
It is clear that one of the $2$-faces of a net is a triangle while
the others are digons.
Suppose that the triangular face is mapped by
$f$ onto the upper half-plane. Then we have three possible
orderings of $f(a_j)$ on the real line. By scaling we may arrange
that $\min_j f(a_j)=0$ and $\max_j f(a_j)=1$, then the intermediate
point of the three $f(a_j)$ serves as a parameter. So
the set of triangles with given integer angles is parametrized by three
intervals.

\hfill$\Box$
\vspace{.1in}

Examples of nets of triangles are shown in Fig.~1c, where the stars
mark the location of the poles of $f$. These three nets correspond to
triangles with the angles $(2\pi,3\pi,4\pi)$.
In Fig.~1d, three
triangles with the angles $(\pi,2\pi,2\pi)$ are shown together with their nets.
\vspace{.1in}

Let us define two types of triangles which we call {\em primitive}.
A primitive triangle of type $A$ has all angles in $(0,\pi)$ and their
sum is $\pi$ (Fig.~1a),
and the primitive triangle of type $B$ has one angle in $[\pi,2\pi)$,
two others in $(0,\pi]$, and the sum of the angles is $3\pi$ (Fig.~1b).

\begin{prop}\label{prop2}
Every BFT can be obtained from one of the two primitive
triangles $A$ or $B$ by gluing half-planes
to the sides.
\end{prop}

{\em Proof.}
This is essentially due to Klein \cite{Klein}, but we sketch a proof. Consider the
case of non-integer angles.
Let $\alpha_i=\{\alpha_i\}+\[\alpha_i\],$ be the decomposition into
fractional and integer parts. 
\vspace{.1in}

\noindent
a) If $\sum\[\alpha_i\]$ is even, then $\sum\{\alpha_i\}$ is odd.
Since the last
sum belongs to $(0,3)$, we must have
\begin{equation}\label{sumone}
\sum\{\alpha_i\}=1.
\end{equation}
Since our triangle is balanced, and we have (\ref{sumone}), we obtain
for all permutations $(i,j,k)$ of $(1,2,3)$:
\begin{eqnarray*}
\[\alpha_i\]&=&\alpha_i-\{\alpha_i\}\leq\alpha_j+\alpha_k-\{\alpha_i\}=
\[\alpha_j\]+\[\alpha_k\]+\{\alpha_j\}+\{\alpha_k\}-\{\alpha_i\}\\
&<&
\[\alpha_j\]+\[\alpha_k\]+1.
\end{eqnarray*}
So $\[\alpha_i\]\leq \[\alpha_j\]+\[\alpha_k\]$. It follows that
the following
quantities are non-negative
integers:
$$x_i=(\[\alpha_j\]+\[\alpha_k\]-\[\alpha_i\])/2,$$
and we have $\[\alpha_i\]=x_j+x_k,$ for all permutations $(i,j,k)$
of $(1,2,3)$.
So we can take a triangle with the angles $\{\alpha_j\}$ of the type A
and glue $x_j$ half-planes to the side opposite to $a_j$, for each $j$.
The resulting triangle has the
same angles as our original triangle, so it is the same triangle by
Proposition~\ref{prop1}.
\vspace{.1in}

\noindent
b) If $\sum\[\alpha_i\]$ is odd, then we decrease one of the
integer parts (for example, the largest one) by $1$, and increase
the corresponding fractional part by $1$. So we set for some $i$
$$p_i=\[\alpha_i\]-1,\quad\alpha_i^{\prime}=\{\alpha_i\}+1,$$
leaving other numbers unchanged ($\alpha_j^\prime=\{\alpha_j\},\;
p_j=\[\alpha_j\]$, and similarly for $k$).
Now we have
\begin{equation}\label{b1}
\alpha_i^{\prime}\in(1,2),\quad\alpha_j^{\prime}\in(0,1),\quad
\alpha_k^{\prime}\in(0,1),
\end{equation}
and since $\sum\{\alpha_\ell\}$ is even and less than $3$ in this case,
it must be $2$, and
thus
\begin{equation}\label{sumthree}
\sum\alpha_\ell^{\prime}=3.
\end{equation}
Now, since our triangle is balanced, and using (\ref{sumthree}) and (\ref{b1})
we obtain
$$p_i=\alpha_i-\alpha_i^{\prime}\leq
\alpha_j+\alpha_k-\alpha_i^{\prime}=p_j+p_k
+\alpha_j^{\prime}+\alpha_k^{\prime}-\alpha_i^{\prime}
\leq p_j+p_k+1,$$
and since the sum of $p_\ell$ is even, we conclude that
$p_i\leq p_j+p_k.$
So we can define $x_i,x_j,x_k$ as in part a) and conclude that our triangle
is obtained from a triangle with the angles $\alpha_\ell^\prime$ of
the type B, by gluing $x_\ell$ half-planes to the side opposite to $a_\ell$.

Triangles with
integer angles can be considered using their nets introduced
in the proof of Proposition~\ref{prop1}.
A net is a chord diagram in a disk with three vertices
on the circle. Evidently each net has one triangular face,
and the rest of the faces are digons. Triangular face corresponds to
the primitive triangle and digons are half-planes.
Since the angle sum is
odd, only case b) can occur, and the primitive triangle
is of type B with the angles $(\pi,\pi,\pi)$.
Thus the primitive triangle in this case
is just a half-plane with three marked boundary points.\hfill$\Box$
\begin{cor}\label{corpol}
Each side of a BFT contains at most one pole, and
the developing map sends each side $(a_i,a_{i+1})$ injectively
either to the interval
\newline
$(f(a_i),f(a_{i+1}))$, or to the complement
of this interval on the line in $\C\cup\{\infty\}$ containing it.\hfill$\Box$
\end{cor}

\noindent
{\bf Remarks.} 1. By analyzing the proof of Proposition~\ref{prop2} one
can obtain the following criterion: the  side opposite to $a_i$
contains a pole if and only if
$$\lceil(\[\alpha_j\]+\[\alpha_k\]-\[\alpha_i\]-1)/2\rceil\quad\mbox{is odd}.$$

\noindent
2. It follows from this Corollary and from (\ref{suma})
that for each $m\geq 2$ we have exactly two possibilities
for the number of sides with a pole. When $m$ is even,
we have either none or 2 sides with poles. If $m$ is odd, we have
ether one or three sides with poles. This shows that 
the set $\T_m$
of all BFT's with the angle sum $\pi(2m+1),\; m\geq 2$ must have at least
four connected components. We show in Section \ref{spaces} that there
are exactly four.
\vspace{.1in}

\noindent
3. The decomposition into a primitive triangle and half-planes
stated in Proposition~\ref{prop2} is canonical when $\sum\[\alpha_j\]$
is even, but not canonical when it is odd.
In the latter case, we can
obtain from one to three such different decompositions, depending
on the number of positive $\[\alpha_j\]$.
\vspace{.1in}

\noindent
4. The primitive triangle $T'$ obtained from a balanced triangle $T$
may be unbalanced. In this case, there is always at least one half-plane
in $T$ glued to the side of $T'$ which is opposite to the largest
angle of $T'$. Indeed, let $\alpha_j^\prime$ be the angles of $T'$,
and $\alpha_1$ is the largest one. Then the angles of
$T$ are $\alpha_i=\alpha_i^\prime+x_j+x_k$, where $(i,j,k)$
is a permutation of $(1,2,3)$, and $x_j$ are the numbers of half-planes
glued to the sides opposite to $\alpha_j^\prime$, $j=1,2,3$. Since
$T$ is balanced, 
$$\alpha_1^\prime+x_2+x_3=
\alpha_1\leq\alpha_2+\alpha_3=\alpha_2^\prime+\alpha_3^\prime+x_2+x_3+2x_1,$$
so $\alpha_1^\prime\leq \alpha_2^\prime+\alpha_3^\prime+2x_1$. Thus if
$T'$
is unbalanced, we have $x_1>0$.
\vspace{.1in}

\noindent
{\em Construction of the map $\Phi$.}
\vspace{.1in}

Let $\T_m$ be the set of all balanced triangles with the sum of the angles
$\pi(2m+1)$.
For every $T\in\T_m$ we define a singular flat torus
$\Phi(T)$ in the following way. We take two copies of $T$
and identify each pair of corresponding sides
by the {\em orientation-reversing} isometry.

Thus all corners of both copies are glued into one point,
and the sides are glued into three simple loops on the torus based
at this point and otherwise disjoint. 

Notice that the resulting torus has an orientation-preserving
isometric involution which interchanges the two triangles.
The four fixed points of this involution are: the conic singularity
(corresponding
to the vertices of the triangle) and the ``midpoints of the sides'' which
are points of order $2$ on the elliptic curve. If
a side is unbounded, its midpoint is a pole.

With such gluing we obtain a flat singular torus
with one singularity 
with the angle $2\pi(2m+1)$. There are $m$ simple poles on the
torus coming from the poles of the metric on the triangle.
An interior pole of $T$ gives two poles on $\Phi(T)$, while
a pole on a side of $T$ gives one pole on $\Phi(T)$. 

Let $T$ be a BFT with corners $a_1,a_2,a_3$ and developing map $f$.
Let $b_j=(f(a_i)+f(a_k))/2$, where $(i,j,k)$ is a permutation of $(1,2,3)$.
We define affine maps $s_j$ to be rotations by $\pi$ about $b_j$.
Consider the group generated by $s_1,s_2,s_3$. It contains a subgroup
$G$ of index $2$ consisting of translations;
elements of $G$ are products of even numbers
of $s_j$. The following proposition is evident:
\begin{prop}\label{prop2.5}
The monodromy group of the developing map of the flat singular torus
$\Phi(T)$ is~$G$.\hfill$\Box$
\end{prop}
In particular, the monodromy consists of collinear translations
if and only if the angles of the triangle are integer multiples of $\pi$. 

Next we address the question when two different BFT's can give
the same (isometric up to a constant factor) tori.
This can happen in at least two ways:
\vspace{.06in}

\noindent
1. The triangles are obtained by a cyclic
permutation of the corners.
\vspace{.1in}

\noindent
2. Some pairs of marginal triangles define the same torus.
\vspace{.1in}

More precisely, we have
\begin{prop}\label{prop3}
For a marginal triangle $T$
with the angles $(\alpha_1,\alpha_2,\alpha_3)$ where
$\alpha_1=\alpha_2+\alpha_3$, and the triangle $T^*$
with the angles $(\alpha_1,\alpha_3,\alpha_2)$ the tori
$\Phi(T)$ and $\Phi(T^*)$ are congruent.
\end{prop}

$T$ and $T^*$ are related by a reflection, an orientation-reversing
isometry.
\vspace{.1in}

{\em Proof of Proposition~\ref{prop3}.} First we notice that a marginal triangle cannot
have integer angles (since the angle sum is odd),
so it is completely determined
by the angles. We use Klein's decomposition described in the proof of
Proposition~\ref{prop2}. Our triangle is obtained from a primitive triangle by gluing
half-planes to the sides.

We claim that for a marginal triangle,
no half-planes are glued to the side opposite to the larger angle.
Indeed, let $\alpha_j'$ be the angles of the primitive triangle,
and $\alpha_j=\alpha_j^\prime+p_j=\alpha_j^\prime+x_i+x_k$,
where $p_j,x_i,x_k$ are non-negative integers,
Then $\alpha_1=\alpha_2+\alpha_3$
implies
\begin{equation}\label{ma}
\alpha_1^\prime=\alpha_2^\prime+\alpha_3^\prime+2x_1,
\end{equation}
and we obtain that $x_1=0$ since
$\alpha_1<2,\; \alpha_2^\prime>0,\alpha_3^\prime>0.$
This proves the claim.

We also conclude from (\ref{ma}) with $x_1=0$ that the primitive triangle $T'$
with the angles $\alpha_1^\prime,\alpha_2^\prime,\alpha_3^\prime$ 
corresponding to our triangle $T$ is also marginal, with the larger angle
$\alpha_1.$ If $T'$ is of the type A, then $\alpha_1^\prime=1/2$.
If $T'$ is of the type B, then $\alpha_1^\prime=3\pi/2$.
Gluing together two copies of $T'$ along the side $(a_3,a_1)$ 
we obtain either a rectangle in the plane, or a complement of
a rectangle on the Riemann sphere. 
It is clear that rectangles obtained from
$T$ and $T^*$ are congruent, and
the corresponding tori $\Phi(T)$ and $\Phi(T^*)$ are also congruent, 
see Fig.~3.\hfill$\Box$
\vspace{.1in}

We give two examples illustrating Theorem~\ref{theorem3} in the simplest cases.
\vspace{.1in}

\noindent
{\bf Example 1.} Consider a triangle $T$ shown in Fig.~1a, whose angle sum is
$\pi$. It is balanced iff all angles are at most $\pi/2$. Gluing two congruent copies
of such a triangle along a common side,
we obtain a parallelogram in the plane. Identifying the opposite sides
of this parallelogram by translations we obtain a flat non-singular
torus $\Phi(T)$ ($m=0$).

Now consider a flat non-singular torus $L$ with a marked point $O$.
Let $f:\C\to L,\; f(0)=O$ be the universal cover. Then there is a
lattice $\Lambda\subset\C$ such that $L=\C/\Lambda$. A fundamental region
$D$ of this lattice can be taken in the form of a parallelogram which can
be normalized so that the shorter side is $[0,1]$. Let $[0,\tau]$ be the
longer side. It is well known that $\tau$ can be always chosen
in the fundamental region of the modular group
$$G=\{\tau:|\tau|\geq 1,|\Rea\tau|\leq 1/2\}.$$
To achieve this one can normalize so that the shortest non-zero element
of $\Lambda$ is $1$, then the shortest non-real element of $\Lambda$
is $\tau$.
For $\tau\in G$, each diagonal of $D$ 
breaks $D$ into a pair
of congruent triangles, such that at least one pair consists of balanced
triangles. Both pairs consists of balanced triangles if and only
if $D$ is a rectangle, in which case the triangles of different
pairs are marginal and are
reflections of each other. It is easy to see
that at least one diagonal breaks $D$ into balanced triangles
if and only if $\tau\in G$. This proves Theorem~\ref{theorem3} for
$m=0$. The proof in the general case is a generalization 
of this argument. \hfill$\Box$
\vspace{.1in}

\noindent
{\bf Example 2.} Consider a triangle $T$ in Fig.~1b.
The angle sum is $3\pi$ ($m=1$). Rotating $T$
by $\pi$ about the point $(a_1+a_3)/2$ we obtain a
congruent triangle $T'$. Gluing $T$ and $T'$ along the side which
contains pole, we obtain the ``exterior parallelogram'' $Q$.
Identifying the opposite sides of $Q$ by translations, we obtain the
flat singular torus $\Phi(T)$ with the angle $6\pi$
at the singularity. It will be proved later that every flat singular
torus with the angle $6\pi$ can be obtained from an exterior 
parallelogram by identifying the opposite sides by translations.
We claim that every exterior parallelogram $Q$ can be obtained
by gluing two {\em balanced} triangles of the type B. Indeed,
consider the diagonals of the parallelogram $\C\backslash Q$, and extend
them to $Q$. Each extension of a diagonal breaks $Q$ into two
triangles of the type B, and it is easy to see that only for
one diagonal these triangles are balanced, unless our $Q$
is an exterior rectangle as in Fig.~3. When $Q$ is an exterior rectangle,
our triangles are marginal, and they are related
by reflection as in Proposition~\ref{prop3}.
\hfill$\Box$
\vspace{.1in}

Now we define an {\em equivalence relation} $\sim$ on $\T_m$,
in accordance with statements 1,~2 before Proposition \ref{prop3}.
Two BFT's $T$ and $T'$
are equivalent if either there is a congruence between them
cyclically permuting the corners, or if they are marginal
and related by a reflection as described Proposition
\ref{prop3}. Then our map $\Phi$ is well defined on the equivalence
classes, and we prove in Section~\ref{sth3} that
the induced map 
$$\Phi_m^*:\T_m^*\,\longrightarrow\L_m,\quad \T_m^*=\T_m/\sim$$
is bijective. We denote by $\Phi^*$ the map defined
by this formula on the disjoint union $\T^*$ of $\T_m^*,\; m\geq 0$,
and mapping it to the disjoint union $\L$ of $\L_m$.

\section{Complex analytic coordinates on $\T_m$}
\label{complexanalytic}

We introduce a complex analytic structure on the set of
BFT. We define functions on the set of BFT:
$$\phi_{i,j,k}=\frac{f(a_i)-f(a_j)}{f(a_k)-f(a_j)}.$$
These functions are locally
injective and we use them as complex coordinates on $\T_m^*$.
The correspondence maps between charts are linear-fractional.
So they define
a complex analytic structure on $\T_m^*$ (and even a projective structure).

To show that the map $\Phi_m^*$ is complex analytic, we recall
Proposition~\ref{prop2.5} which implies that $\phi_{i,j,k}$
are ratios of periods of the differential in
Proposition~\ref{trans-structures}. It is clear
that periods and their ratios are analytic on $\L_m$ (which can
be identified with a space of proportionality classes of differentials).
We refer to a much more general statement of this kind in
\cite[Corollary 2.3]{Ba}.  

So $\phi_{i,j,k}$ are local complex analytic
coordinates on $\T_m^*$, and we have

\begin{prop}\label{prop4}
The map $\Phi_m^*:\T_m^*\to \L_m$ 
is complex analytic.\hfill$\Box$
\end{prop}
Next we prove 
\begin{prop}\label{prop5}
The maps $\Phi_m^*$ are proper. 
\end{prop}

{\em Proof.} We are going to show that a degenerating sequence
of triangles gives a degenerating sequence of tori. First we clarify
the notions of degenerating sequences.

For every flat singular torus $L$ with one conic singularity at $O$,
we define the set $G$ of all simple geodesics loops based at $O$.
Some  of these loops  may pass through a pole and thus
have infinite length.
The {\em systole} $s_1$ is the minimal length of all elements of $G$.
The {\em second systole} $s_2\geq s_1$ is the minimal length of
all elements of $G$ whose homotopy class is not a multiple
of the class of some element of $G$ of length $s_1$.
Since we identify tori with proportional metrics,
only the ratio $s_2/s_1$ is defined as a function on $\L_m$.
It is clear that $s_2/s_1$ is bounded on every compact
subset of $\L_m$.

{\em If $s_2/s_1\to\infty$ for a sequence of tori in $\L_m$, then
this sequence diverges in $\L_m$, and
we say that tori of this sequence degenerate.}

Proposition~\ref{prop2} defines for every balanced $T$ a primitive
triangle $T'$ which is of type A or B (Fig.~1a,b). 
When the sum of $\lfloor\alpha_j\rfloor$ is even, 
then $T'$ is of type $A$ and it is defined
uniquely. When this sum is odd, we may have up to three choices for $T'$;
they are of type $B$
(see Remark~3 after Corollary~\ref{corpol}).
We pick one of them, as we did
in the proof of Proposition~\ref{prop2}.

This primitive triangle $T'$ may be unbalanced. If this is the case,
then there is at least one half-plane in $T$ attached
to the side of $T'$
opposite to the largest angle of $T'$ (Remark 4 after Corollary~\ref{corpol}).
Then we denote by $T''$ the union (more precisely the result of gluing) of
$T'$ with this half-plane. When $T'$ is balanced we set $T''=T'$.
We call $T''$ the {\em balanced extension} of $T'$.
There are $4$ types of $T''$: balanced of type $A$ or $B$
and balanced extensions of unbalanced $A$ and $B$. We call
these last two types $A''$ and $B''$.

The numbers $|f(a_i)-f(a_j)|$ are all the same for
$T, T'$ and $T''$. Consider the ratios
$$\frac{f(a_i)-f(a_j)}{f(a_k)-f(a_j)}$$
It is easy to see that if a sequence of triangles in $\T_m$ leaves
every compact in $\T_m$ then one of these ratios tends to infinity.
We call triangles of such a sequence {\em degenerate},
and notice that {\em $T,T'$ and $T''$ degenerate simultaneously.}
(Unbalanced triangle of type $A$
can degenerate in a different way:
when a vertex tends to an interior point of the opposite side. But
its balanced extension does not degenerate in this case.)

We claim that $s_1,s_2$ for $\Phi(T'')$ are the same as for
$\Phi(T)$. Indeed, attaching $n$ half-planes to a side of
a triangle $T''$ results
in attaching a digon $D$ (with poles)
with angles $2\pi n$ at at its two corners to the torus $\Phi(T)$.
Every curve in such a digon with endpoints at $\partial D$ is
at least as long as the segment between its endpoints.
On the other hand, there are four
types of $T''$ of which $A$ and $B''$ have all sides bounded,
while the complementary segment of the unbounded side of $A''$ or $B$
is at least $s_2$ for $T''$.

It remains to consider the tori $\Phi(T'')$, and to show that
when $T''$ degenerates then $s_2/s_1\to\infty$ for these tori.
There are 4 cases to consider:
\vspace{.06in}

\noindent
a) $T'$ is balanced, then $T''=T'$. Then $T'$ is of type A or B,
and balanced,
and this is essentially the cases of Examples 1, 2 in the end of Section 3.
\vspace{.1in}

\noindent
b) $T'$ is unbalanced, $T''$ is the balanced extension
of $T'$ of the type $A$ or $B$.
One has to find the first and second systole of such tori $\Phi(T'')$.

We denote $A''$ and $B''$ the classes of triangles consisting
of  balanced extensions
of unbalanced triangles of type
$A$ and $B$. 
\vspace{.1in}

We are interested in the length spectrum of the set of geodesic loops
in $\Phi(T_1)$ where $T_1$ is a balanced triangle of one of the types 
$A,B,A'',B''$.
To find it we describe all geodesic loops of finite length
in $G$ for each type.
\begin{lemma}\label{u1}
For $T_1$ of type $A$, all elements of $G$ have finite length,
and they are in bijective
correspondence with indivisible elements\footnote{An element of $\Z^2$
is called indivisible
if it is not an integer multiple of any other element.}
of the lattice generated by
the sides of $T_1$.
So $s_1$ and $s_2$ are the smallest and the second smallest
lengths of sides of $T_1$.

For $T_1$ of type $B$, there are two elements in $G$ of finite length.
They correspond to the two sides of $T_1$ of finite length.

For $T_1$ type $A''$, let $a_3$ be the corner opposite to the side
of infinite length, and denote $a=f(a_2)-f(a_1),\; b=f(a_3)-f(a_1)$
Then the length spectrum is 
\begin{equation}\label{spectrum}
\{|na+b|:n\in\Z\}
\end{equation}
and $\{ s_1,s_2\}=\{|b|,|a-b|\}.$

For $T_1$ of type $B''$, 
there are exactly three elements of $G$ of
finite length. They correspond to the sides of $T_1$.
\end{lemma}

{\em Proof.} We look at all geodesics starting at
a vertex of $T_1$ with all possible slopes.

If $T_1$ is of type $A$, then the torus $\Phi(T_1)$
is a flat non-singular torus from Example~1
in Section~\ref{sBFT}.
Elements
of $G$ are in bijective correspondence with segments whose endpoints
are elements of the lattice and which contain no lattice points in their
interiors.

For $T_1$ of type $B$, the torus $\Phi(T_1)$ can be represented
as an exterior parallelogram $Q$ whose opposite sides
are identified by translations. Any geodesic starting from a vertex
which is not a side of $Q$, visits the pole and thus
has infinite length. The set $G$ contains two elements of finite length
(corresponding to two pairs of parallel sides of $Q$).

Let $T_1$ be the balanced extension of an unbalanced triangle
$T_2$ with angle sum $\pi$, and $f$ the developing map of $T_1$.
We normalize so that $f(T_2)=(0,1,b),$ so that $a=1$,
and $|b|\leq 1.$
The torus $\Phi(T_1)$
is partitioned into two triangles
congruent to $T_2$ and two half-planes. Suppose that a
geodesic starting from
a vertex $v$ is not a side of $T_2$. Then it visits the regions of
our partition, one after another. Once the geodesic enters a half-plane
it must stay there until it hits a pole, so the length of such geodesic is
infinite. A geodesic in $G$ of finite length must cross
the two finite sides of $T_1$ alternatively, and its length is given by
(\ref{spectrum}). See Fig.~4, which shows the images
under the developing map:
$f(T_1)$ is dark, the image of one congruent copy of $f(T_1)$
is grey, and images other congruent copies of $T_2$ are white triangles.
Images of several geodesics issued from one corner are shown:
the dotted lines are images of geodesics of
infinite length which contain poles, and dashed lines are images of
some
geodesics in $G$ of finite length. Since $|b|\leq 1$, the smallest
and the second smallest elements of (\ref{spectrum})
are $|b|$ and $|1-b|$.

Now consider a torus $\Phi(T_1)$ where $T_1$ is of type $B''$,
that is $T_1$ is a balanced extension of some unbalanced $T_2$ of type $B$.
Then $T_1$ is an exterior triangle (exterior of a bounded
triangle with angle sum $\pi$).
The torus $\Phi(T_1)$ 
is obtained by gluing two copies of exterior
triangles $T_1$. Every geodesic on this torus
which does not correspond to a side of
$B''$ passes through a pole.\hfill$\Box$

We conclude from Lemma~\ref{u1}:
\vspace{.1in}

\noindent
{\em For a balanced triangle $T$, the first and second
systoles of $\Phi(T)$ are two of the three numbers
$|f(a_i)-f(a_j)|$.}
\vspace{.1in}

It follows that
when a balanced triangle $T$ degenerates, then
the torus $\Phi(T)$ also degenerates, therefore the map $\Phi$ is proper,
and thus $\Phi^*$ is proper as well.
This completes the proof of
Proposition~\ref{prop5}.\hfill$\Box$

\section{Surjectivity of $\Phi^*$}\label{sth3}

We recall that $\T_m$ is the set of BFT with the angle sum $\pi(2m+1)$,
and $\T^*_m=\T_m/\sim$ is the quotient by the following equivalence relation:

(i) we identify triangles obtained from each other by a cyclic permutation
of the three vertices, and

(ii) we identify pairs of triangles described in Proposition \ref{prop3}.

In this section we prove that the map $\Phi^*_m:\T_m^*\to\L_m$
is surjective establishing the first part of
Theorem \ref{theorem3}. Injectivity will be proved in Section~\ref{spaces}.
\vspace{.1in}

{\em Proof of Surjectivity of $\Phi^*$.}
\vspace{.1in}

The plan of the proof is the following: for a given flat singular torus $L$
we find two special geodesic loops whose complement
is a quadrilateral $Q$. Then we construct cell decompositions $C_2,C_4$
of $Q$, and reassembling certain cells of $C_4$ we obtain
a decomposition of $L$ into two congruent balanced triangles.

Let $L$ be a torus with the singular point $O$.
Consider the germ at $O$ of the developing map $f:L\to\bC$, $f(O)=0$.
Let $g:\C\to L$ be a universal covering with $g(0)=O$. Then the
composition\footnote{Many authors call this $F$ a developing map.}
$$F=f\circ g$$
has a meromorphic continuation
to the whole plane. This meromorphic function
satisfies
\begin{equation}\label{period3}
F(z+\omega)=F(z)+\eta,\quad \omega\in\Lambda.
\end{equation}
Here $\omega\mapsto\eta(\omega)$ is a group homomorphism $\Lambda\to\C$,
and there are
two possibilities:
\vspace{.1in}

\noindent
{\em a) Generic case.} The image of $\Lambda$ is 
another lattice $\Lambda'\subset \C$, of rank $2$, and $F:\Lambda\to\Lambda'$
is an isomorphism, or
\vspace{.1in}

\noindent
{\em b) Degenerate case.} The image of $\Lambda$ belongs to a line through the
origin.
\vspace{.1in}

The pull back the flat metric via $F$
has the length element
$$ds=|F'(z)||dz|.$$
This metric has conic singularities at the critical points
of $F$ which are the points of
$\Lambda$, and some poles.

Let $\gamma_1$ be a shortest curve among all curves from $0$ to some point
$\omega\in\Lambda\backslash\{0\}$. We denote its endpoint other than $0$
by $\omega_1$. It is clear that $\gamma_1$ is a simple curve; $F(\gamma_1)$
is a segment $[0,\eta_1]$, where $\eta_1\in\Lambda'\backslash\{0\}$,
and the map $F:\gamma_1\to[0,\eta_1]$ is a homeomorphism.

Let $\gamma_2$ be a shortest of all curves from $0$ to some point
$\omega_2\in \Lambda\backslash\{\Z\omega_1\}$. The following lemma implies
that
$\gamma_1$ and $\gamma_2$ are disjoint except their common endpoint
at $0$.
\begin{lemma}\label{lemma31}
The curves $g(\gamma_1)$ and $g(\gamma_2)$ in $L$ intersect only at $O$.
\end{lemma}

{\em Proof.} Suppose that this is not so, and let $p\neq O$
be a point of intersection. Since $p$ is not a conic singularity,
and our curves are geodesic, they must make a non-zero angle st $p$.
It follows that the ratio of the periods of the differential $df$
over $g(\gamma_1)$ and $g(\gamma_2)$ is not real.

Now we construct a loop $\Gamma$ in $L$ based at $O$ which is shorter than
$g(\gamma_2)$ and whose homology class is not a multiple of $g(\gamma_1)$.
The point $p$ breaks $g(\gamma_1)$ into two arcs, and we denote the shorter
of these arcs by $I_1$. Similarly $I_2$ is the shorter of the two
arcs into which $p$ breaks $g(\gamma_2)$. Let $\Gamma$ be the concatenation
of $I_1$ and $I_2$. From our observations on the periods of $df$ we conclude
that 
$$\int_\Gamma df\neq 0,$$
therefore $\Gamma$ is non-trivial. Here we used that $df$ has no residues
at the poles. Moreover, this integral cannot be a real multiple of
the integral over $g(\gamma_1)$, so $\Gamma$ is not a
multiple of $g(\gamma_1)$. Finally, the length of $\Gamma$ is at most
the length of $g(\gamma_2)$, but $\Gamma$ can be shortened since it
has a non-zero angle at $p$, so we obtain a contradiction.\hfill$\Box$
\vspace{.1in}

The loops $g(\gamma_1)$ and $g(\gamma_2)$ cut the torus into
a quadrilateral. Preimage of this quadrilateral under $g$ is a quadrilateral
in the plane 
bounded by $\gamma_1,\gamma_2$ and their shifts $\gamma_1+\omega_2$
and $\gamma_2+\omega_1$. From Lemma~\ref{lemma31}
we conclude that all four curves are pairwise
disjoint except their endpoints. 

Thus we obtain a Jordan quadrilateral that will be called
$Q$ (the boundary is included).

Since the curves $g(\gamma_1)$ and $g(\gamma_2)$ have intersection index
$\pm1$, they generate the fundamental group of the torus, and it
follows that $Q$ does not contain other lattice points except
$0,\omega_1,\omega_2$ and $\omega_1+\omega_2$. 
The image $F(\partial Q)$ consists of $4$ straight segments
which form a parallelogram in the plane in the non-degenerate case.
In the degenerate case these $4$ segments belong to the same line.
Next we study 
\vspace{.1in}

{\em  Topology of the map $F:Q\to\bC$.} 
\vspace{.1in}

The following argument is purely topological, so we consider
an arbitrary Jordan quadrilateral $Q$ in the plane
(a closed disk with four
distinct marked boundary points
$a_1,\ldots,a_4$, which we call {\em corners},
enumerated
according to the standard orientation). The boundary arcs $(a_i,a_{i+1})$,
where $i$ is a residue modulo 4,
are called the {\em sides}.

Let $F:Q\to\bC$ be a continuous function which is a local homeomorphism
on the complement of the corners,
and topologically holomorphic\footnote{Topologically equivalent to
$z\mapsto z^{n_i}$, $1\leq i\leq 4$.}
at the corners.

About the boundary behavior we make
one of the two assumptions:
\vspace{.1in}

\noindent
a) Generic case:
$F(\partial Q)$ is a Jordan curve $\gamma$ and $F:\partial Q\to\gamma$
is a homeomorphism, or
\vspace{.1in}

\noindent
b) Degenerate case: the restrictions of $F$ to the sides
are homeomorphisms onto the image of each side, and
these images are segments of the same straight line $\ell$ in $\C$. 
The images of opposite sides have equal length.
\vspace{.1in}

We want to obtain a topological description
of possible partitions of $Q$
by $F^{-1}(\gamma)$ in case a) and by $F^{-1}(\ell)$ in case b).

First we address the generic case a). Consider the cell decomposition $C_1$
of the Riemann sphere which has two $2$-cells:
the interior $I$ and the exterior $E$
of
$\gamma$ (Fig.~5d).
The $0$-cells are $F(a_j)$ and $1$-cells are the four arcs
into which $F(a_j)$ divide $\gamma$.
We assign the labels to $0$- and $1$-cells by the
following rules: $F(a_j)$ has label $j$; the arc $(F(a_j),F(a_{j+1}))$
has label $j$.

Now consider the preimage $C_2=F^{-1}(C_1)$ in $Q$. 
Our assumptions about $F$ imply that $C_2$ is a finite cell decomposition
of $Q$.
It is called the {\em net} of $F$. 
Closures of the cells of $C_2$ are mapped onto the closures
of the cells of $C_1$ homeomorphically,
and we label cells of $C_2$ by their images.
Since $F$ is a local homeomorphism
on $Q\backslash\{ a_j\}$, the $1$-skeleton of $C_2$ consists of simple curves
which can meet only at the corners. We call the intersections
of these curves with the interior of $Q$ {\em arcs} and define the
length of an arc as the number of $1$-cells that it contains. An example
of the cell
decomposition $C_2$ is shown in Fig.~5a, where the black dots are $0$-cells.

The faces of $C_2$ are quadrilaterals, and we classify them as follows:

A face is called {\em lateral} if its boundary consists of one arc
of length $1$ and one arc of length $3$, both arcs having
as endpoints two adjacent corners of $Q$.

A face is called {\em diagonal} if its boundary consists of two arcs
of length $2$ both having as endpoints two opposite corners of $Q$.

A face is called {\em triangular} if its boundary consists of two arcs
of length $1$ and one arc of length $2$, arcs of length $1$ connecting
pairs of adjacent corners, while the arc of length 2
connects opposite corners.

A face is called {\em quadrilateral} if its boundary consists of $4$
arcs of length $1$, each connecting a pair of adjacent corners of $Q$.

Fig.~5a contains 8 lateral, 1 diagonal and 2 triangular faces.

Let us show that
this classification exhausts all possibilities
for the faces of $C_2$.
A face of $C_2$ cannot have all $4$ boundary vertices
in the interior of $Q$, since then there would be an adjacent face which
is not simply connected. Neither a face of $C_2$ can have two vertices
at the same corner, because the restriction of $f$ on the boundary
of a face is a homeomorphism onto $\gamma$. A face of $C_2$ cannot have
only one
vertex at a corner, because if this were the case, an adjacent $2$-cell
will have all its 4 boundary edges the same as the original face, which is
impossible.

\begin{lemma}\label{lemma32}
Under the assumption a) there are the
following possibilities:

(i) The net contains one quadrilateral face and some (possibly none)
lateral
faces.

(ii) The net contains at least one diagonal face, two triangular faces
and  several (possibly none) lateral faces.
All diagonal faces share the same
opposite corners on their boundaries.
\end{lemma}

This lemma and its proof are illustrated in Fig.~5. In Fig.~5b
case (i) is illustrated ($C_2$ is shown with bold lines).
Fig.~5a is an example of case (ii).
\vspace{.1in}

{\em Proof.} Notice that lateral faces come in pairs, so the number
of lateral faces
sharing two given corners $a_i,a_{i+1}$ on their boundaries must be even.
So
the innermost arc in $Q$, connecting $(a_i,a_{i+1})$, has length $1$.
Removing all lateral faces, we obtain a smaller quadrilateral
$Q'$, and a cell decomposition $C_3$ of it which has no lateral faces.
The restriction of $f$ to $Q'$ satisfies the
same conditions as $f$ on $Q$: the boundary $\partial Q$ is mapped
onto $\gamma$ homeomorphically.

If $C_3$ consists of a single face, we are in case (i).
If $C_3$ contains a diagonal face, suppose it has $a_1$ and $a_3$ on
the boundary. Then all diagonal faces must have $a_1$ and $a_3$
on their boundaries. Removing all of them, we obtain two triangular faces,
so we are in case (ii).

If $C_2$ contains no diagonal faces, then there are no
triangular faces. Indeed, suppose that the cell decomposition
of $Q'$ consists of just two triangles.
The $1$-cells on the boundary of each triangle have $4$ distinct labels,
and two of these $1$-cells are in the common boundary of these two
triangles.
But the $1$-cells on the boundary of $Q'$ also have $4$ distinct labels,
and this is evidently impossible.
Thus if there are no diagonal faces in $Q'$, then there are also no
triangular faces, and we are in case (i).\hfill$\Box$
\vspace{.1in}

{\em Transformation of the
cell decomposition $C_2$ into another cell decomposition
$C_4$ of $Q$.}
\vspace{.1in}

The edges of $C_4$ are defined as follows. First, they are
arcs of length $1$ of $C_2$. Then we discard all arcs of length at least $2$,
and add new edges by the following rules:

Suppose that $C_2$ has a diagonal face $G$.
We recall that cells of $C_2$ are labeled by their images in $C_1$.
If $G$ is a cell of $C_2$ is labeled $I$, we
draw the {\em diagonal}: the $F$-preimage in $G$ of that diagonal of
the parallelogram $I$ which has two corners of $Q$ as its extremities.
If $G$ is labeled $E$,
we use one of the two
{\em exterior diagonals} of the parallelogram $I$. An exterior
diagonal is the complement to a diagonal in the line
which contains this diagonal.
In Fig.~5b the added diagonal is red (dotted), and the discarded arcs are grey.

If $C_2$ has no diagonal faces, then it has one quadrilateral face.
If this quadrilateral face is labeled $I$,
we break this quadrilateral face by  the preimage of the
{\em shorter} diagonal of the
parallelogram $I$. If the quadrilateral face of $C_2$ is labeled $E$
we break this quadrilateral face by the preimage of the
exterior diagonal of $I$ which connects the two vertices of
this parallelogram with the larger exterior corner (Fig.~5g).
Partition of a quadrilateral face is shown if Fig.~5c.

By these rules, we obtain a cell decomposition $C_4$ of $Q$ which has
no interior vertices. This decomposition contains two congruent triangles
and a number of digons. Each digon is mapped to the sphere
with a cut along a segment. We break it into two digons by the $F$-preimage
of the complement of this segment in the line that contains it.
These lines are shown as red/dotted lines in Fig.~5b,c.
After these cuts are made, the number of digons in every ``bunch''
becomes even. 
Adding half of them to the adjacent side of the triangle
we obtain a decomposition of our torus into
two triangles. 

The final decomposition of $Q$ into two primitive triangles
and digons isometric to half-planes
is shown in red/dashed and bold black lines in Figs.~5b,c in two cases:
4b) when $C_2$ has a diagonal face, and 4c) when it does not.

Now we show that these two triangles are balanced.
We refer to the decomposition of a singular triangle described
in Proposition~\ref{prop2}.

Gluing any numbers of half-planes to the sides of a balanced triangle
results in a balanced triangle. Primitive triangles are
balanced in the following cases. Primitive triangle of the type A 
is balanced if all angles are less than $\pi/2$. If the greater angle
is $>\pi/2$ and at least one half-plane is glued to the opposite side,
the resulting triangle is balanced. If the cell decomposition $C_2$
contains a diagonal face, this implies that at least one half-plane
was glued opposite the largest angle of the triangular face.
If the triangular face
is of the type A, then its longest side is the diagonal, so
the largest angle is opposite to it.
If this triangle is of type B, then it is
balanced (a triangle of this type is always balanced).

If there was no diagonal in $C_2$, then we obtained triangular faces of
$C_4$ by drawing either the smaller diagonal in a parallelogram,
or the exterior diagonal in its exterior which has endpoints
at the bigger
exterior angles. In both cases the triangle is balanced.
\vspace{.1in}

So we obtained a partition of $Q$ into two balanced triangles.
We can re-assemble it by moving digons adjacent to a side of $Q$
to the opposite side to make the two balanced triangles congruent.
This completes the proof in the non-degenerate case.
\vspace{.1in}

Now we consider the degenerate case. $F:Q\to\bC$ maps the sides
of $Q$ into a line $\ell$, and we assume without loss of generality
that $\ell=\R\cup\{\infty\}$. The images of sides occupy some segment
$(a,b)\in\R$, where $a<b$. It is evident that $a$ and $b$ are the
images of two opposite corners of $Q$. Without loss of generality,
these corners are $a_1$ and $a_3$. 
The preimage
$F^{-1}(\R)$ defines a cell decomposition $C_5$ of $Q$.
It is exactly of the same type as nets studied in \cite{EG}:
they consist of simple curves with endpoints at the corners and
disjoint interiors,
and each curve is mapped homeomorphically onto its image.

\begin{lemma}\label{lemma33} Under these assumptions $C_5$ contains
a curve from $a_1$ to $a_3$. So the faces of $C_5$ are two triangles
and several (possibly none) digons.
\end{lemma}

{\em Proof.} Any component of $F^{-1}((\R\backslash[a,b])\cup\{\infty\})$
must be a curve from $a_1$ to $a_3$ in the interior of $Q$. This proves
the lemma.
\vspace{.1in}

Each digon of $C_5$ is mapped by $F$ to $\bC$ with a cut
(bounded or
containing $\infty$). We partition digons into two halves by complements
of these cuts to the $\R\cup\{\infty\}$. Then we split these half-planes
in each ``bunch'' into two equal parts and add them to the corresponding
sides of triangular faces. This defines a decomposition of
our torus into two triangles. That these triangles are balanced
is proved in the same way as in the non-degenerate case.
\hfill$\Box$

This completes the proof of surjectivity of $\Phi^*$.\hfill$\Box$
\vspace{.1in}

\section{The spaces $\A_m$, $\T_m$ and $\T_m^*$}\label{spaces}
\subsection{Connected components}

To visualize Proposition~\ref{prop1}, we introduce
the {\em space of angles} $\A_m$. In the intersection
of the plane 
\begin{equation}\label{plane}
P=\{\alpha\in\R^2:\alpha_1+\alpha_2+\alpha_3=2m+1\},
\end{equation}
with the open first octant in $\R^3$ (Fig.~9) we consider the triangle
$\Delta_m$ defined by the inequalities 
$$0<\alpha_j\leq\alpha_i+\alpha_k\quad\mbox{for all permutations}\quad (i,j,k);$$
it is shaded in Fig.~9.
The vertices of $\Delta_m$ are
$$(m+1/2,m+1/2,0), \quad (m+1/2,0,m+1/2), \quad (0,m+1/2,m+1/2).$$
Notice that the vertices do
not belong to $\Delta_m$ but the sides do belong, so $\Delta_m$
is neither open nor closed.

To obtain $\A_m$ we remove from $\Delta_m$ all lines where some $\alpha_j$
is an integer, and add all points where
all three $\alpha_j$ are integers. 
 The intersections of lines $\alpha_j=k$ with
$\Delta_m$ will be called {\em segments}. A segment is called even
or odd depending on the parity of $k$. There are three
families of parallel segments, each containing $m$ segments. Spaces of
angles for $m=1,\ldots,5$ are shown in Figs.~6, 7.

The set $\A_m$ has a natural partition
into open topological disks (faces) open intervals (edges)
and points (vertices): the faces are components
of the interior of $\A_m$ (they are triangles or quadrilaterals), the vertices
are the points where all $\alpha_j$ are integers, and the edges
are open intervals in $\A_m\cap\partial\Delta_m$. The set of
vertices of $\A_m$ will be denoted by $\V_m$.

We have a natural projection 
\begin{equation}\label{phi}
\varphi:\T_m\to\A_m
\end{equation}
which to every
balanced triangle with angles $(\pi\alpha_1,\pi\alpha_2,\pi\alpha_3)$
puts into correspondence the point
$(\alpha_1,\alpha_2,\alpha_3)\in\A_m\subset\R^3$. It follows
from Proposition~\ref{prop1} that this correspondence maps
the part of $\T_m$ where $\alpha_j$ are not integers bijectively
onto $\A_m\backslash\V_m$.
Triangles with integer angles
are mapped to $\V_m$
and the preimage of each point in $\V_m$
consists of three intervals.

This induces a partition of $\T_m$: the faces of $\T_m$ are $\varphi$-preimages
of the faces of $\A_m$, the edges are of two types: interior edges
which are mapped by $\varphi$ to the vertices of $\A_m$ and boundary edges
which are mapped bijectively onto intervals of $\A_m\cap\partial \Delta_m$.
There
are no vertices in this partition of $\T_m$.

The faces of $\T_m$ are adjacent when their images in $\A_m$
share a boundary vertex and
{\em their angles at this vertex are vertical\,\footnote{Opposite angles
among the four angles made by crossing of two lines.}.}
Notice that the map $\phi$ switches the orientation when one passes through
any interior edge of $\T_m$ from a face to an adjacent face.
This can be seen from the explicit formula for the angles in terms of
the conformal coordinate $z=\phi_{i,j,k}$ introduced in Section~\ref{complexanalytic}:
in the chart where $f(a_1)=0,\; f(a_2)=1,\; f(a_3)=z=x+iy$ we have
$$\alpha_1=p+\arctan(y/x),\quad \alpha_2=q+\arctan(y/(1-x)),$$
where $p,q$ are integers. Assuming that $x\in(0,1)$ we compute the Jacobian
and see that it switches sign simultaneously with $y$. 
\vspace{.1in}

\noindent
{\bf Remark.} Gluing of two $2$-cells along their common
boundary $1$-cell corresponding to a vertex of $\A_m$ 
reverses the natural
orientation of these $2$-cells induced from the
$(\alpha_1,\alpha_2)$-plane. Nevertheless,
it is easy to check that the surface $\T_m$ is orientable. To
do this one paints the $2$-cells of $\T_m$ into two colors,
so that each two $2$-cells with a common vertex have different
colors. It is clear that such a coloring is possible, see Fig.~8.
\vspace{.1in}

To study the surface $\T_m$ we introduce the graph $\Gamma_m$, which will
be called the {\em nerve}.
Examples of these graphs are shown in figures 6, 7.
Their vertices correspond to $2$-cells of $\T_m$ (or faces of $\A_m$)
and two vertices of $\Gamma_m$ are
connected by an edge if the corresponding two $2$-cells of $\T_m$ share
an edge or, which is the same, if the corresponding $2$-cells
of $\A_m$ share a vertex and their angles at this vertex are vertical.
Then we find 

\begin{prop}\label{prop6}
When $m\geq 2$, the graph $\Gamma_m$ has
$4$ connected components. Exactly one of them, $\Gamma_m^\prime$ is
invariant under the order $3$ rotation about the center
of $\A_m$. Three others are permuted by this rotation.
$\Gamma_1$ has only three components, permuted by
the rotation. $\Gamma_0$ consists of one vertex only.\hfill$\Box$
\end{prop}

In Fig.~6, $\Gamma_m^\prime$ is blue/dotted, while in Fig.~7
one of the three components permuted by the order $3$ rotation
is red, any of these three components is called $\Gamma_m^{\prime\prime}$
(they are isomorphic graphs embedded in the plane).

We give first a geometric sketch which makes our proposition evident,
and then a more formal proof.
\vspace{.1in}

{\em Sketch of a proof of Proposition~\ref{prop6}.} Let us consider the plane
$P$ which is defined in (\ref{plane}).
Intersections of $P$ with the planes $\{\alpha_j=\mbox{integer}\}$ break $P$
into triangles. 
Connecting the centers of pairs of triangles which share a vertex
and whose angles at this vertex are vertical, we obtain four honeycomb
structures $X_j$ with disjoint vertices. See Fig.~8 which shows two
of these honeycombs.
Choosing one vertex of one honeycomb as a center, we see that this honeycomb
is invariant under rotation by $120^o$ about this vertex,
while the other three are permuted.
\vspace{.1in}

{\em Proof of Proposition~\ref{prop6}.}
\vspace{.1in}

Let $\Delta_m^\prime$ be the intersection of the plane
$$\alpha_1+\alpha_2+\alpha_3=2m+1$$
with the closed first octant
$\alpha_j\geq 0,\; 1\leq j\leq 3.$
The {\em segments}
$$\{(\alpha_1,\alpha_2,\alpha_3)\in\Delta_m^\prime:\alpha_j=k\},
\quad 1\leq j\leq 3,\quad 0\leq k\leq 2m$$
divide $\Delta_m^\prime$ into open triangles which we call {\em faces}
of $\Delta_m^\prime$. There are three families of 
these segments, depending on the value of $j$ which were discussed
in the beginning of this section.
We recall that a segment with even/odd $k$ is
called {\em even/odd}.
Since $\alpha_1+\alpha_2+\alpha_3$ is odd, among three segments intersecting
at an integer vertex, either one or all of them are odd. This implies
that a face has either one or three sides on even segments.

We classify faces of $\Delta_m^\prime$ into four types:
\vspace{.1in}

\noindent
type $I$, if all three sides of the triangular face belong to even segments,
\vspace{.1in}

\noindent
type $II_j,\; j\in\{1,2,3\}$
if one side belongs to an even segment of family $j$, while two
other sides belong to odd segments of the other two families.
\vspace{.1in}

Two faces sharing an integer vertex are called {\em vertical}
if their
angles at this vertex are vertical (opposite).

We claim that vertical faces are of the same type.
Indeed if the two segments bounding the vertical angles are both
even, then each of our two faces must have all three sides on even segments,
thus both faces are of type $I$. If exactly one of the segments
is even, and belongs to family $j$, then both faces are of the type
$II_j$. If both segments are odd, then the sides of our faces
opposite to the considered vertex are even and parallel, so they
are in the same family and our two faces are in the same family.
This proves the claim.

Now we claim that faces of the same type cannot have a common side. 
If two faces have a common side on an even segment then one of them is
of type $I$ and another is of type $II$. If the common side is
on an odd segment then both faces are of type $II$ and their sides
on even segments 
are not parallel. Thus they belong to different types.
This proves the claim.

Our next claim is that the closure of the union of faces of the same type
is connected. Consider faces of type $I$. The closure of their union
consists of the faces themselves and all even segments. It is clear
that the union of all even segments is connected.

The proof for other types is similar:  the closure of the
union of faces of type $II_j$ consists of the faces themselves,
the even segments of family $j$,
and odd segments of two other families. The union of these 
segments is connected.

If we restrict now to $\A_m$ and consider the union of those faces
of a single type which intersect $\A_m$ and their vertices,
this union is still connected. Indeed if two faces of one family
share a vertex and both intersect $\A_m$ then this common vertex belongs
to $\A_m$.

This proves Proposition~\ref{prop6}.\hfill$\Box$
\vspace{.1in}

Now consider the map $\varphi:\T_m\to\A_m$. Component $\L_m^I$ of $\T_m$
consists of preimages of faces of type $I$ of $\A_m$ and common
edges of pairs of these preimages that project to the vertices of
faces of type $I$ of $\A_m$.
Similarly for $i\in\{1,2,3\}$, components $II_i$ of $\T_m$,
consist of preimages of faces of type $II_i$ and common edges of pairs of
these preimages which project to the vertices of faces of type $II_i$
of $\A_m.$ So we obtain

\begin{cor}
When $m\geq 2$, $\T_m$ consists of four connected components
\newline
$I,II_1,II_2,II_3$.
Cyclic permutation of vertices preserves component $I$ and permutes
components $II_i$. As a consequence, $\T_m^*$ has two components, $I$ and $II$.
These components
are distinguished by the number of sides with poles:

When $m$ is even,
$\L_m^I$ consists of BFT with no poles on the sides, and 
$\L_m^{II}$ consists of BFT with two poles on the sides. 

When $m$ is odd, $\L_m^I$ consists of BFT with 3 poles on the sides,
and
$\L_m^{II}$ consists of BFT with one pole on
the side. \hfill$\Box$ 
\end{cor}

\noindent
{\bf Example.}
Figure 2 shows all types of BFT for $m=2$ (angle sum $5\pi$).
Triangles of types a), b), c) belong to $\L_m^I$.
Suppose for example that triangle a) is deformed so that the top
vertex moves towards the opposite (horizontal) side. Eventually we obtain
a triangle b) with the angles $2\pi,2\pi,\pi$. If the middle vertex of b)
continues
moving downwards, we obtain triangle c). Its developing map
is $2$-to-$1$ onto the darkly shaded region and $1$-to-$1$ onto the lightly
shaded region. There are three types of such triangles c) if the
vertices are labeled, but only one type with unlabeled vertices.

Triangles of types d), e), f) belong to $\L_m^{II}$. It is easy to
visualize how they are deformed to each other.

Triangles a), b), c) have one pole inside, while triangles d), e), f)
have 2 poles, one on each of the two unbounded sides, and the third side
is free of poles. Fig.~2 should be compared with Figs.~6, 7, $m=2$:
the set $\A_2$ shows the location of all these triangles
in the parameter space.

\subsection{Proof of injectivity of $\Phi^*$}

We established in Section~\ref{complexanalytic}
that $\Phi_m^*:\T_m^*\to\L_m$
is a proper holomorphic map between punctured Riemann surfaces. So to prove injectivity
it is sufficient to show that every component of $\T^*_m$ contains an interval $I$
such that for every $T\in I$,
$T$ is the unique $\Phi^*$-preimage of $\Phi^*(T)$.

Assume that $m\geq 1$. Then every component contains a triangle
with integer angles.
Let $T$ be a triangle with integer angles.
According to Proposition~\ref{prop2.5}
the monodromy group of $\Phi^*(T)$ is a subgroup of a line,
so for the torus $\Phi^*(T)$ alternative
b) (degenerate case) holds in the proof of surjectivity.
It follows that any $\Phi^*$-preimage of $\Phi^*(T)$ also
has all integer angles (see Proposition~\ref{prop2.5}).

Recall that a triangle with integer angles is obtained from a triangle
in Fig.~1b in which $a_1,a_2,a_3$ belong to the same line, by gluing
half-planes to the sides (see Proposition~\ref{prop2}).
Let us normalize so that $a_3=0$, $a_1=1$, and we choose our interval $I$
so that
$a_2:=a\in(0,1/2)$. Then the following properties
of $\Phi^*(T)$ are evident:

The shortest non-trivial loop
$\gamma_1$ based
at the singularity
has length $a_2$. We define the orientation
of this loop by orienting $(0,a)$ from $a$ to $0$
on the boundary of the reduced triangle $T'$.
We recall that the reduced triangle $T'$ is the upper half-plane
with corners at $0,a,1$.

So the parameter $a$ is uniquely defined by $\Phi^*(T)$,
as the shortest length of a loop based at $O$ on $\Phi^*(T)$.
Now we define $\gamma_2$ as the loop corresponding to
$(a,1)$ in $\partial T'$, oriented from $a$ to $1$.
This loop $\gamma_2$ is characterized as the shortest loop
whose class does not belong to $\Z\gamma_1$.

Each side of $T'$
defines a homotopy class of loops based at $O$. Two of them
are $\gamma_1$ and $\gamma_2$. 
Suppose that $m_1$ half-planes were glued to the side $(0,a)$,
and $m_2$ half-planes were glued to the side $(a,1)$.
Then the torus $\Phi^*(T)$ contains $m_1+1$ disjoint (except the base point)
geodesic loops in the class $[\gamma_1]$, and  $m_2+1$ disjoint
geodesic loops in the class $[\gamma_2]$. 

This implies that the angles $\pi\alpha_i$ of $T$ are defined by
the properties of the torus $\Phi^*(T)$, namely $\alpha_i=1+m_j+m_k$.
This proves injectivity of
the map $\Phi^*$ and completes the proof of Theorem~\ref{theorem3}.

\section{Euler characteristics of components of $\L_m$
and completion of the proof of Theorem~\ref{theorem1}}\label{ident}

Theorem~\ref{theorem3} reduces the study of topology of $\L_m$ to the study of
topology of $\T_m^*$. 

It is convenient to use the nerves
$\Gamma_m^\prime$ and $\Gamma_m^{\prime\prime}$ 
introduced in Section~\ref{spaces}. 

First we recall that $\Gamma_m/\Z_3$ consists of $2$ components. One of them
comes from the component $\Gamma_m^\prime$ which is invariant with
respect to the $\Z_3$ action.
This is our component $\L_m^I$. Component $\L_m^{II}$
comes from
the three components of $\Gamma_m^{\prime\prime}$ which are permuted by the $\Z_3$ action.
See
Figs. 6, 7. 
\vspace{.1in}

\noindent
{\em Computation of the Euler characteristic for component $\L_m^I$.}
\vspace{.1in}

Let $\Gamma_m^\prime$ be the component of $\Gamma_m$ which is invariant
with respect to the $\Z_3$ action.
The numbers $\epsilon_0,\epsilon_1$ are defined in (\ref{epsilon0})
and (\ref{epsilon2}).
 We add to them
$$\epsilon_2=\left\{\begin{array}{ll}0,&\mbox{if}\; m\equiv 1\; (\mod 2)\\
1,&\mbox{if}\; m\equiv 0\; (\mod 2).\end{array}\right.$$
%$\epsilon_2\in\{0,1\},$
and interpret these numbers in terms of 
$\A_m$: 
\vspace{.1in}

$\epsilon_0=1$ if the center of $\A_m$ belongs to a $2$-cell.
Equivalent condition is that a vertex of $\Gamma_m^\prime$
is fixed by the $\Z_3$
action.
\vspace{.1in}

$\epsilon_1=1$ if there is a vertex of $\Gamma_m^\prime$ representing
a face which has the
middle of the side of $\Delta_m$ on the boundary. 
\vspace{.1in}

$\epsilon_2=1$ if there is a vertex of $\Gamma_m^\prime$ representing
a face which has a corner
of $\Delta_m$ on the boundary. 
\vspace{.1in}

\vspace{.1in}

We introduce further notation:

$V_1,V_2,V_3$ are the numbers of vertices of $\Gamma_m^\prime$ of
degrees $1,2,3$, and $V$ is the total number of vertices.

$E$ is the number of edges of $\Gamma_m^\prime$.

Taking into account all identifications on $\Gamma_m^\prime$, we obtain
the following formula for the Euler characteristic:
\begin{equation}\label{chiI}
\chi(\L_m^I)=
V_3/3+V_1/6+V_2/6-E/3+2\epsilon_0/3+(\epsilon_1/2-\epsilon_2)/2.
\end{equation}
%Begin insert2
To explain this formula, we compute the contributions to the
(ordinary) Euler characteristic of $\L_m^I$.

The group $\Z_3$ fixes the center of $\A_m$ which belongs to $\Gamma_m^\prime$
if and only if $\epsilon_0=1$, and acts freely on the set of remaining
$3$-valent vertices, so the number of $2$-cells in $\L_m^I$ corresponding to
$3$-valent vertices is $V_3/3+2\epsilon_0/3$.
Also, the group $\Z_3$ acts freely on
the set of edges, and thus there are $E/3$
corresponding $1$-cells in $\L_m^I$.

For a $1$-valent or $2$-valent vertex $v$ of $\Gamma_m^\prime$, its
class under the equivalence relation consists of
\begin{itemize}
\item[a)]
6 elements, if $v$ is neither a corner of $\Delta_m$ nor a midpoint
of a side of $\Delta_m$;
\item[b)]
3 elements, if $v$ is a corner of $\Delta_m$;
\item[c)]
3 elements, if $v$ is the midpoint of a side of $\Delta_m$.
\end{itemize}
If we have $V^a$ vertices of type a), such vertices correspond to
$V^a/6$ $2$-cells in $\L_m^I$.`The number $V^b$ of vertices
of type b) can be $0$ or $3$: the latter case occurs if and only if 
$\epsilon_2=1$, in which case such vertices correspond to a punctured disk in
$\L_m^I$. In both cases, vertices of type b) contribute $V_b/6-\epsilon_2/2=0$
to the Euler characteristic. Similarly, the number $V^c$ of vertices
of type c) can be $0$ or $3$: the latter case occurs
if and only if $\epsilon_1=1$, in which case such vertices correspond to a disk
in $\L_m^I$. In both cases, vertices of type c)
contribute $V^c/6+\epsilon_1/2=1$
to the Euler characteristic. Adding these contributions, we obtain (\ref{chiI}).
%End insert2

It is easy to see that for odd $m$
$$E=3(m^2-1)/8,\quad V_1=3(m-1)/2,\quad V_2=0,\quad V=(m^2+4m-5)/4,$$
and $V_3$ can be computed by the formula $V_3=V-V_1-V_2$.
This gives the formula for $\chi(\L_m^I)$ when $m$ is odd.
When $m$ is even, we have
$$E=3(m^2+2m)/8,\quad V_1=3,\quad V_2=3(m-1)/2,\quad V=(m/2+1)^2,$$
and again $V_3=V-V_1-V_2$. This gives the formula for $\chi(\L_m^I)$ when $m$
is even. The resulting formulas for $\chi$ in terms of
$m$ are written in the Appendix. Expressing $m$ in terms
of $d_m^I$ in (\ref{degFI})
and subtracting the orbifold correction we obtain
$\chi^O(\L_m^I)=-(d_m^I)^2/6.$
\vspace{.1in}

{\em Computation of the Euler characteristic of component $\L_m^{II}$}
\vspace{.1in}

Let $\Gamma_m^{\prime\prime}$ be one of the three components of $\Gamma_m$
which are permuted by the $\Z_3$ action. We use the following notation

$E$ is the number of edges of $\Gamma_m^{\prime\prime}$

$V_1,V_2,V_3$ and $V$ are the numbers of vertices of $\Gamma_m^{\prime\prime}$
of orders $1,2,3$ and the total number of vertices.

$\epsilon_1$ and $\epsilon_2$ have the same meaning as before.
%Begin insert 3
A computation analogous to that for $\L_m^I$ gives
%For the Euler characteristic we have in this case
%$$\chi(\L_m^{II})=E+(3/2)V_1+V_2/2-2V+(\epsilon_2-\epsilon_1)/2.$$
$$\chi(\L_m^{II})=V_3+V_1/2+V_2/2-E+(\epsilon_2-\epsilon_1)/2,$$
as vertices of type b) (respectively, of type c)) belong to
$\L_m^{II}$ if and only if $\epsilon_2=0$ (respectively $\epsilon_1=0$).
%End insert 3
When $m$ is odd,
$$E=(3m^2+4m+1)/8,\quad V_1=(m+3)/2,\quad V_2=m-1,\quad V=(m^2+4m+3)/4.$$
When $m$ is even,
$$E=(3m^2+2m)/8,\quad V_1=m,\quad V_2=m/2,\quad V=(m/2+1)^2-1.$$
This gives the formulas for $\chi(\L_m^{II})$
in terms of $m$ and $\epsilon_1$
(written in the Appendix). Expressing $m$ in terms of $d_m^{II}$
from (\ref{degFII}) and subtracting the orbifold correction
we obtain $\chi^O(\L_m^{II})=-(d_m^{II})^2/18$.
\vspace{.1in}

{\em Computation of the number of punctures.}
\vspace{.1in}

Consider a small simple loop around a puncture of $\T^*_m$. This loop projects
to a contour in $\A_m$ which goes close to 
the lines $\alpha_j=k$, switching the side at each integer point.
For component $\L_m^I$, the contour goes near lines with the same even $k$,
and $j=1,2,3$ and closes. See Fig.~10. So there is a $1-1$ correspondence between these
contours and triples of segments (one in each family) with even $k$.
So there are $\lfloor m/2\rfloor$ of such loops. In addition,
when $m$ is even there is a puncture corresponding to the vertices
of $\Delta_m$. Thus the total number of punctures on $\L_m^I$
is $m/2+1$ when $m$ is even and $(m-1)/2$ when $m$ is odd.
In other words, the number of punctures on Component $\L_m^I$ equals 
\begin{equation}\label{holesI}
h_m^I:=d_m^I,
\end{equation}
where $d_m^I$ was defined in (\ref{degFI}).

For component $\L_m^{II}$, the computation is similar, see Fig.~11.
Each contour goes either near an even segment, in which case it closes
after describing three segments, one of each family. If a contour accompanies
an odd segment, it ends on the other side of the odd segment after
describing three segments. So the total number of contours is
$m$ when $m$ is even and $m+1$ when $m$ is odd (the extra puncture for odd $m$
coming from the corners of $\Delta_m$), in other words
\begin{equation}\label{holesII}
h_m^{II}:=2\lceil m/2\rceil=2d_m^{II}/3.
\end{equation}
\begin{center} Component $\L_m^I$
\vspace{.1in}\nopagebreak

\nopagebreak
\begin{tabular}{|c|c|c|c|c|c|c|c|c|c|c|c|c|}\hline
$m$&$\epsilon_0$&$\epsilon_1$&$\epsilon_2$&$V_1$&$V_2$&$V_3$&$E$&$V$&$\chi$&$h$&$g$&$d$\\
2&1&0&1&3&0&1&3&4&0&2&0&2\\
3&1&1&0&3&0&1&3&4&1&1&0&1\\
4&0&1&1&3&3&3&9&9&-1&3&0&3\\
5&1&0&0&6&0&4&9&10&0&2&0&2\\
6&1&0&1&3&6&7&18&16&-2&4&0&4\\
7&0&1&0&9&0&9&18&18&-1&3&0&3\\
8&1&1&1&3&9&13&30&25&-3&5&0&5\\
9&1&0&0&12&0&16&30&28&-2&4&0&4\\
10&0&0&1&3&12&21&45&36&-6&6&1&6\\
11&1&1&0&15&0&25&45&40&-3&5&0&5\\
12&1&1&1&3&15&31&63&49&-7&7&1&7\\
13&0&0&0&18&0&36&63&54&-6&6&1&6
\end{tabular}
\end{center}
\vspace{.1in}

\begin{center} Component $\L_m^{II}$\nopagebreak 
\vspace{.1in}\nopagebreak
\nopagebreak

\begin{tabular}{|c|c|c|c|c|c|c|c|c|c|c|c|}\hline
$m$&$\epsilon_1$&$\epsilon_2$&$V_1$&$V_2$&$V_3$&$E$&$V$&$\chi$&$h$&$g$&$d/3$\\
1&0&0&2&0&0&1&2&0&2&0&1\\
2&0&1&2&1&0&2&3&0&2&0&1\\
3&1&0&3&2&1&5&6&-2&4&0&2\\
4&1&1&4&2&2&7&8&-2&4&0&2\\
5&0&0&4&4&4&12&12&-4&6&0&3\\
6&0&1&6&3&6&15&15&-4&6&0&3\\
7&1&0&5&6&9&22&20&-8&8&1&4\\
8&1&1&8&4&12&26&24&-8&8&1&4\\
9&0&0&6&8&16&35&30&-12&10&2&5\\
10&0&1&10&5&20&40&35&-12&10&2&5\\
11&1&0&7&10&25&51&42&-18&12&4&6\\
12&1&1&12&6&30&57&48&-18&12&4&6
\end{tabular}
\end{center}
We include two tables for $1\leq m\leq 13$. 
Notation, besides that already introduced is: $g$ for the genus,
$h$ for the number of punctures, $d$ for the degree of the forgetful map
as in (\ref{degFI}), (\ref{degFII}). 
The formulas for the degrees follow from \cite[sections 23.21-23.24]{WW}.
\vspace{.1in}

{\bf Remarks.} There is an alternative method of counting the punctures,
based on the description on compactifications of the spaces
of Abelian differentials in \cite{Ba,Ba2}. In recent preprints
\cite{M1}, \cite{M2} a general method of computation
of Euler's characteristics
for spaces of Abelian differentials with prescribed multiplicities
of zeros and poles is developed. However, our results do not
follow from the results stated in these preprints, mainly because
of the additional condition that residues vanish.
\vspace{.1in}

{\em Orbifold points.}
\vspace{.1in}

By definition, an orbifold point in $\L_m$ is a point which corresponds
to a flat singular torus with a non-trivial automorphism. 
An automorphism here means an orientation-preserving isometry.
The trivial automorphism is the involution which exists on every flat singular
torus. There are two types of tori with non-trivial automorphisms:
hexagonal ones with an automorphism
of order $3$, and square ones, with non-trivial automorphism of order $4$.
In the representation of tori as $\Phi^*(T)$, hexagonal tori correspond to
triangles whose all angles are equal, while square tori correspond to
marginal triangles whose two smaller angles are equal.
So in the space of angles $\A_m$, the hexagonal torus arises from the center
of $\Delta_m$ when this center belongs to $\A_m$, and the square torus
corresponds to the middles of the sides of $\Delta_m$. In Figs.~6, 7,
these points are denoted by little circles in the center of the picture,
and little black triangles in the middles of the sides.

In the next section we will use the following
\begin{prop}\label{special}
In the Lam\'e equation (\ref{lame1}) or (\ref{lame2})
corresponding to a hexagonal
or square torus
(in the metric sense), the accessory parameter $\lambda$ is
equal to $0$ (see
the text after (\ref{scaling})).
\end{prop}

{\em Proof.} Since a metric automorphism is also a conformal automorphism,
it corresponds to an automorphism of
the Lam\'e equation, that is to a fixed point of
transformation (\ref{scaling}). For both fixed points we have $\lambda=0$.
\hfill$\Box$

\section{Theorem \ref{theorem2} and Maier's conjecture}\label{nons}

To prove Theorem~\ref{theorem2} and its corollaries we first state
the exact relation between $\L_m$ and $\H_m$.

A {\em marked} elliptic curve is an elliptic curve on which the three
points of (exact) order $2$ are labeled. Legendre's family (\ref{le})
parametrizes marked elliptic curves: the labels are $0,1,a$.
The permutation group $S_3$ acts on the space of marked elliptic curves
by permuting the labels. Explicitly, the orbit of $a$ under this
action is
\begin{equation}\label{cross}
a,\quad 1-a,\quad 1/a,\quad 1-1/a,\quad 1/(1-a),\quad a/(a-1).
\end{equation}
This action lifts to the moduli space $\C\times\C_a$ of Lam\'e equations
in the form of Legendre: the generators $a\mapsto 1-a$ and $a\mapsto1/a$
lift to
$$(B,a)\mapsto\left(-B-m(m+1),\; 1-a\right),\quad\mbox{and}$$
$$(B,a)\mapsto\left(B/a,\;1/a\right).$$
To obtains these two transformations,
one changes the  independent variable in
(\ref{legendre}) 
$z\mapsto 1-z$ and $z\mapsto z/a$, respectively.
Taking the quotient of the space $\C\times\C_a$ of equations
(\ref{legendre}) by this $S_3$ action
we obtain an orbifold covering $\Psi_m$ of degree $6$
from the moduli space
of equations (\ref{legendre}) to the moduli space $\Lame_m$,
such that the following diagram is commutative:
\begin{equation}\label{cd}
\begin{CD}
\H_m^j    @>{\displaystyle\Psi_m^j}>>   \L_m^K\\
@VV{\displaystyle\sigma_m} V           @VV{\displaystyle\pi_m} V\\
\C_a      @>{\displaystyle\psi}>>         \C_J
\end{CD}
\end{equation}
Here $\Psi_m^j$ are restrictions of $\Psi_m$ on $\H_m^j$, and $K=I$
for $j=0$, $K=II$ for $j\in\{1,2,3\}$.
(We have not proved yet that $\H_m^j$ are irreducible; this will be
done only in the end of this section).

The explicit expression of $\psi$ is in (\ref{Ja}). To obtain an explicit
expression of $\Psi_m$ we change the independent variable $z$ in the equation
(\ref{legendre}) to $z-(1+a)/3$. Then we easily obtain $\Psi_m=(R_1,R_2,R_3)$
modulo scaling (\ref{scaling}), where
\begin{eqnarray}\label{R}
\lambda&=&R_1(B,a):=B+m(m+1)(a+1)/3,\\
\nonumber
g_2&=&R_2(B,a):=4(a^2-a+1)/3,\\
\nonumber
g_3&=&R_3(B,a):=8(a^3-3a^2/2-3a/2+1)/27.
\end{eqnarray}
We define compact Riemann surfaces $\overline{\L}_m^K$, $\overline{\H}_m^j$,
$\bC_J$ and $\bC_a$ by filling the punctures.
Later we will endow them
with orbifold structures.
The forgetful maps $\pi_m,\sigma_m$ and maps $\psi,\; \Psi$ extend uniquely to these compactifications.
\def\bL{\overline{\L}}
\def\bH{\overline{\H}}
\begin{definition}\label{def3} A point $x\in\bL_m$ is called special if
$\pi_m(x)\in\{0,1,\infty\}$. A point $x\in\bH_m$ is called special
if
$$\sigma_m(x)\in\{0,1,\infty,2,1/2,-1,(1\pm i\sqrt{3})/2\}.$$
\end{definition}
We will later show (in the proof of Corollary 1.1 in this section)
that $\Psi_m^j:\H_m^j\to\L_m^j$ are orbifold coverings,
so the maps $\Psi_m^j:\bH_m^j\to\bL_m^K$, as maps between Riemann
surfaces, can be ramified only at special
points.

Next we study ramification properties of forgetful maps at the special
points. For this we need two lemmas,
the first one is classical, see for example
\cite[Ch. II, \S 1, Thm 1]{GK}:
\begin{lemma}\label{lemma1}
 Let $A=(a_{i,j})$ be an $n\times n$ matrix with
$a_{i,i+1}>0,\; 1\leq i\leq n-1,$ and $a_{i,i-1}>0,\;2\leq i\leq n,$
the rest of the entries are zeros. Then all roots of the characteristic
polynomial are real and simple. The characteristic polynomial
is either even or odd, in other words it has the form
$\lambda^kP(\lambda^2),$ where $k\in\{0,1\}$, and $P$ is a real
polynomial.\hfill$\Box$
\end{lemma}

The second lemma was communicated to us by V. Tarasov; it is inspired by
\cite[Prop. 3]{ST}:
\begin{lemma}\label{lemma2}
Let $A=(a_{i,j})$ be an $n\times n$ matrix with
$a_{i,i+1}>0,\; 1\leq i\leq n-1$ and $a_{i,i-2}>0,\; 3\leq i\leq n$,
the rest of the entries are zeros. Then all roots of the characteristic
polynomial, except possibly $0$, are simple and their arguments
are of the form $2\pi k/3,\; k\in\{0,1,2\}$. In fact this characteristic
polynomial has the form $\lambda^kP(\lambda^3)$, where $k\in\{0,1,2\}$
and $P$ is a real polynomial with all roots positive.
\end{lemma}

A proof of Lemma~\ref{lemma2} will be given in the next section.

The following proposition lists ramification of
forgetful maps over special points. We use the word ``ramification''
in the sense of maps between Riemann surfaces, not orbifolds.
\begin{prop}\label{prop-ramJ}
\noindent
{\bf 1.} Ramification of $\pi_m^K$ over special
points is the following:

\noindent
Over $J=0$ there are $\lfloor d_m^K/3\rfloor$ triple points,
and one additional point $x$ when $d$ is not divisible by $3$.
This additional point $x$ is the orbifold point of order $3$,
and $\pi_m^K$ has $x$ as a double point when $d\equiv 2\; (\mod 3),$
and a simple point when $d\equiv1\; (\mod 3)$. 

\noindent
Over $J=1$ there are $\lfloor d_m^K/2\rfloor$ double points,
and one simple point when $d_m^K$ is odd. This simple point is
the orbifold point of order $2$.

\noindent
Over $J=\infty$ there are $d_m^{II}/3$ double points when $K=II$. The rest
$d_m^{II}/3$ points are simple. For $K=I$ all points over $\infty$
are simple.
\vspace{.1in}

\noindent
{\bf 2.} Ramification of $\sigma_m^j$ over special points
is the following:
Over each $a=1/2\pm i\sqrt{3}/2$,
there is one double point when $d_m^j\equiv 2\; (\mod 3).$
There is no other ramification over special points.
\end{prop}

{\em Proof.} For component $\L_m^I$ with even $m$ and $J=0$, we consider
polynomial solutions $Q$ of equation (\ref{lame1}) with
$g_2=0,\; g_3=1$, that is
$$(4x^3-1)Q''+6x^2Q'-m(m+1)xQ=\lambda Q.$$
The matrix of the linear operator in the left-hand side in the
basis of monomials has the form as in Lemma \ref{lemma2}.
Therefore the characteristic polynomial of this
matrix has the form $\lambda^kP(\lambda^3)$.
This has a root of multiplicity $2$ at $0$ when $d\equiv 2\; (\mod 2)$.
Other roots come in triples, each triple
lies on the same orbit under the $\C^*$ action (\ref{scaling}),
so we have $\lfloor d_m^K/3\rfloor$
triple points.

Similar considerations apply to other special points.

For component $\L_m^I$ with odd $m$ and $J=0$, we consider
solutions of (\ref{lame1})
of the form $\sqrt{4x^3-g_2x-g_3}\, Q(x)$, where $Q$ is a polynomial.
The equation for $Q$ becomes
$$(4x^3-1)Q''+18x^2Q'+\left(12-m(m+1)\right)x Q=\lambda Q,$$
and this leads to a matrix of the same form
described in Lemma~\ref{lemma2},
so the same argument as in the case of even $m$ applies.

For component $\L_m^I$ with even $m$ and $J=1$, we set $g_2=1,g_3=0$, and
obtain 
$$(4x^3-x)Q''+(6x^2-1/2)Q'-m(m+1)xQ=\lambda Q$$
which leads to a matrix described in Lemma~\ref{lemma1}.
The characteristic polynomial is of
the form $\lambda^kP(\lambda^2),\; k\in\{0,1\}$
which has one simple root $\lambda=0$ when $k=1$ and other roots come
in pairs which are on the same orbit under the $\C^*$ action.

For component $\L_m^I$ with odd $m$ and $J=1$ we obtain the equation
$$(4x^3-x)Q''+(18x^2-3/2)Q'+\left(12-m(m+1)\right)x Q=\lambda Q$$
which again leads to a matrix described in Lemma~\ref{lemma1}.
The conclusion
is similar.

For component $\L_m^{II}$ we use the Legendre's
form of Lam\'e equation (\ref{legendre}).
When $m$ is odd, and $J=1$, we plug the solution of the form
$\sqrt{z}\, Q(z)$ and obtain
\begin{eqnarray*}
&&4z(z-1)(z-a)Q''+(10z^2+8z(1+a)+6a)Q'\\
&-&\left((m^2-m-2)z+1+a+B\right)Q=0.
\end{eqnarray*}
When $m$ is even, and $j=1$, we plug the solution of the form
$\sqrt{z(z-1)}\, Q(z)$ and obtain
\begin{eqnarray*}
&&4z(z-1)(z-a)Q''+(14z^2-(12a+8)z+6a)Q'\\
&-&\left((m^2+m-6)z+B+4a+1\right)Q=0.
\end{eqnarray*}
Both these equations lead to Jacobi matrices as in Lemma~\ref{lemma1}.

To study ramification at the punctures, we use again Legendre's form.
Take, for example, $a=0$.
The matrix of the operator in the left-hand side of
(\ref{lame2}) is triangular, with distinct eigenvalues.
So $\sigma_m$ is unramified at a point $x$ with $\sigma_m(x)=0$.
Now we have $\deg_0\psi=2$, so by (\ref{cd})
$$\deg_{\Psi_m(x)}(\pi_m)\cdot\deg_x\Psi_m=2,$$
thus each multiple is either $1$ or $2$. But we know the total number
of points in $\L_m^K$ over $J=\infty$ (punctures) and this implies
the statement of Proposition~\ref{prop-ramJ} for $J=\infty$. 

That $\pi_m$ is an orbifold map follows from the identification
of the orbifold points in $\L_m^K$ in Proposition~\ref{special}.\hfill$\Box$
\vspace{.1in}

The difference between $\pi_m$ and  $\sigma_m$
is that there is no $\C^*$ action in the second case.

Proposition~\ref{prop-ramJ} together with relation (\ref{cd}) and known
ramification of $\psi$ allows us to define the orbifold structure
on the compactified spaces, so that the $\psi$ and $\Psi_m$ extend
to
orbifold coverings of these compactifications. 

For what follows we define compactifications of our orbifolds:
$$\bC_J=\bC(0(3),1(2),\infty(2)),\quad \bC_a=\bC.$$
Then $a\mapsto J=\psi(a)$ which is defined in (\ref{Ja}) is an orbifold
covering. Then we define ${\overline{\L}}_m^K$ by adding the punctures
$x, \;\pi_m^K(x)=\infty$, and defining $n(x)=1$ if $\deg_x(\pi_m^K)=2$
and $n(x)=2$ when $\deg_x(\pi_m^K)=1$.
Finally we define $\bC_a$ as the Riemann sphere with $n(a)=1$ for all $a$,
and define ${\overline{\H}}_m^j$ as $\H_m^j$ with filled punctures.
The orbifold structure on ${\overline{\H}}_m^j$ is trivial: $n(x)\equiv 1.$
With these definitions Theorem~\ref{theorem1} gives:
\begin{equation}\label{chiO}
\chi^O({\overline{\L}}_m^K)=\chi^O(\L_m^K)+d_m^K/2.
\end{equation}
\begin{prop}\label{total-ram} The following diagram is commutative:
\begin{equation}\label{cd1}
\begin{CD}
\overline{\H}_m^j    @>{\displaystyle\Psi_m^j}>>   \overline{\L}_m^K\\
@VV{\displaystyle\sigma_m} V           @VV{\displaystyle\pi_m} V\\
\bC_a      @>{\displaystyle\psi}>>         \bC_J
\end{CD}
\end{equation}
Here all four spaces are orbifolds, with orbifold functions just
defined,
the horizontal arrows
are orbifold coverings, and vertical arrows are maps of orbifolds.
We have
$$\deg\Psi_m^0=6\quad\mbox{and}\quad\deg\Psi_m^j=2,\quad j\in\{1,2,3\}.$$
Furthermore,
\begin{equation}\label{chiH0}
\chi({\overline{\H}}_m^0)=\chi^O({\overline{\H}}_m^0)=6\chi^O({\overline{\L}}_m^I),
\end{equation}
\begin{equation}\label{chiHj}
\chi(\overline{\H_m^j})=\chi^O({\overline{\H}}_m^j)=
2\chi^O({\overline{\L}}_m^{II}).
\end{equation}
\end{prop}

{\em Proof.} Since in the diagram (\ref{cd}) the horizontal arrows
are orbifold coverings and vertical arrows are orbifold maps,
it remains to check the points over $J=\infty$ and
over $a\in\{0,1,\infty\}$.
That $\psi:\bC_a\to\bC_J$ is an orbifold covering is well known
and follows from the explicit formula (\ref{Ja}).

Let $x\in\H_m^j$, $\sigma_m(x)\in\{0,1,\infty\}$.
By Proposition~\ref{prop-ramJ}, $\deg_x(\sigma_m)=1$ and we know that
$\deg_{\sigma(x)}(\psi)=2.$ Therefore, 
$$\deg_x(\Psi_m^j)\cdot\deg_{\Psi_m^j(x)}(\pi_m)=2,$$
thus $\deg_x\Psi_m^j$ is either $1$ or $2$, and the definition
of $n(\Psi_m^j(x))$ ensures that $\Psi_m^j$ is an orbifold
covering.  

That the vertical arrows are orbifold maps follows from
Proposition~\ref{prop-ramJ}. Formulas (\ref{chiH0}), (\ref{chiHj})
follow from (\ref{rh}).\hfill$\Box$

Now we are ready to prove Theorem~\ref{theorem2}.
\begin{prop}\label{generaH}
The Riemann surfaces $\bH_m^j$ are connected,
and their images in $\CP^2$ are non-singular.
\end{prop}

{\em Proof.} 
Using (\ref{rh}), (\ref{chiO}) and (\ref{chiH0}), we obtain
\begin{eqnarray*}
2-\chi({\overline{\H}}_m^0)&=&2-6\chi^O({\overline{\L}}_m^I)
=2-6\chi^O(\L_m^I)-3d_m^I\\
&=&2+(d_m^I)^2-3d_m^I
=(d_m^I-1)(d_m^I-2).
\end{eqnarray*}
Similarly, using (\ref{rh}), (\ref{chiO}) and (\ref{chiHj}),
we obtain
\begin{eqnarray*}
2-\chi({\overline{\H}}_m^j)&=&2-2\chi^O({\overline{\L}}_m^j)
=2-2\chi^O(\L_m^j)-d_m^{II}\\
&=&2+(d_m^{II})^2/9-d_m^{II}
=(d_m^{II}/3-1)(d_m^{II}/3-2),\quad j\in\{1,2,3\}.
\end{eqnarray*}
Therefore, in any case we have
\begin{equation}\label{aa}
2-\chi(\bH_m^j)=(\deg\bH_m^j-1)(\deg\bH_m^j-2).
\end{equation}
Suppose that for some $j$ and $m$,
$\bH_m^j$ has $N$ irreducible components of degrees $d_k$ genera $g_k$
and degrees $d_k$ for $1\leq k\leq N$.
Then
\begin{equation}\label{bb}
\deg\bH_m^j=\sum_{k=1}^Nd_k,\quad\chi(\bH_m^j)=\sum_{k=1}^N\chi_k,
\end{equation}
\begin{equation}\label{cc}
\chi_k=2-2g_k,\quad\mbox{and}\quad 2g_k\leq(d_k-1)(d_k-2);
\end{equation}
the last inequality follows from (\ref{gd}). Substituting the expressions
$\deg\bH_m^j$ and $\chi(\bH_m^j)$ from (\ref{bb}) to (\ref{aa})
and using (\ref{cc}) we obtain after simple manipulation
$$\left(\sum_{k=1}^Nd_k\right)^2\leq\sum_{k=1}^Nd_k^2;$$
since all $d_k\geq 1$, this is
possible only when $N=1$. Thus $\bH_m^j$ is irreducible.
Then from (\ref{aa}) we obtain its genus,
$$g(\bH_m^k)=(2-\chi(\bH_m^j))/2=(\deg\bH_m^j-1)(\deg\bH_m^j-2)/2,$$
so it is non-singular since it satisfies (\ref{gd})
with equality.\hfill $\Box$

This proposition 
completes
the proof of Theorem~\ref{theorem2}.
\vspace{.1in}

{\em Proof of Corollary \ref{corollary1}.}
\vspace{.1in}

Consider the map $R:\H_m\to\{ F_m(\lambda,g_2,g_3)=0\}$
defined in (\ref{R}). We will show that it is transversal
to the orbits
of the $\C^*$ action (\ref{scaling})
at non-special points. 

A trajectory of restriction of this action onto the $(g_2,g_3)$ plane
has the form $(g_2,g_3)=(t^2,ct^3),\, t\in\C^*$, so the tangent vector is
$(2t,3ct^2)$ which is parallel to $(2/g_3,3/g_2)=(2/R_3,3/R_2)$.
If the vectors $(2/R_3,3/R_2)$ and $(R_2^\prime,R_3^\prime)$ are
collinear, we must have
$$S:=R_2R_3(3R_2^\prime/R_2-2R_3^\prime/R_3)=0.$$
But an explicit computation shows that
$$S=-16a(a-1)/3,$$
which can be zero only at the special points.

Therefore the maps $\Psi_m^j$ are ramified only at the special points.
Diagram (\ref{cd1}) is clear from the definition. \hfill$\Box$
\vspace{.1in}

\noindent
{\em Proof of Corollary \ref{corollary2}.}
\vspace{.1in}

We compute the ramification of the forgetful map $\pi$ and then
make the correction for special points. 

The usual (not orbifold) Euler characteristic of the compactification
of $\L_m^{I}$ is $\chi(\bL_m^I)=\chi(\L_m^I)+h$, where $h$ equals
the number of punctures. So from Theorem~\ref{theorem1} for $\L_m^I$
and Riemann--Hurwitz formula
the total ramification of $\pi$ is 
\begin{eqnarray*}
&&2d-\chi(\bL_m^I)=2d-\chi(\L_m^I)-h=d-\chi(\L_m^I)\\
&=&d+d^2/6-(4\epsilon_0+3\epsilon_1)/6
=\lfloor (d^2-d+4)/6\rfloor+2\lfloor d/3\rfloor+\lfloor d/2\rfloor,
\end{eqnarray*}
where $d=\deg_\lambda F_m^{I}$. The first summand is the degree of
the Cohn polynomial, and the other two reflect the additional
ramification over the special points $0$ and $1$ (Proposition~\ref{prop-ramJ}).
Indeed, since the only singularities of
the surface $F_m^I(\lambda,g_2,g_3)=0$ lie over $g_2=0$ and $g_3=0$,
the zeros of the Cohn polynomial at all points $J\in\C\backslash\{0,1\}$
come from ramification points of $\pi$. For $J=0$, our curve
has the form $\lambda^kP(\lambda^3)$, where $P$ has only simple zeros,
so only $\lambda=0$ is a multiple zero when $k=2$. Other ramification
points of $\pi$ over $J=0$ do not contribute to zeros of Cohn's polynomial.

Similarly, $J=1$ is not a zero of Cohn's polynomial.

For $\L_m^{II}$ the total ramification is
$$d+d^2/18-(1-\epsilon_1)/2=(d/3)(d/3-1)/2+2d/3+\lfloor d/2\rfloor.$$
where $d=\deg_\lambda F_m^{II}.$ 
Again, the first summand corresponds to the degree of the Cohn polynomial
while the other two reflect additional ramification over the special points
(Proposition~\ref{prop-ramJ}).\hfill$\Box$
\vspace{.1in}

Next we briefly describe an alternative approach to our main results.
The following remarks are not necessary for understanding the
rest of the paper.
\vspace{.1in}

\noindent
{\bf Remarks on parametrization of $\H_m$ by a space of triangles}
\vspace{.1in}

In the beginning of the previous section we mentioned
that $\H_m$ represents the space of marked singular tori.

In view of the above interpretation of $\mathbf{H}_m$,
we can construct a natural lift of the isomorphism
$\Phi^*:\T_m^*\to\L_m$
to an isomorphism $\widehat{\Phi}^*:\widehat{\T}_m^*\rightarrow\H_m$ that makes
the following diagram commutative:
$$
\begin{CD}
\widehat{\T}_m^* @>{\displaystyle \Sigma} >> \T_m^*\\
@VV{\displaystyle\widehat{\Phi}^*} V @ VV{\displaystyle\Phi^*} V\\
\H_m @ >{\displaystyle\Psi} >> \L_m
\end{CD}
$$
Here $\widehat{\mathbf{T}}^*_m$ is a suitable space of flat singular triangles.
The point is that one could fully describe the
topology of the open Riemann surface $\widehat{\mathbf{T}}^*_m$
and so of $\mathbf{H}_m$ in a direct way,
and from that deduce the topology of $\mathbf{L}_m$
by taking the quotient by a natural $S_3$-action described below.

Recall that in a balanced flat singular triangle
$(D,\{a_i\},f)\in \T_m$,
the cyclic order of the corners $(a_1,a_2,a_3)$ on $\partial D$
matches the orientation induced by the developing map $f$,
and that we coordinatize
$\T_m$ using the functions $\phi_{i,j,k}$.
Denote by $-\T_m$ the space of balanced
flat singular triangles $(D,\{a_i\},f)$
for which the cyclic order $(a_1,a_3,a_2)$
matches the orientation induced by $f$, and coordinatize $-\T_m$
using the complex conjugates $\bar{\phi}_{i,j,k}$.
Moreover, let
$\widehat{\T}_m$ be the disjoint union of $\T_m$ and $-\T_m$.
The permutation group 
$S_3$ acts
on $\widehat{\T}_m$ by relabeling the corners.

By identifying each marginal triangle
$(D,\{a_i\},f)\in\widehat{\T}_m$
with its conjugate $(D,\{a_i\},\bar{f})$, we obtain a space $\widehat{\T}_m^*$.
It is immediate that the charts defined above
for $\widehat{\T}_m$ induce
a complex structure on $\widehat{\T}_m^*$, and that the $S_3$-action
descends to $\widehat{\T}_m^*$.
The quotient space $\widehat{\T}_m^*/S_3$ is naturally identified
with $\T_m^*$ and such identification induces the map $\Sigma$.

Moreover, to a triangle $T$ in $\widehat{\T}_m$
we can associate the torus $\Phi(T)$ with the marking $t(T)=(t_1,t_2,t_3)$,
where $t_1,t_2,t_3$ are the midpoints of
$[a_1,a_2]$, $[a_1,a_3]$, $[a_2,a_3]$ respectively.
Thus we can define $\widehat{\Phi}^*(T)=(\Phi(T),t(T))$.
Since the diagram is manifestly commutative, 
the map $\widehat{\Phi}^*:\widehat{\T}_m\to\H_m$
is an isomorphism of Riemann surfaces
(with trivial orbifold structure on both surfaces),
and $\Sigma$ is an orbifold cover.

Connected components of $\mathbf{H}_m$ can be
studied by analyzing $\widehat{\T}^*_m$, instead of exploring
the orbifold cover $\Psi_m$ as we did in the previous section.
Similarly to what we did with $\T_m$,
we can construct a nerve graph $\widehat{\Gamma}_m$
analogous to $\Gamma_m$. It consists of the disjoint union of
two isomorphic components: $\Gamma_m$ associated
to $\T_m$ and $-\Gamma_m$ associated to $-\T_m$.
Note that, if $v$ is a {\em lateral vertex} of $\Gamma_m$,
namely a vertex corresponding to a face of $\A_m$
adjacent to a boundary edge in $\partial\Delta_m$,  then
the preimage of $v$ and the preimage of $-v$ belong
to the same connected component
of $\widehat{\T}^*_m$.
Since each component of $\Gamma_m$ has a lateral
vertex, it follows that every component
of $\Gamma_m$ exactly corresponds to
a component in $\widehat{\T}^*_m$,
which thus has four connected components.

Similarly to what was done in Section~\ref{ident}
one can also compute the genera of $\H_m^j$
using the parametrization $\widehat{\Phi}:\widehat{\T}_m^*\to\H_m$
and obtain an alternative proof of our results about
$\H_m$ in this section. Once the number of components, their
genera and non-singularity are established for $\H_m$, one can
obtain Theorem ~\ref{theorem1} via (\ref{cd1}).

\section{Proof of Lemma~\ref{lemma2}}\label{tarasov}

In this section, we give a proof of Lemma~\ref{lemma2}
which is due to V. Tarasov.
It is inspired by an argument from \cite[Proposition 3]{ST}
which was brought
to our attention  by Eduardo Chavez Heredia from the University of Bristol.

The proof is a generalization of the classical arguments, going back
to Ch. Sturm, which are used in the proof of Lemma~\ref{lemma1}.

Let $D=\{ d_{i,j}\}$ be an $n\times n$ matrix
with entries
\begin{eqnarray*}
d_{i,i}=s,&&i=1,\dots, n,\\
d_{i,i+1}=a_i,&& i=1,\dots, n-1,\\
d_{i,i-2} = b_i,&& i=3,\dots, n.
\end{eqnarray*}
and consider the principal minors
$$D_k = \det(d_{ij})_{i,j=1,\dots, k}$$
They satisfy the recurrences
$$D_{k+3} = sD_{k+2} + c_k D_k,$$
where
$$c_k=a_{k+1}a_{k+2}b_{k+3}>0,$$
with the initial conditions
$$D_0=1,\quad D_1=s ,\quad D_2=s^2.$$
Then
$$D_{3j}(s)=P_j(s^3),\quad D_{3j+1}(s)=sQ_j(s^3),\quad D_{3j+2}(s)=
s^2R_j(s^3)$$
where the polynomials $P_j,Q_j,R_j$ satisfy
the recurrences
\begin{equation}\label{11t}
P_{j+1}=sR_j+A_jP_j,
\end{equation}
\begin{equation}\label{12}
Q_{j+1}=P_{j+1}+B_j Q_j,
\end{equation}
\begin{equation}\label{13}
R_{j+1}= Q_{j+1}+C_j R_j,
\end{equation}
where
$$A_j = c_{3j},\quad B_j=c_{3j+1},\quad C_j=c_{3j+2}$$
and the initial conditions
\begin{equation}\label{init}
P_0 = Q_0 = R_0 = 1.
\end{equation}
Lemma~9.2 follows from
\begin{prop}\label{prop-tarasov}
Let us define polynomials $P_j(s), Q_j(s), R_j(s)$ by recurrences
(\ref{11t}), (\ref{12}), (\ref{13}) where all
$A_j,B_j,C_j$ are strictly positive, and initial conditions (\ref{init}).
Then all polynomials $P_j, Q_j$ and $R_j$, $j=1,\dots$ are monic,
have degree $j$ and positive coefficients, and
all their roots are negative and simple. Moreover, if
$p_{j1}>\ldots>p_{jj}$, $q_{j1}>\ldots>q_{jj}$, $r_{j1}>\ldots>r_{jj}$,
are respective roots of the polynomials $P_j,Q_j,R_j$, then
$\;p_{jk}>q_{jk}>r_{jk}$ for all $\;k=1,\dots,j$, and $\;r_{jk}>p_{j,k+1}$
for all $\;k=1,\dots,j-1$.
\end{prop}

\noindent
{\em Proof.} We prove this by induction.
It is evident that our polynomials are monic and have positive coefficients.
Therefore their real roots are negative. To find the number of real
roots and to show that they are interlacent we look at the signs
of our polynomials at the roots of other polynomials using
(\ref{11t})--(\ref{13}).
\smallskip
The base of induction is given by $j=1$ and is clear.
The induction procedure is as follows.

\smallskip
By (\ref{11t}) and the induction assumption,
$$
P_{j+1}(0)\,P_{j+1}(p_{j1})\,=\,p_{j1} A_j P_j(0)\,R_j(p_{j1})<0\,,
$$
$$
P_{j+1}(r_{jk})\,P_{j+1}(p_{j,k+1})\,=\,
p_{j,k+1} A_j P_j(r_{jk})\,R_j(p_{j,k+1})<0\,,
$$
for all $k=1,\dots,j-1$, and
$$
(-1)^j P_{j+1}(r_{jj})=(-1)^j A_j P_j(r_{jj})>0\,.
$$
This implies that $P_{j+1}$ has roots $p_{j+1,k}$, $\,k=1,\dots, j+1$,
located as follows:
\begin{equation}
\label{pr}
0>p_{j+1,1}>p_{j1}\,,\quad
r_{j,k-1}>p_{j+1,k}>p_{jk}\,,\;\;k=2,\dots,j\,,\quad r_{jj}>p_{j+1,j+1}\,.
\end{equation}
Thus all roots of $P_{j+1}$ are negative and simple.

\smallskip
By the induction assumption and (\ref{pr}),
\begin{equation}
\label{pq}
p_{j+1,1}>q_{j1}\,,\quad
q_{j,k-1}>p_{j+1,k}>q_{jk}\,,\;\;k=2,\dots,j\,,\quad q_{jj}>p_{j+1,j+1}\,.
\end{equation}
Then by (\ref{12}) and (\ref{pq}),
$$
Q_{j+1}(p_{j+1,k})\,Q_{j+1}(q_{jk})\,=\,B_j Q_j(p_{j+1,k})\,P_{j+1}(q_{jk})<0\,,
$$
for all $k=1,\dots,j$, and
$$
(-1)^j Q_{j+1}(p_{j+1,j+1})=(-1)^j B_j Q_j(p_{j+1,j+1})>0\,.
$$
This implies that $Q_{j+1}$ has roots $q_{j+1,k}$, $\,k=1,\dots, j+1$,
located as follows:
\begin{equation}
\label{qp}
p_{j+1,k}>q_{j+1,k}>q_{jk}\,,\;\;k=1,\dots,j\,,\quad
p_{j+1,j+1}>q_{j+1,j+1}\,.
\end{equation}
Thus all roots of $Q_{j+1}$ are negative and simple.

\smallskip
By the induction assumption, (\ref{qp}), and (\ref{pr}),
\begin{equation}
\label{qr}
q_{j+1,1}>r_{j1}\,,\quad
r_{j,k-1}>q_{j+1,k}>r_{jk}\,,\;\;k=2,\dots,j\,,\quad r_{jj}>q_{j+1,j+1}\,.
\end{equation}
Then by (\ref{13}) and (\ref{qr}),
$$
R_{j+1}(q_{j+1,k})\,R_{j+1}(r_{jk})\,=\,C_j R_j(q_{j+1,k})\,Q_{j+1}(r_{jk})<0\,,
$$
for all $k=1,\dots,j$, and
$$
(-1)^j R_{j+1}(q_{j+1,j+1})=(-1)^j C_j R_j(q_{j+1,j+1})>0\,.
$$
This implies that $R_{j+1}$ has roots $r_{j+1,k}$, $\,k=1,\dots, j+1$,
located as follows:
\begin{equation}
\label{rq}
q_{j+1,k}>r_{j+1,k}>r_{jk}\,,\;\;k=1,\dots,j\,,\quad
q_{j+1,j+1}>r_{j+1,j+1}\,.
\end{equation}
Thus all roots of $R_{j+1}$ are negative and simple.

\smallskip
The inequalities $\;p_{j+1,k}>q_{j+1,k}>r_{j+1,k}$ for all $\;k=1,\dots,j+1$,
follow from (\ref{qp}) and (\ref{rq}), and the inequalities
$\;r_{j+1,k}>p_{j+1,k+1}$ for all $\;k=1,\dots,j$, follow from
(\ref{rq}) and (\ref{pr}). This completes the proof.\hfill$\Box$

\section{Projective monodromy of the Lam\'e equation}\label{mono}

In this section we prove Theorem~\ref{LWm}.

Lam\'e equation with integer $m$ has trivial local monodromy about the origin.
Since the fundamental group of the torus is $\Z^2$, the projective
monodromy
is represented by a pair of commuting elements $PSL(2,\C)$.

So we investigate the set of pairs of commuting elements of $PSL(2,\C)$
modulo conjugation. Every such pair $(A,B)$ can be
conjugated to one of the following forms:
\begin{equation}
\label{diag}
(z\mapsto \mu_1z,\;z\mapsto\mu_2z),\quad(\mu_1,\mu_2)\in (\C^*)^2,
\end{equation}
\begin{equation}\label{transl}
(z\mapsto z+a_1,\; z\mapsto z+a_2),\quad(a_1,a_2)\in\C^2,
\end{equation}
or
\begin{equation}\label{klein}
(z\mapsto -z,\; z\mapsto 1/z).
\end{equation}
It is proved in \cite[Theorem 2.2]{L1} that the third possibility (\ref{klein})
cannot happen for
the projective monodromy of Lam\'e equations with integer $m$.
Notice that $PSL(2,\C)$ representations
(\ref{diag}) and (\ref{transl}) can be lifted to $SL(2,\C)$. The
pair $(A,B)=(\id,\id)$ is also excluded.

Conjugacy classes of pairs of the form (\ref{diag}) are parametrized by
$(\C^*)^2$. Two pairs of the form (\ref{transl}) are conjugate iff 
$(a_1:a_2)=(a_1^\prime:a_2^\prime)$, so they are parametrized
by the projective line $\CP^1$ and one point $(0,0)$.

Suppose that a sequence of pairs $(A,B)$ of type (\ref{diag}) converges
in $PSL(2,\C)$
to a pair of type (\ref{transl}).
To figure out how $(a_1,a_2)$ in (\ref{transl}) are related to $(\mu_1,\mu_2)$
we consider
commuting pairs of linear-fractional transformations $(\phi_1,\phi_2)$,
$$\phi_j(z)=\frac{(1+f_j)z+(a_j+p_j)}{q_jz+(1+g_j)},\quad j\in\{1,2\},$$
where $f_j,g_j,p_j,q_j$ are small numbers,
and $a_j$ are constants not simultaneously equal to $0$.

The condition that matrices representing $\phi_j$ have determinant $1$ gives
\begin{equation}\label{11}
f_j+g_j\equiv a_jq_j,
\end{equation}
where $\equiv$ means that we neglected the terms of
order $2$ and higher.
The condition that $\phi_1$ and $\phi_2$ commute implies by comparing
the diagonal elements of the product matrices
\begin{equation}\label{comm}
a_1q_2\equiv a_2q_1.
\end{equation}
Now it follows from (\ref{11}) and (\ref{comm}) that
\begin{equation}
\label{cond1}
(f_1+g_1):(f_2+g_2)\to (a_1^2:a_2^2).
\end{equation}
Similar equation can be obtained when the eigenvalue of
one or both limit matrices is $-1$.
In other words, if $A$ and $B$ are $SL(2,\C)$ matrices
representing $\phi_1$ and $\phi_2$ tending to parabolic $\phi_1^*$
and $\phi_2^*$, then we have
\begin{equation}\label{cond2}
\lim\frac{\tr^2 A-4}{\tr^2 B-4}=(a_1^2:a_2^2),
\end{equation}
where $(a_1:a_2)$ is the ``ratio invariant'' of the pair of commuting
parabolic transformations.

So we obtain that the space of projective monodromy representations
for Lam\'e equations is the blow up of $(\C^*)^2$ at the point $(1,1)$.
Monodromies of unitarizable Lam\'e equations form the real torus
$$\{(\mu_1,\mu_2)\in(\C^*)^2:|\mu_1|=|\mu_2|=1\}\backslash(1,1),$$
and the boundary of this real torus in the blow up is the real projective line.
Since $A$ and $B$ are elliptic,
the left hand side of (\ref{cond2}) is positive, so the ratio
$(a_1:a_2)$ is real.

Thus we obtain
\begin{prop}\label{real}
Abelian integrals
arising as limits of developing maps of spherical metrics
have real period ratios.
\end{prop}

So we have
$$\partial\,\Sph_{1,1}(2m+1)\subset\LW_m.$$
Let us prove the opposite inclusion.
We recall that $\Lame_m$ is biholomorphic to $\CP(1,2,3)$,
and each point of $\Lame_m$ corresponds to a Schwarz equation (\ref{schwa})
whose solution is a developing map of a translation structure
on a torus with one conic singularity with angle $2\pi(2m+1)$.
Let us fix $L\in\LW_m$. We want to prove that there exists
a sequence $L_n\in\Lame_m$, $L_n\to L$ such that $L_n$
have unitarizable monodromy.

Since $L\in\LW_m$, is obtained as $\Phi(T)$ from some BFT,
we can normalize $T$ so that
the images of corners under the developing map are
$$(-\delta,k\delta,\delta),$$
where $k\in(-1,1)$ depends on the similarity class of $T$ and thus
is fixed, while $\delta>0$ is arbitrary.
We glue a congruent copy $-T$ to $T$ by the side
with endpoints $\{-\delta,\delta\}$. Notice that the gluing map $z\mapsto-z$
is the orientation-reversing isometry of this side, with respect to
both flat and spherical metric. The result of this gluing is a
``parallelogram'' $Q$, which can be represented by
$(\Delta,\{ a_j\},f)$, where $\Delta$ in the unit disk, $a_1,a_2,a_3,a_4$
are boundary points (corners), and $f$ is a meromorphic function,
a developing map of $Q$. This developing map is a local homeomorphism
at all points of $\overline\Delta\backslash\{ a_j\}$ while at the corners
it behaves like an integer power. The torus $\Phi(T)$
is obtained by gluing the
opposite sides of $Q$, so we have
\begin{equation}\label{71}
f(p(z))=t_i\circ f(z),\quad z\in\partial\Delta, \quad i=1,\ldots,4,
\end{equation}
where $p(z)$ is an involution of $\partial\Delta$ defining the gluing
of opposite sides, $t_j$ are translations, and $i$ depends
on the side to which $z$ belongs. Subscript $i$ has to
be interpreted similarly in all following formulas (\ref{72}), (\ref{75}), (\ref{73}).
Each of these translations $t_i$
is a product of two Euclidean rotations by $\pi$ about the ``middles'' of
the sides (see Proposition~\ref{prop2.5}).

Notice that $T$ can be also considered as a spherical triangle, since its
sides belong to the real line, which is both spherical and Euclidean geodesic.
So we can produce from $T$ a spherical torus with one conic singularity
with the same angle as in $\Phi(T)$, by identifying the opposite sides
of $Q$ by spherical isometries $s_i,\; 1\leq i\leq 4$. Similarly to $t_i$,
each $s_i$ is a product of two spherical rotations by $\pi$ about
the ``spherical middles'' of the sides.

Notice the following: {\em let $a,b$ be two points in $[-\delta,\delta]$.
There exists a unique Euclidean rotation $e$ by angle $\pi$,
and a unique spherical
rotation $s$ by angle $\pi$ which interchange $a$ and $b$.
Moreover, when $\delta$ is
small, then $e$ and $s$ are close to each other
in the sense that $s=e\circ\phi$, where
$\phi$ is a linear fractional homeomorphism of $R$ which
is $\epsilon$-bi-Lipschitz with respect to the spherical metric, where
 $\epsilon=\epsilon(\delta)\to 1$ as $\delta\to 0$.}
We call such homeomorphisms $\phi$
{\em close to the identity.} Using this statement,
we obtain that
\begin{equation}\label{72}
s_i=t_i\circ\phi_i,
\end{equation}
where $\phi_j$ are close to the identity. Now we define homeomorphisms
$\psi_i$ of the sides of $\Delta$ and the involution
$q$ of the unit circle $\partial\Delta$ by
\begin{equation}\label{75}
f\circ\psi_i=\phi_i^{-1}\circ f,\quad q:=\psi_i^{-1}\circ p,
\end{equation}
and the function $g:\partial\Delta\to\R$,
\begin{equation}\label{73}
g=f\circ\psi_i.
\end{equation}
It is easy to see that $\psi_i$ are close to the
identity: the only points where $f$ is not bi-Lipschitz are
the corners, where $f$ behaves like $z^n$, and conjugating
a diffeomorphism which is $C^1$-close to the identity
by $z^m$ result in a diffeomorphism which is $C^1$ close to the identity.
 Now we use the theorem from \cite[Thms. 5.3, 5.17]{Tu} which says that
every $\epsilon$-bi-Lipschitz homeomorphism of the unit circle
with $\epsilon$ sufficiently close to $1$ has an $\epsilon_1$-bi-Lipschitz
extension to a homeomorphism of the unit disk. So we extend $\psi$ and
define $g$ in $\Delta$ by (\ref{73}). This extended $g$
satisfies
$$g(q(z))=s_i\circ g(z),\quad z\in\partial\Delta,$$
which is analogous to (\ref{71}), and follows from (\ref{71}), (\ref{72}),
(\ref{75}) and (\ref{73}).
So $g$ can be considered
as a developing
map of a spherical metric on a torus obtained by identifying the opposite
sides
by the involution $q$. The map $g$ is not conformal with respect to
the standard conformal structure on $\Delta$ but conformal in a new
structure defined by $\psi$. This new structure is close to the standard
one in the sense of the theory of quasiconformal mappings \cite{Ahlfors}.
Indeed, a bi-Lipschitz map with constant close to $1$ is quasiconformal
with dilatation close to $1$. Let $L_\delta$ be the spherical torus
obtained by this construction, and let $F$ and $G_\delta$ be the
lifts the developing maps
of the tori $\Phi(T)$ and $L_\delta$ to the universal cover.
Then $F$ is a meromorphic function, and $G_\delta$
a quasiregular one, conformal
with respect to the non-standard conformal structure in $\C$.
We normalize $F$ and $G_\delta$ so that the preimages
of the conic singularities
form lattices $\Lambda_F$ and $\Lambda_{G_\delta}$,
with primitive periods $1$ and $\omega_F,\omega_{G_\delta}$.

To compare $F$ and $G_\delta$ we find a homeomorphism $\eta_\delta:\C\to\C$
such that
$F_\delta:=G_\delta\circ\eta_\delta$ is a meromorphic function,
and $\eta_\delta$ is normalized by $\eta_\delta(0)=0,\;
\eta_\delta(1)=1$. This homeomorphism is obtained as a normalized solution
of the Beltrami equation \cite[Chap. V]{Ahlfors},
and since the quasiconformal dilatation of $\eta_\delta$ is small,
$\eta_\delta\to\mathrm{id},$ as $\delta\to 0$.
We conclude that the lattices of poles $\Lambda_{F_\delta}$
and $\Lambda_F$ are close to each other 
uniformly on compacts, 
and $F_\delta\to F$ as $\delta\to 0$ uniformly
(with respect to the spherical
metric) on compact subsets of the plane. This means that when $\delta\to 0$,
the lattice invariants $g_j(\Lambda_{F_\delta})\to g_j(\Lambda_F),\;
j=2,3,$ and the Schwarzian of $F_\delta$ converges
to the Schwarzian of $F$. Since the Schwarzians
of $F$ and $F_\delta$ are both of the form
$$m(m+1)\wp(z,\Lambda)+\lambda,$$
where $\lambda=\lambda(F)$ or $\lambda=\lambda(F_\delta)$,
we conclude that $\lambda(F_\delta)\to\lambda(F)$ as $\delta\to 0$,
and this proves
that $L_\delta\to L$ in $\Lame_m$. This completes the proof of 
Theorem~\ref{LWm}.

\section{Lin--Wang curves}\label{linwang}

Proposition~\ref{real} gives
the following
characterization of Lin--Wang curves: {\em $J\in\pi_m(\LW_m)$ if
some Lam\'e equation with invariant $J$ has translational monodromy
which belongs to a straight line.}
In particular, there is a flat singular
torus corresponding to $(J,\lambda)$, and this torus is of the
form $\Phi^*(T)$ where $T$ is a BFT with integer angles.
\vspace{.1in}

{\em Proof of Theorem~\ref{theorem4}.}
\vspace{.1in}

The set $\LW_m$ is the image
of the set of triangles with angles integer multiples of $\pi$ under
the map $\Phi_m^*$. Triangles with angles integer multiples of $\pi$
form straight intervals in the local coordinates $\phi_{i,j,k}$ introduced
in Section~\ref{complexanalytic}, and the map $\Phi^*_m$ is biholomorphic.

We conclude that the set $\LW_m$ consists of analytic (non-singular)
curves. There are three such curves corresponding to each integer point
in $\A_m$. So all together we have $m(m+1)/2$ Lin--Wang curves.
\hfill$\Box$
\vspace{.1in}

Proposition~\ref{real} shows how to find
equations of Lin--Wang curves. They are curves in $\L_m$, the
moduli space of Lam\'e functions,
where the ratio of periods of the integral (\ref{abel}) is real.
According to \cite{EMP}, to each triple of integers satisfying
the triangle inequalities corresponds a component of the space
$\Sph_{1,1}(2m+1)$. This component is parametrized by an open triangle,
and the sides of this triangle correspond to three interior edges of $\T_m$
which are mapped by $\phi$ in (\ref{phi})
to integer points of $\A_m$ by the map $\phi$
in (\ref{phi}). These three edges parametrize 
Lin--Wang curves by the map $\pi_m\circ\Phi_m^*$.

Thus each integer point in $\A_m$ corresponds to a component
of the moduli space of $\Sph_{1,1}(2m+1)$ of spherical tori. The boundary of
this component consists of one or three Lin--Wang curves: when the integer
point is the center of $\A_m$, which happens when $m\equiv1\, (\mod 3)$,
there is one curve, otherwise there are three of them. If there is one
curve, it belongs to $\L_m^{II}$, as it happens for $m=1$ and $m=4$.

For integer points other than the center of $\A_m$ we have three curves which can be all
on component $\L_m^{II}$, or one of them can be on $\L_m^I$ and two on
$\L_m^{II}$.
Figs.~6, 7 and similar figures for other $m$ permits to determine this
for every integer point on $\A_m$.

We give explicit formulas for Lin--Wang curves for $1\leq m\leq 3$.
We recall that the developing
map is given by (\ref{Ab}). Since working with
elliptic functions
is easier than with Abelian integrals (this was the primary reason
for introducing elliptic functions), we pass to the universal covering. 

In what follows, $g$ is the integrand in (\ref{Ab}). It is an even
elliptic function
with a single zero of multiplicity $2m$ at the origin, and double poles
with vanishing residues. The general form of such function is
$$g(z)=c_0+\sum_{j=1}^mc_j\wp(z-a_j).$$
By Abel's theorem, $2(a_1+\ldots+a_m)\equiv 0$, and we want to choose
$c_j$ so that $g$ and its first $2m-1$ derivatives vanish at the origin.
So
$$c_0=-\sum_{j=1}^mc_j\wp(a_j),$$
and for the rest of $c_j$ we have a system of equations
$$\sum_{j=1}^mc_j\wp^{(k)}(-a_j)=0,\quad 1\leq k\leq 2m-1,$$
which has a non-trivial solution if the matrix of this system has rank
at most $m-1$.

Once $g$ is found, we are interested in the ratio of periods
of the integral (\ref{Ab}). Lin--Wang curves make the locus of points where
this ratio is real. Below we give the results of computation for $1\leq m\leq3$.
We use the standard notation of the theory of elliptic functions \cite{A}:
$\wp(\omega_j)=e_j,\; 1\leq j\leq 3$, so $\omega_j$ are {\em half-periods},
$\zeta$ is the Weierstrass zeta function, $\zeta'=-\wp$, satisfying
$$\zeta(z+2\omega_j)=\zeta(z)+2\eta_j.$$
\vspace{.1in}

\noindent
{\em Case $m=1$}. $g(z)=\wp(z-\omega_j)-e_j,\; j\in\{1,2,3\}$. In this case,
$\L_1=\L_1^{II}$. 
 The equation of Lin--Wang curves is
$$\Ima\frac{\eta_1+\omega_1e_j}{\eta_2+\omega_2e_j}=0,\quad j\in\{1,2,3\}.$$
This defines three curves in the fundamental region of the modular group
in the $\tau$-half-plane, Fig.~12,
which correspond to one curve in the $J$-plane
(Fig.~13).
\vspace{.1in}

\noindent
{\em Case $m=2$}. For $\L_2^I$, we obtain
$$g(z)=\wp(z+a)+\wp(z-a)-2\wp(a),\quad\mbox{where}\quad \wp(a)=\sqrt{g_2/12}.$$
The equation of the Lin--Wang curve is
$$\Ima\frac{\eta_1+\omega_1\sqrt{g_2/12}}{\eta_2+\omega_2\sqrt{g_2/12}}=0.$$

For $\L_2^{II}$, we have
$$g(z)=\wp^{\prime\prime}(\omega_j)\left(\wp(z+\omega_k)-e_k\right)-
\wp^{\prime\prime}(\omega_k)\left(\wp(z+\omega_j)-e_j\right),$$
and the equation of Lin--Wang curve is
$$\Ima\frac{\omega_1(e_j\wp''(\omega_k)-e_k\wp''(\omega_j))+
\eta_1(\wp''(\omega_k)-\wp''(\omega_j))}{\omega_2(e_j\wp''(\omega_k)-e_k\wp''(\omega_j))+\eta_2(\wp''(\omega_k)-\wp''(\omega_j))}=0.$$
These curves in the $\tau$-half-plane are shown in Fig.~14 and their images
in the $J$-plane are in Fig.~15. These are the three curves bounding a single
triangle
which is the moduli space $\Sph_{1,1}(5)$. 
One of these curves, (which has a loop
in Fig. 12) belongs to $\L_2^I$, other two belong to $\L_2^{II}$.
Shading in Fig.~15 is the {\em hypothetical} projection
of a component of $\Sph_{1,1}(5)$ by the forgetful map.
We do not know whether the restriction of the
forgetful map on $\Sph_{1,1}(2m+1)$
is open. So we don't know that the boundary of this projection
is contained in Lin--Wang curves.
\vspace{.1in}

\noindent
{\em Case $m=3$}. For $\L_3^I$ we obtain:
$$g(z)=c_0+c_1\wp(z+\omega_1)+c_2\wp(z+\omega_2)+c_3\wp(z+\omega_3),$$
where
$$c_0=-\sum_{j=1}^3c_je_j,$$
$$c_k=(6e_{k+2}^2-g_2/2)(6e_{k+1}^2-g_2/2)(e_{k+2}-e_{k+1}),\quad k\in\{1,2,3\},$$
where the subscripts are understood as residues $\mod 3$, but we use $3$
instead of $0$ to prevent the confusion with previous formula.
Setting $B=c_1+c_2+c_3$, the equation of the Lin-Wang curve is
$$\Ima\frac{\omega_1c_0-\eta_1B}{\omega_2c_0-\eta_2B}=0.$$
This curve in the $\tau$-plane is shown 
in Fig.~16.

For component $\L_3^{II}$ we introduce the notation for $k\in\{1,2,3\}$:
$$P_k^\pm=-e_k/5\pm\sqrt{3(5g_2/4-3e_k^2)},$$
$$c_{1,k}^\pm=-2\frac{6(P_k^\pm)^2-g_2/2}{6e_k^2-g_2/2},$$
$$c_{0,k}^\pm=-c_{1,k}^\pm e_k-2P_k^\pm.$$
Then
$$g(z)=c_{0,k}^\pm+c_{1,k}^\pm\wp(z+\omega_k)+\wp(z+a_k^\pm)+\wp(z-a_k^\pm),$$
where $\pm a_k^\pm$ are solutions of the equation $\wp(z)=P_k^\pm$.
Then we have six Lin--Wang curves
$$\Ima\frac{c_{0,k}^\pm\omega_1-(c_{1,k}^\pm+2)\eta_1}{c_{0,k}^\pm\omega_2-
(c_{1,k}^\pm+2)\eta_2}=0.$$
These curves in the $\tau$-half-plane are shown in Fig. 17, and their image
in the $J$-plane in Fig. 18, where the curve from $\L_3^I$
is also included (it is the one which looks like a vertical line in
the middle). Fig. 18 shows in the right-hand side
the detail which looks like
a tripod in the left-hand side. Figure 17 contains 16 curves which give $5$
images in Fig. 18. Three of these $5$ curves in Fig. 18
constitute the full boundary
of one triangle of the moduli space for spherical tori, and the remaining
three curves in Fig. 18, including that one curve which comes from Fig. 16
constitute the boundary of the second triangle in the moduli space for 
spherical tori.
\newpage

\noindent
\section*{Appendix. List of formulas}
\vspace{.2in}

Polynomials $F_m^I$:
\vspace{.2in}

\begin{tabular}{r|l}
2&$ \lambda^2-3g_2$\\ \\
3&$ \lambda$\\ \\
4&$ \lambda^3-52g_2\lambda+560g_3$\\ \\
5&$ \lambda^2-27g_2$\\ \\
6&$ \lambda^4-294g_2\lambda^2+7776g_3\lambda+3465 g_2^2$\\ \\
7&$ \lambda^3-196g_2\lambda+2288g_3$\\ \\
8&$ \lambda^5-1044g_2\lambda^3+48816g_3\lambda^2+112320g_2^2\lambda
-4665600g_2g_3$\\ \\
9&$ \lambda^4-774g_2\lambda^2+21600g_3\lambda+41769 g_2^2.$
\end{tabular}
\vspace{.2in}

Polynomials $F_m^{II}$:
\vspace{.2in}

\begin{tabular}{r|l}
1&$ 4\lambda^3-g_2\lambda-g_3$\\ \\
2&$ 4\lambda^2-9g_2\lambda+27 g_3$\\ \\
3&$ 16\lambda^6-504g_2\lambda^4+2376g_3\lambda^3+4185g_2^2\lambda^2-
36450g_2g_3\lambda-3375g_2^3+91125g_3^2$\\ \\
4&$ 16\lambda^6-1016g_2\lambda^4+8200g_3\lambda^3+10297g_2^2\lambda^2
-41650g_2g_3\lambda$\\
&$-27783g_2^3-42875g_3^2.$
\end{tabular}
\newpage

Polynomials $H_m^0$:
\vspace{.2in}

\begin{tabular}{r|l}
2&$ B^2+4(a+1)B+12a$\\ \\
3&$ B+4(a+1)$\\ \\
4&$ B^3+20(a+1)B^2+(64a^2+336a+64)B+640(a^2+a)$\\ \\
5&$ B^2+20(a+1)B+64(a^2+1)$\\ \\
6&$ B^4+56(a+1)B^3+(784a^2+2744a+784)B^2$\\
&$+(2304a^3+29472a^2+29472a+2304)B$\\
&$+48384a^3+152208a^2+48384a$
\end{tabular}
\vspace{.2in}

Polynomials $H_m^1$:
\vspace{.2in}

\begin{tabular}{r|l}
1&$ B+a+1$\\ \\
2&$ B+4a+1$\\ \\
3&$ B^2+10(a+1)B+9a^2+78a+9$\\ \\
4&$ B^2+(20a+10)B+64a^2+136a+9$\\ \\
5&$ B^3+35(a+1)B^2+(259a^2+1046a+259)B$\\
&$+225a^3+5235(a^2+a)+225$\\ \\
6&$ B^3+(56a+35)B^2+(784a^2+1568a+259)B$\\
&$+2304a^3+13008a^2+7464a+225$
\end{tabular}
\vspace{.2in}

Degrees of forgetful maps:

\begin{eqnarray}
d_m^I&:=&\left\{\begin{array}{ll}m/2+1,& m\equiv0\;
(\mod 2)\\
(m-1)/2,& m\equiv 1\; (\mod 2),\end{array}\right.\\
d_m^{II}&:=&3\lceil m/2\rceil.
\end{eqnarray}
\newpage

Euler characteristics:

$$\chi(\L_m^I)=\left\{\begin{array}{ll}\displaystyle
-\frac{(m+2)^2}{24}+\frac{4\epsilon_0+3\epsilon_1}{6},&m\equiv 0\;(\mod 2),\\ \\
\displaystyle
-\frac{(m-1)^2}{24}+\frac{4\epsilon_0+3\epsilon_1}{6},&m\equiv 1\;(\mod 2).
\end{array}\right.$$
$$\chi(\L_m^{II})=\left\{\begin{array}{ll}\displaystyle
-\frac{m^2}{8}+\frac{1-\epsilon_1}{2},& m\equiv0\;(\mod 2),\\ \\
\displaystyle
-\frac{(m+1)^2}{8}+\frac{1-\epsilon_1}{2},& m\equiv 1\;(\mod 2).
\end{array}\right.,$$
where
\begin{eqnarray*}
\epsilon_0&=&\left\{\begin{array}{ll}0,&\mbox{if}\; m\equiv1\;(\mod 3),\\
1&\mbox{otherwise,}\end{array}\right.\\ \\
\epsilon_1&=&\left\{\begin{array}{ll}0,&\mbox{if}\; m\in\{1,2\}\;(\mod 4),\\
1&\mbox{otherwise.}\end{array}\right.
\end{eqnarray*}

Numbers of punctures: $h_m^I=d^I_m,\; h_m^{II}=2d^{II}_m/3.$

Genera in terms of $d_m^K$:
\begin{equation}\label{gI}
g(\L^I_m)=1+\frac{(d_m^I)^2}{12}-\frac{d_m^I}{2}-
\frac{4\epsilon_0+3\epsilon_1}{12},\end{equation}
\begin{equation}\label{gII}
g(\L^{II}_m)=1+\frac{(d_m^{II})^2}{36}-\frac{d_m^{II}}{3}-
\frac{1-\epsilon_1}{4}.
\end{equation}

Genera in terms of $m$:
\begin{equation}\label{gIe}
g(\L_m^I)=\frac{m^2-8m+28}{48}-\frac{4\epsilon_0+3\epsilon_1}{12},
\quad m\equiv0\; (\mod 2),
\end{equation}
\begin{equation}\label{gIo}
g(\L_m^I)=\frac{m^2-14m+61}{48}-\frac{4\epsilon_0+3\epsilon_1}{12},
\quad m\equiv1\; (\mod 2),
\end{equation}
\begin{equation}\label{gIIe}
g(\L_m^{II})=\frac{m^2-8m+16}{16}-\frac{1-\epsilon_1}{4},\quad
m\equiv0\; (\mod 2),
\end{equation}
\begin{equation}\label{gIIo}
g(\L_m^{II})=\frac{m^2-6m+9}{16}-\frac{1-\epsilon_1}{4},
\quad m\equiv 1\; (\mod 2).
\end{equation}

Degrees of Cohn's polynomials
$$\deg C_m^I=\left\lfloor \frac{(d_m^I)^2-d_m^I+4)}{6}\right\rfloor,$$
$$\deg C_m^{II}=\frac{d_m^{II}(d_m^{II}-1)}{2}.$$

\vspace{.2in}

{\em A. E. and A. G.: Purdue University, West Lafayette, IN 47907 USA

eremenko@purdue.edu, agabriel@math.purdue.edu
\vspace{.1in}

G. M.: ``Sapienza'' Universit\`a di Roma, Rome, Italy,
mondello@mat.uniroma1.it
\vspace{.1in}

D. P.: Kings College, London UK, dmitri.panov@kcl.ac.uk}
\begin{figure}\label{fig1}
\begin{center}
\includegraphics*[width=5.0in]{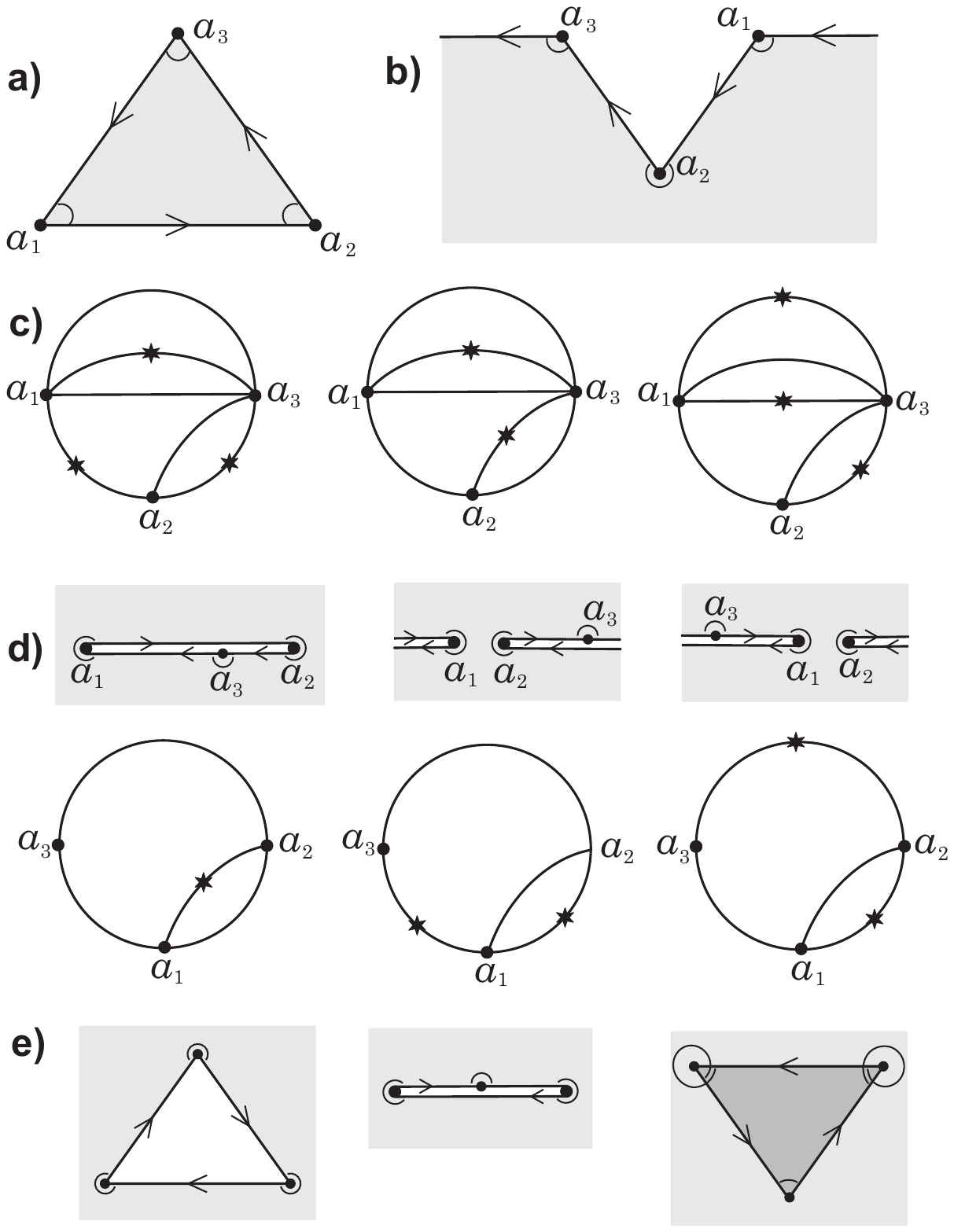}
\caption{a),b) - primitive triangles; c) - examples of nets with the angles
$(3\pi,2\pi,4\pi)$; d) - triangles
with the angles $(2\pi,\pi,2\pi)$ and their nets; e) - deformation of triangle
with the angle sum $5\pi$.}
\end{center}
\end{figure}
\begin{figure}\label{fig2}
\begin{center}
\includegraphics*[width=4.8in]{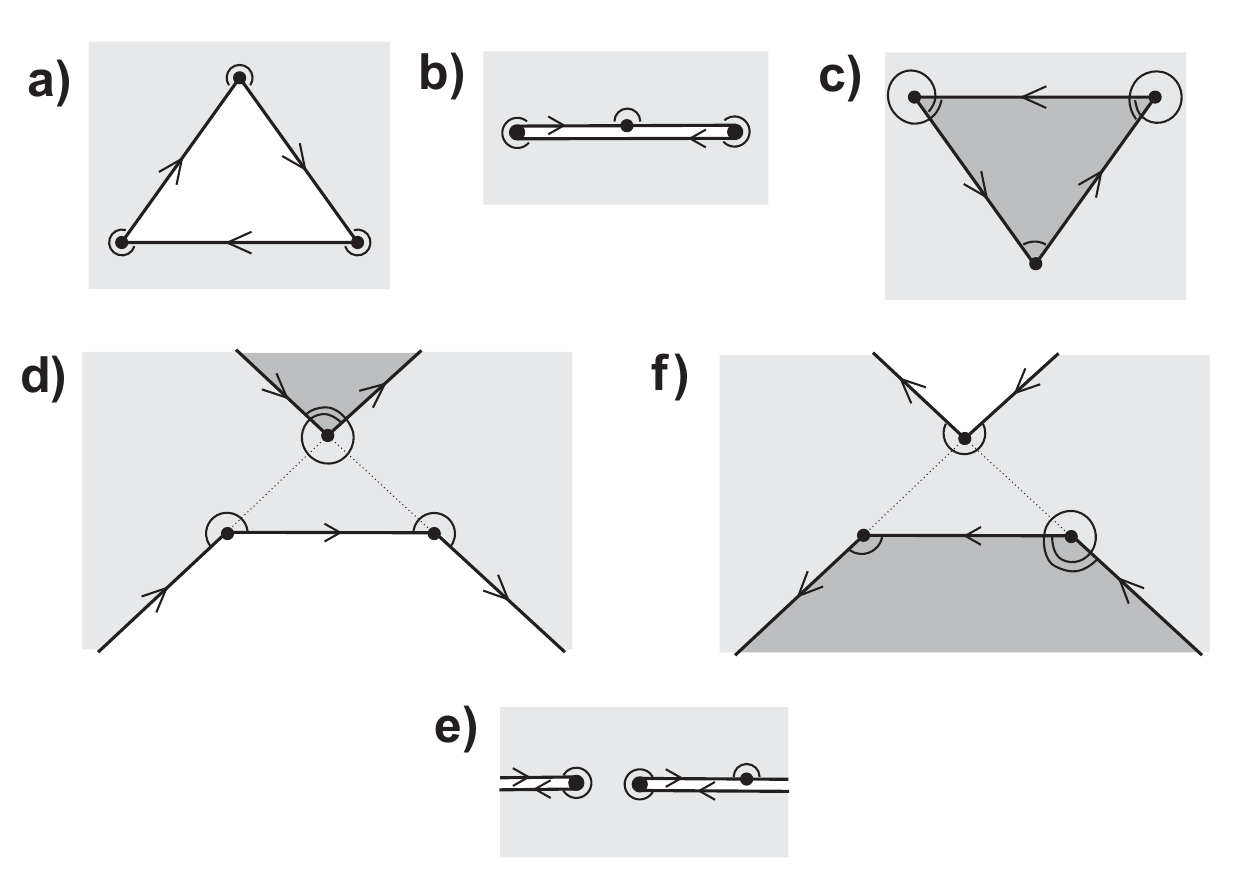}
\caption{All types of BFT for $m=2$.}
\end{center}
\end{figure}
\begin{figure}\label{fig3}
\begin{center}
\includegraphics*[width=4.8in]{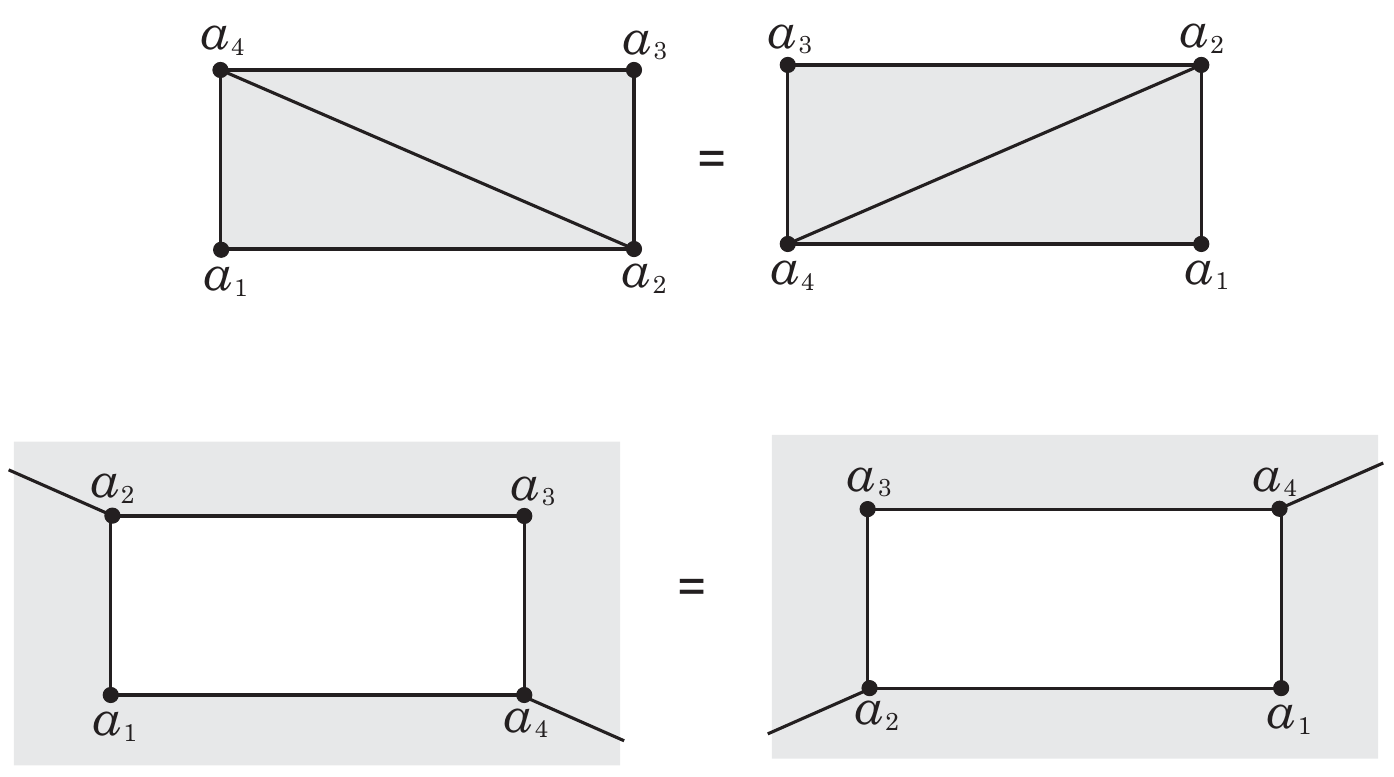}
\caption{To the proof of Proposition \ref{prop3}.}
\end{center}
\end{figure}
\begin{figure}\label{fig4}
\begin{center}
\includegraphics*[width=5in]{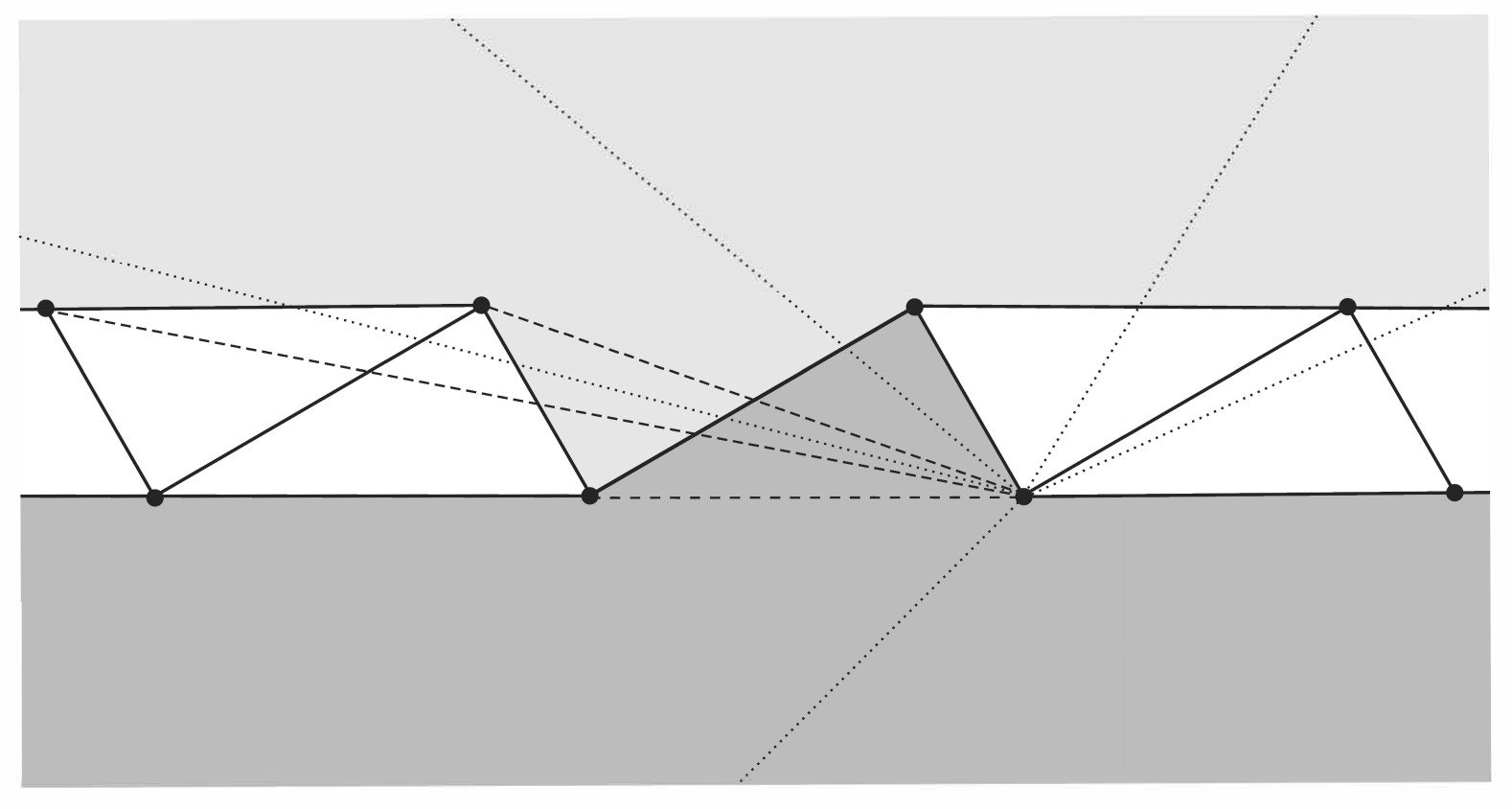}
\caption{To the proof of Lemma \ref{u1} for a triangle of type $A''$.}
\end{center}
\end{figure}

\begin{figure}\label{fig5}
\begin{center}
\includegraphics*[width=5in]{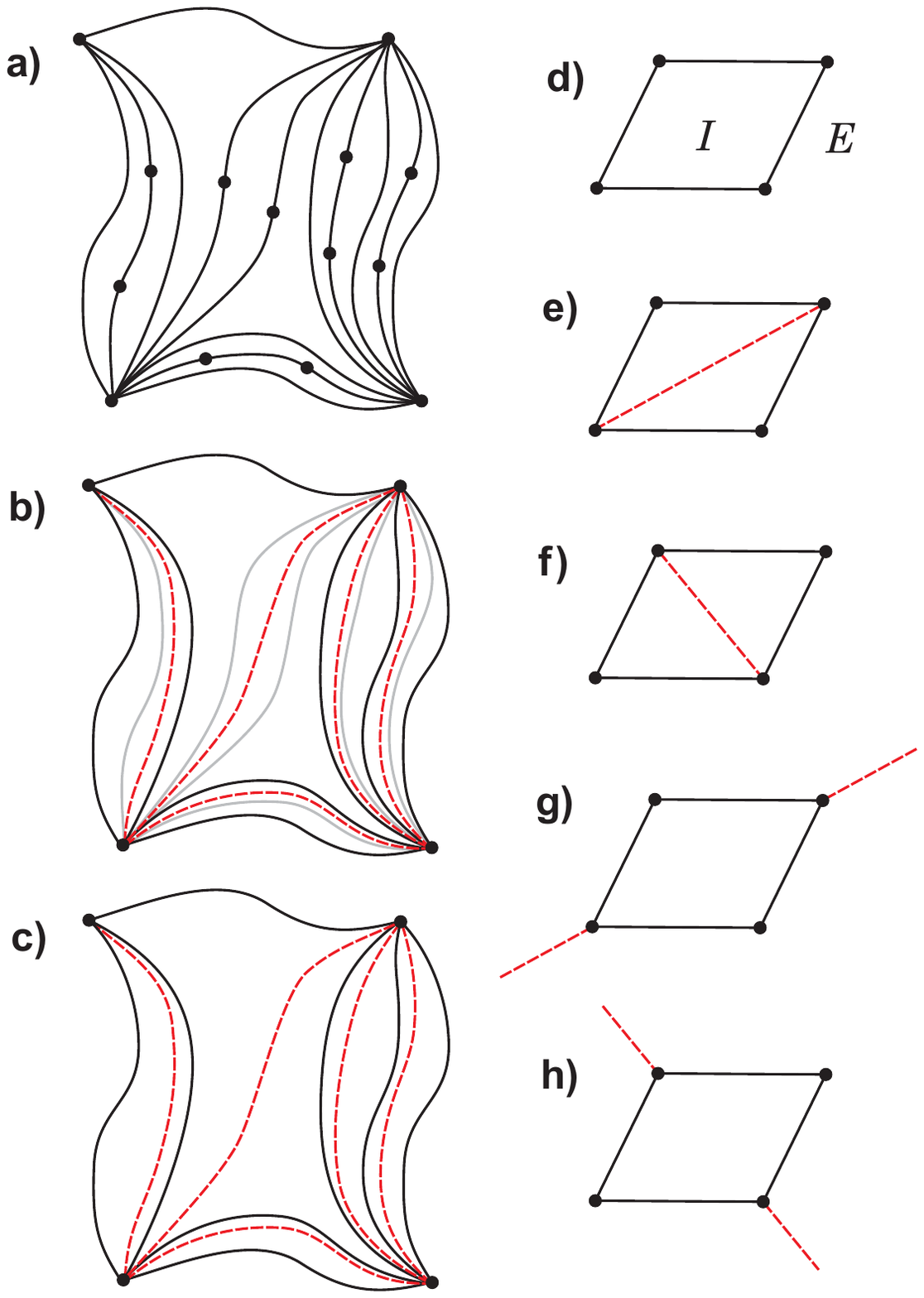}
\caption{To the proof of Lemma \ref{lemma32}.}
\end{center}
\end{figure}
\begin{figure}\label{6}
\begin{center}
\includegraphics*[width=5in]{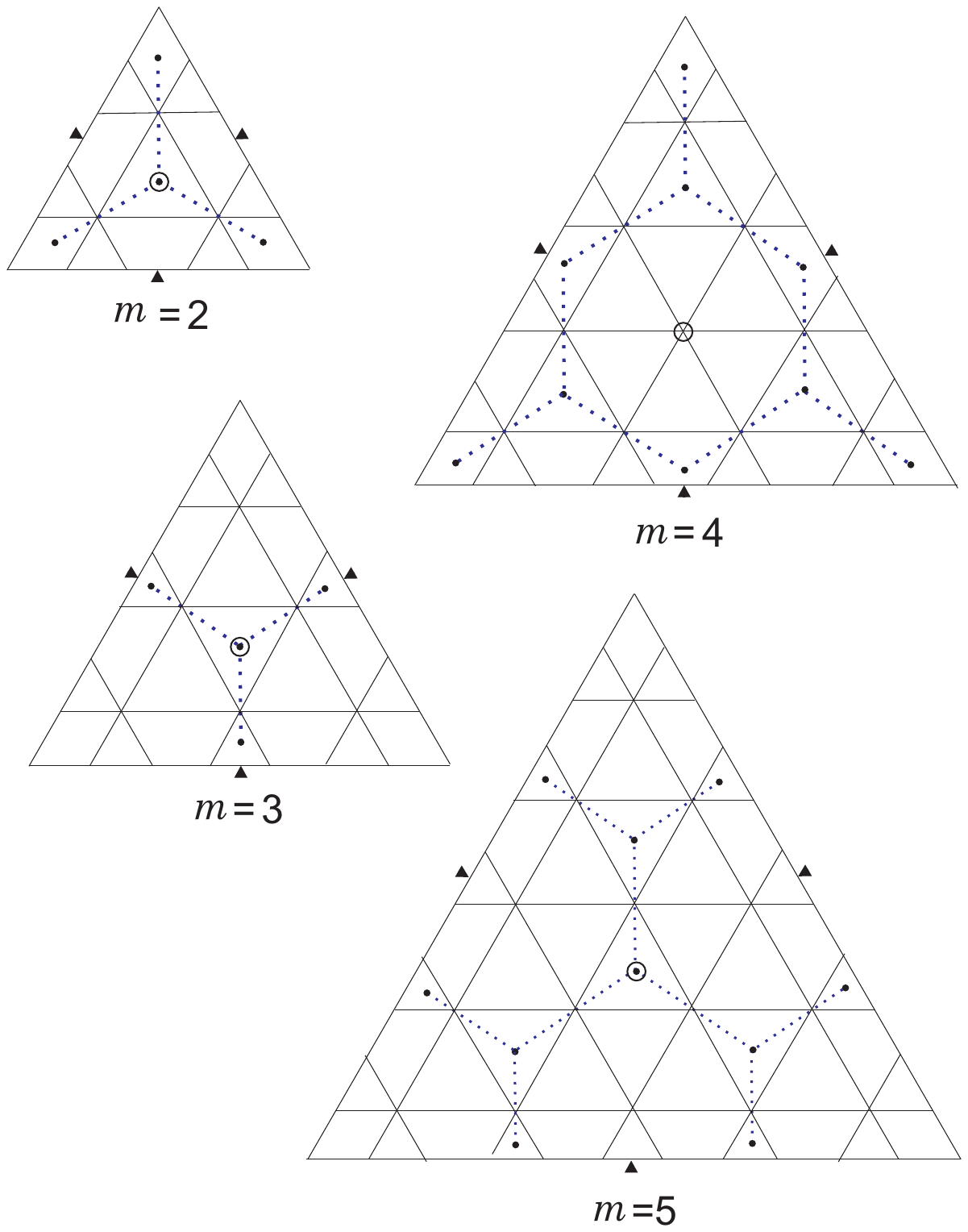}
\caption
{Spaces of angles $\A_m$ for $m\leq 5$. The nerve of component $\L_m^I$
 is blue/dotted.}
\end{center}
\end{figure}
\begin{figure}\label{fig7}
\begin{center}
\includegraphics*[width=5in]{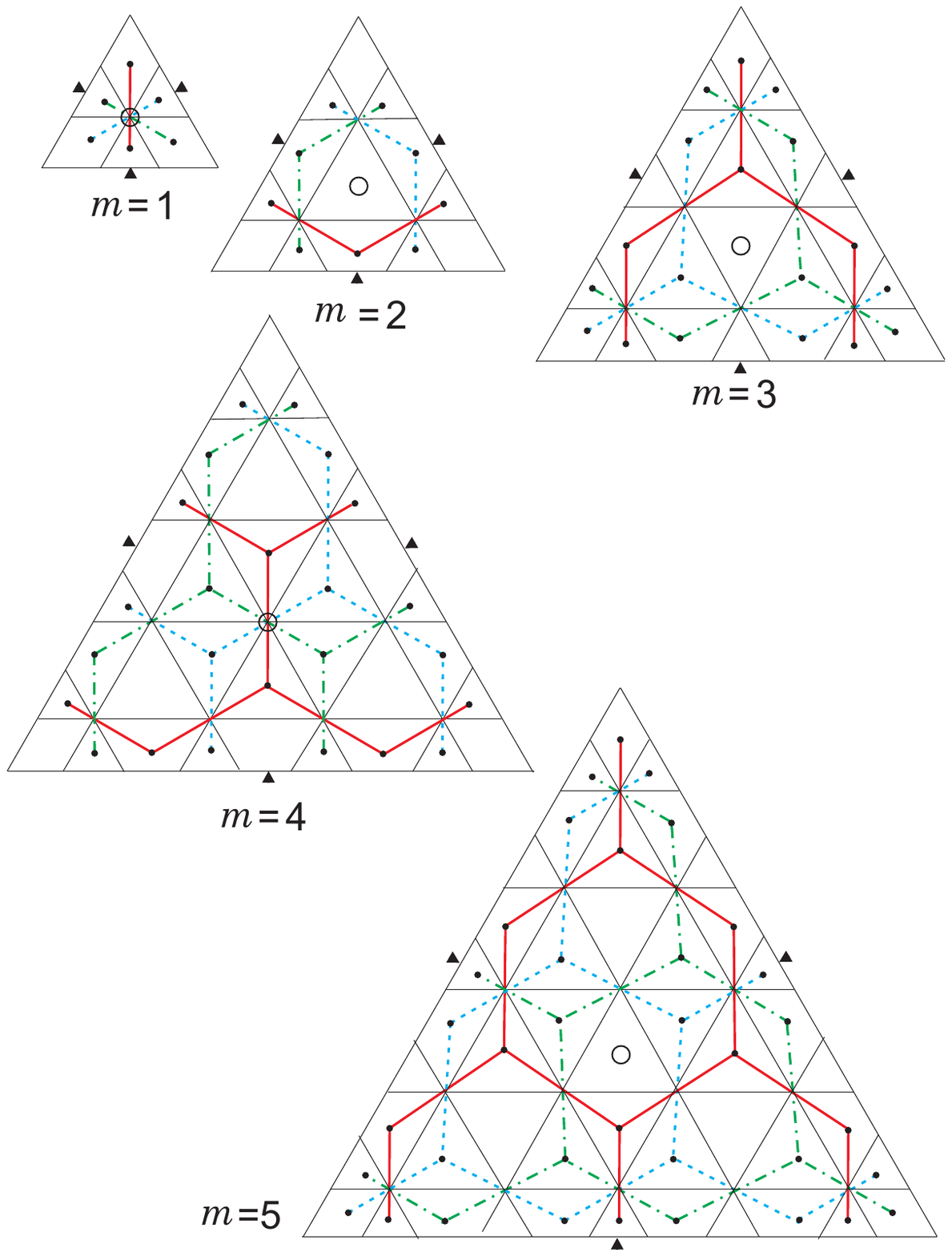}
\caption{Spaces of angles $\A_m$ for $m\leq 5$. The nerves of
components $II_1,II_2,II_3$ are shown.}
\end{center}
\end{figure}
\begin{figure}\label{fig8}
\begin{center}
\includegraphics*[width=5in]{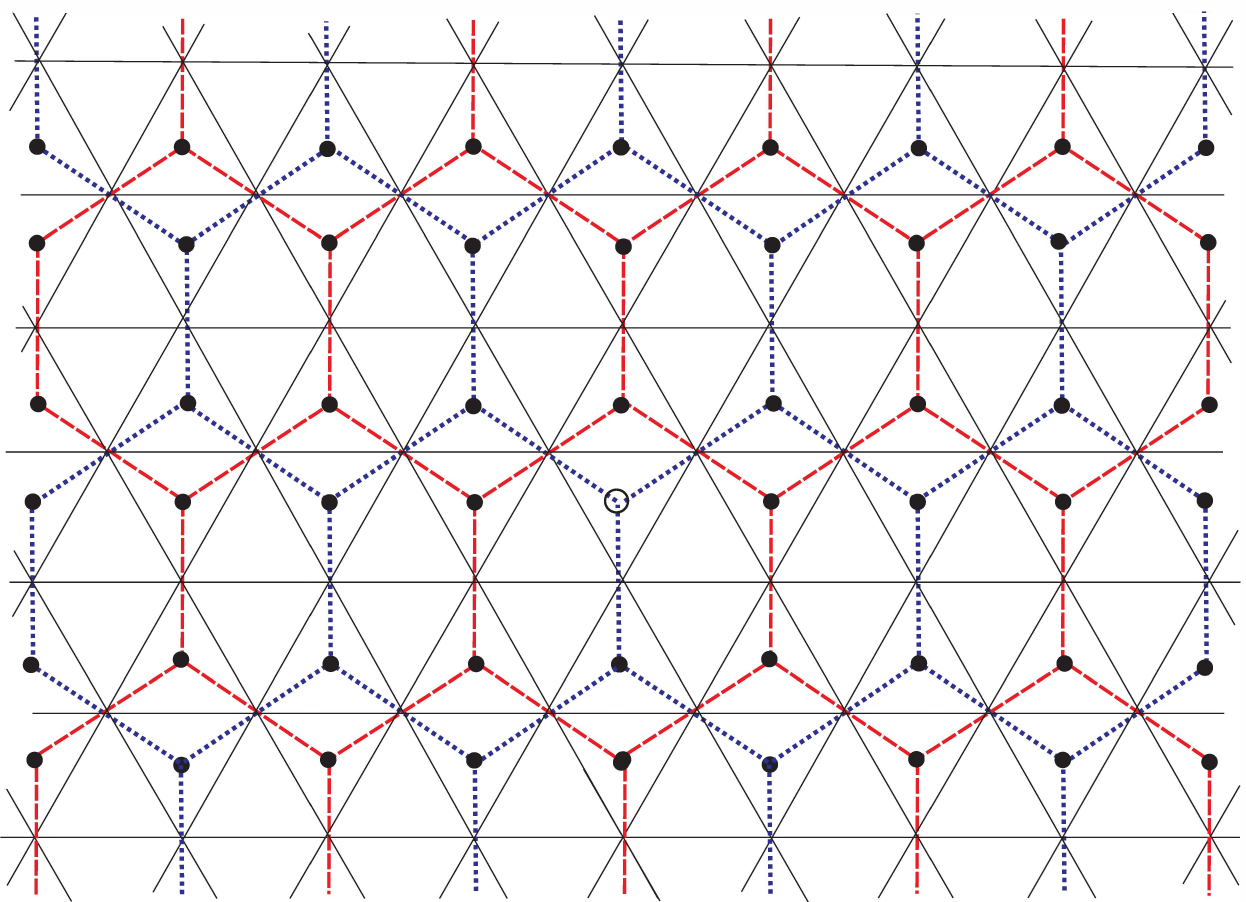}
\caption{Two components in the space of flat singular triangles.}
\end{center}
\end{figure}
\begin{figure}\label{fig9}
\begin{center}
\includegraphics*[width=5in]{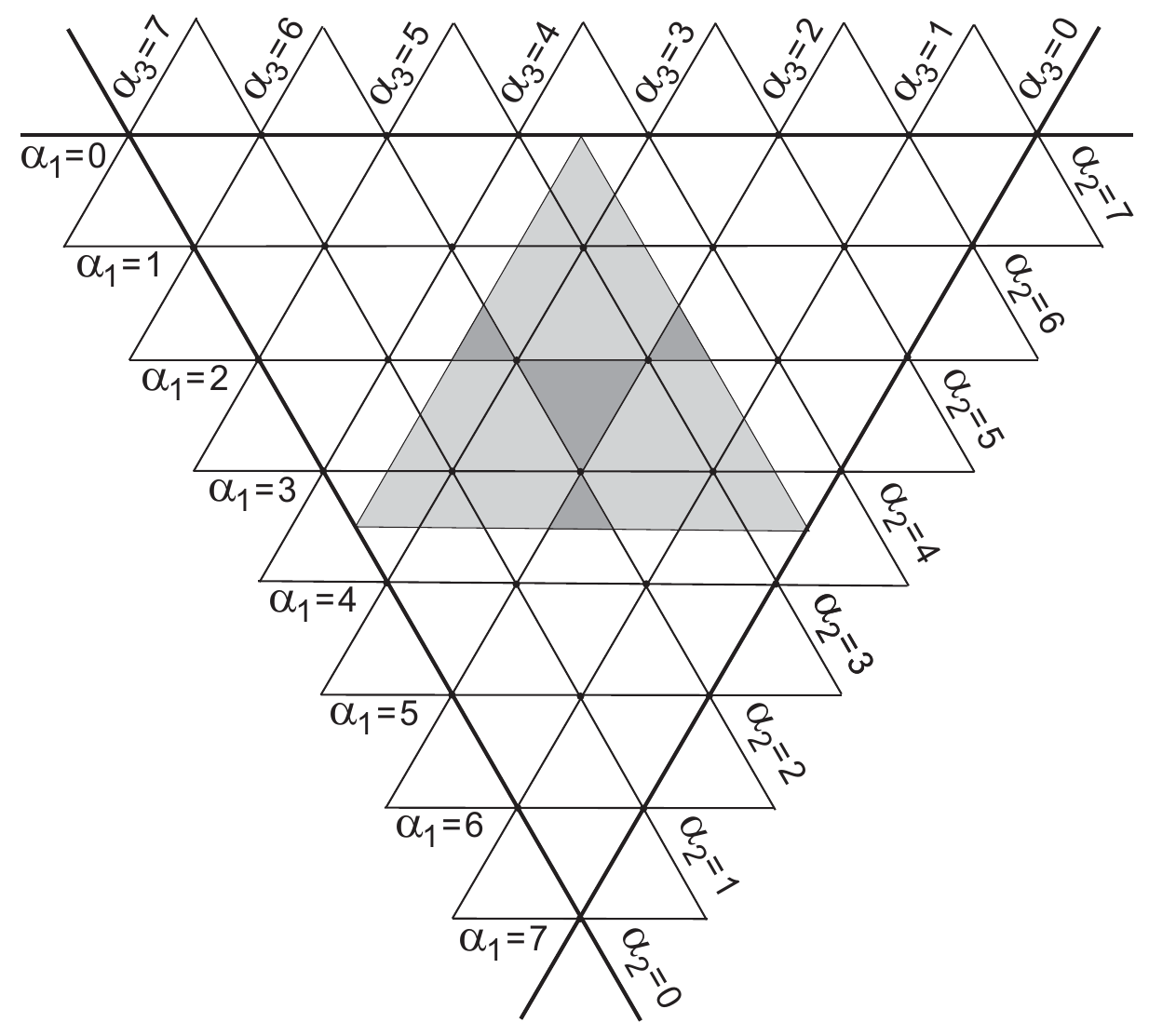}
\caption{Triangle $\Delta_3$ (shaded) in the intersection
of the plane $\alpha_1+\alpha_2+\alpha_3=7$ with the first octant.
Faces of type $I$ have darker shading.}
\end{center}
\end{figure}
\begin{figure}\label{fig10}
\begin{center}
\includegraphics*[width=5in]{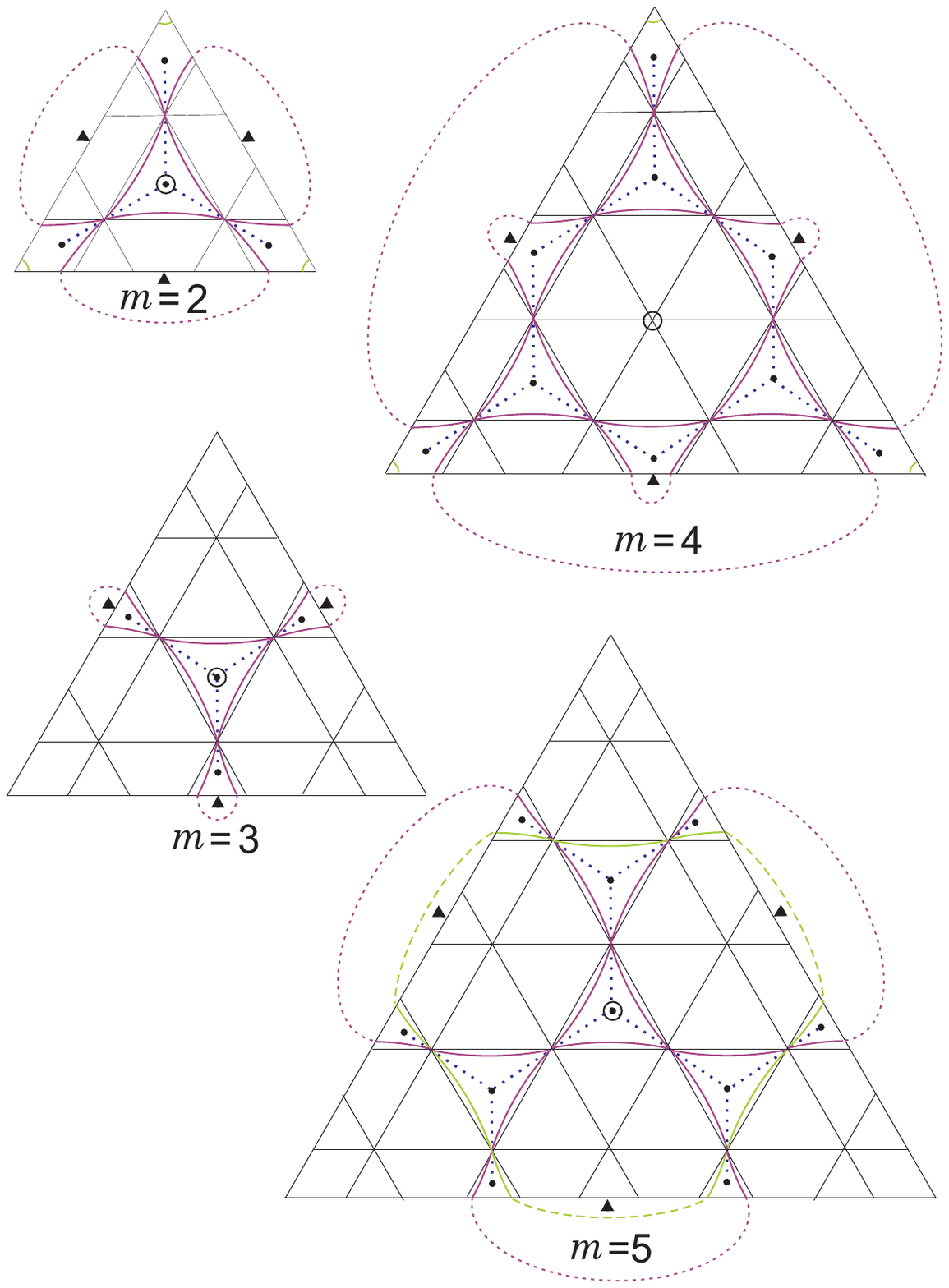}
\caption{Counting punctures for $\L_m^I$.}
\end{center}
\end{figure}
\begin{figure}\label{fig11}
\begin{center}
\includegraphics*[width=5in]{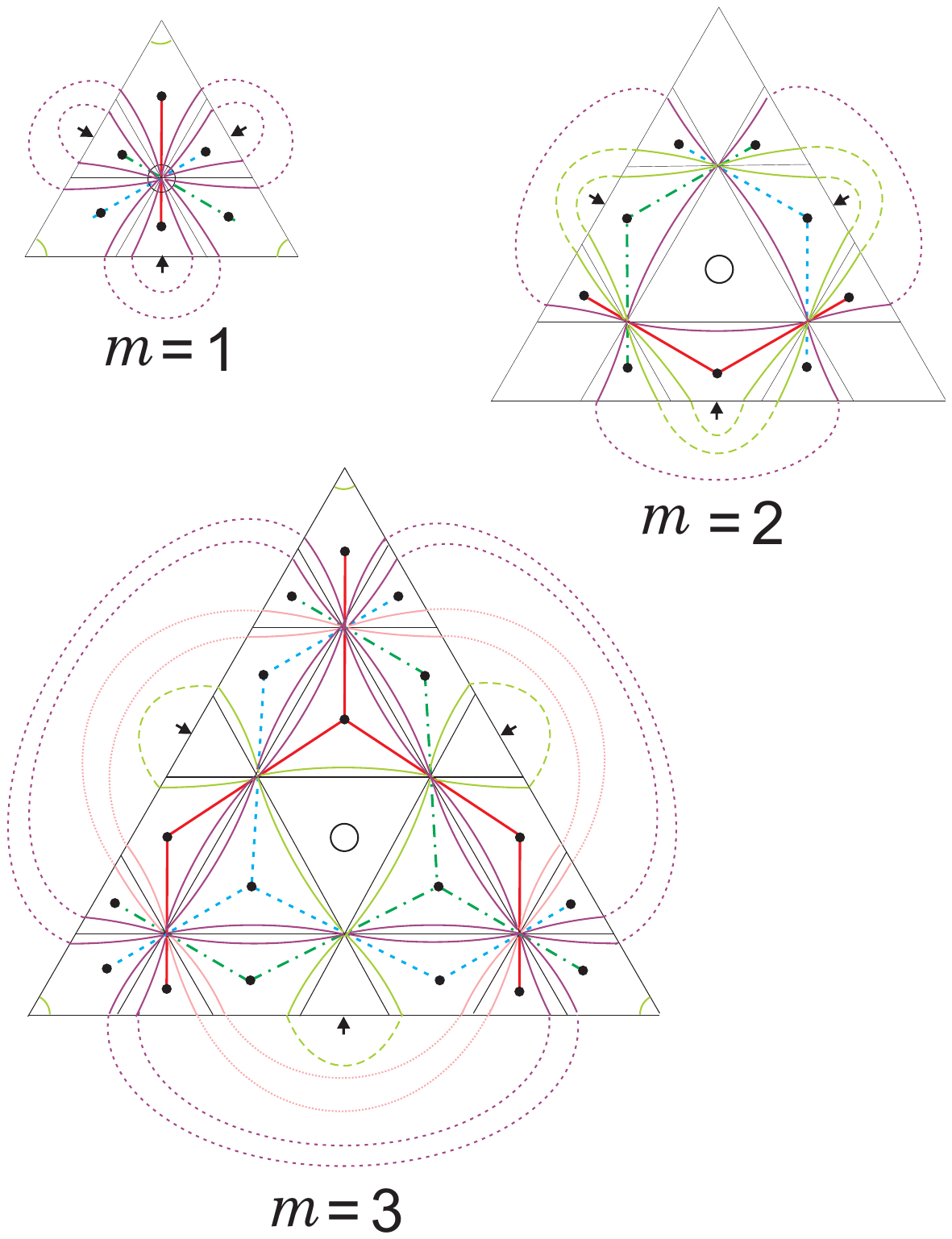}
\caption{Counting punctures for $\L_m^{II}$. Dotted lines represent gluings.}
\end{center}
\end{figure}
\begin{figure}\label{fig12}
\begin{center}
\includegraphics*[width=4.5in]{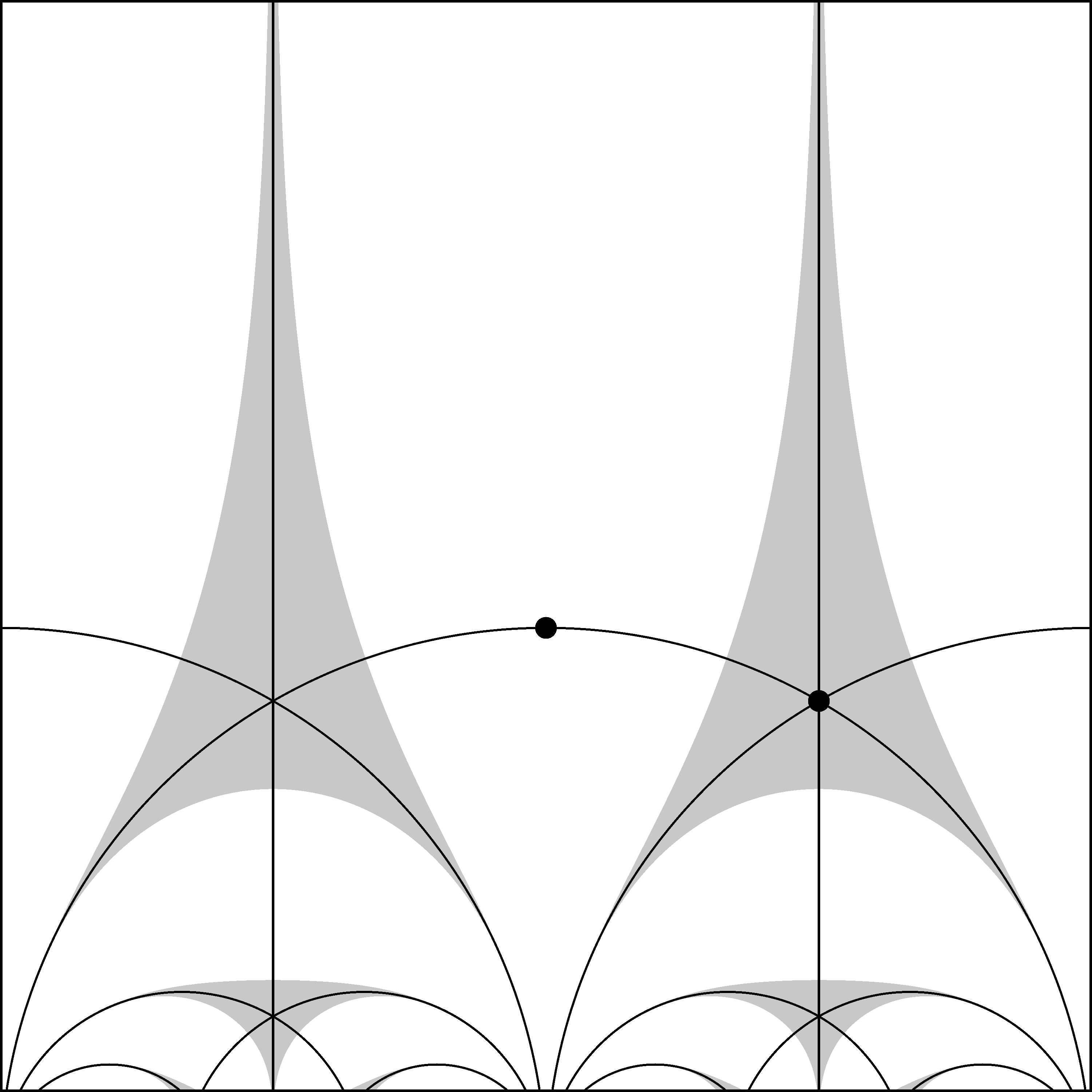}
\caption{Lin--Wang curves for $m=1$, $\tau$-half-plane. Shaded area
corresponds to $\Sph_{1,1}(3)$.}
\end{center}
\end{figure}
\begin{figure}[htb]\label{fig13}
\centering
\includegraphics[width=4.5in]{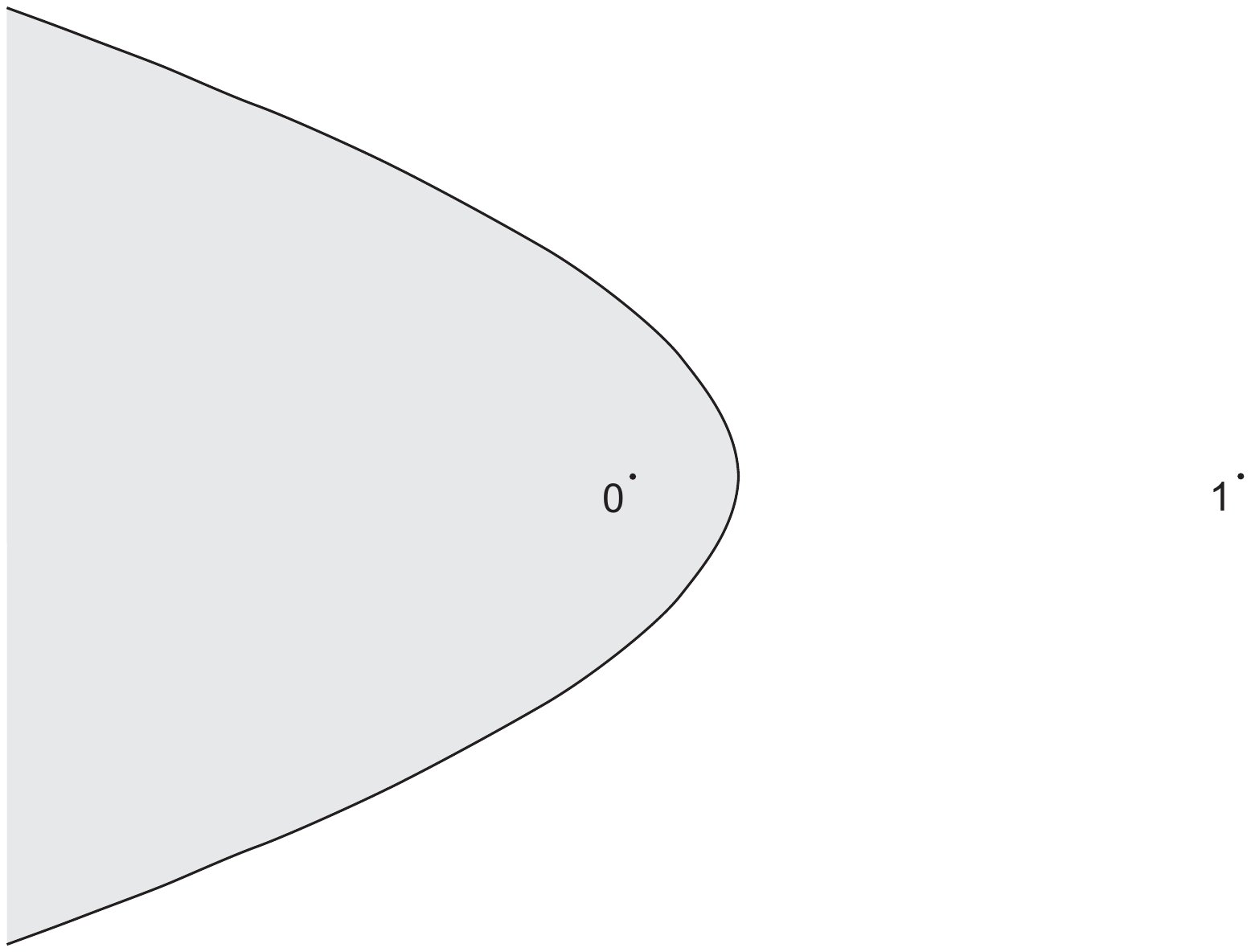}
\caption{Lin--Wang curve $m=1$, $J$-plane. Projection
of $\Sph_{1,1}(3)$ is shaded.}
\end{figure}
\newpage

\begin{figure}\label{fig14}
\begin{center}
\includegraphics*[width=4.5in]{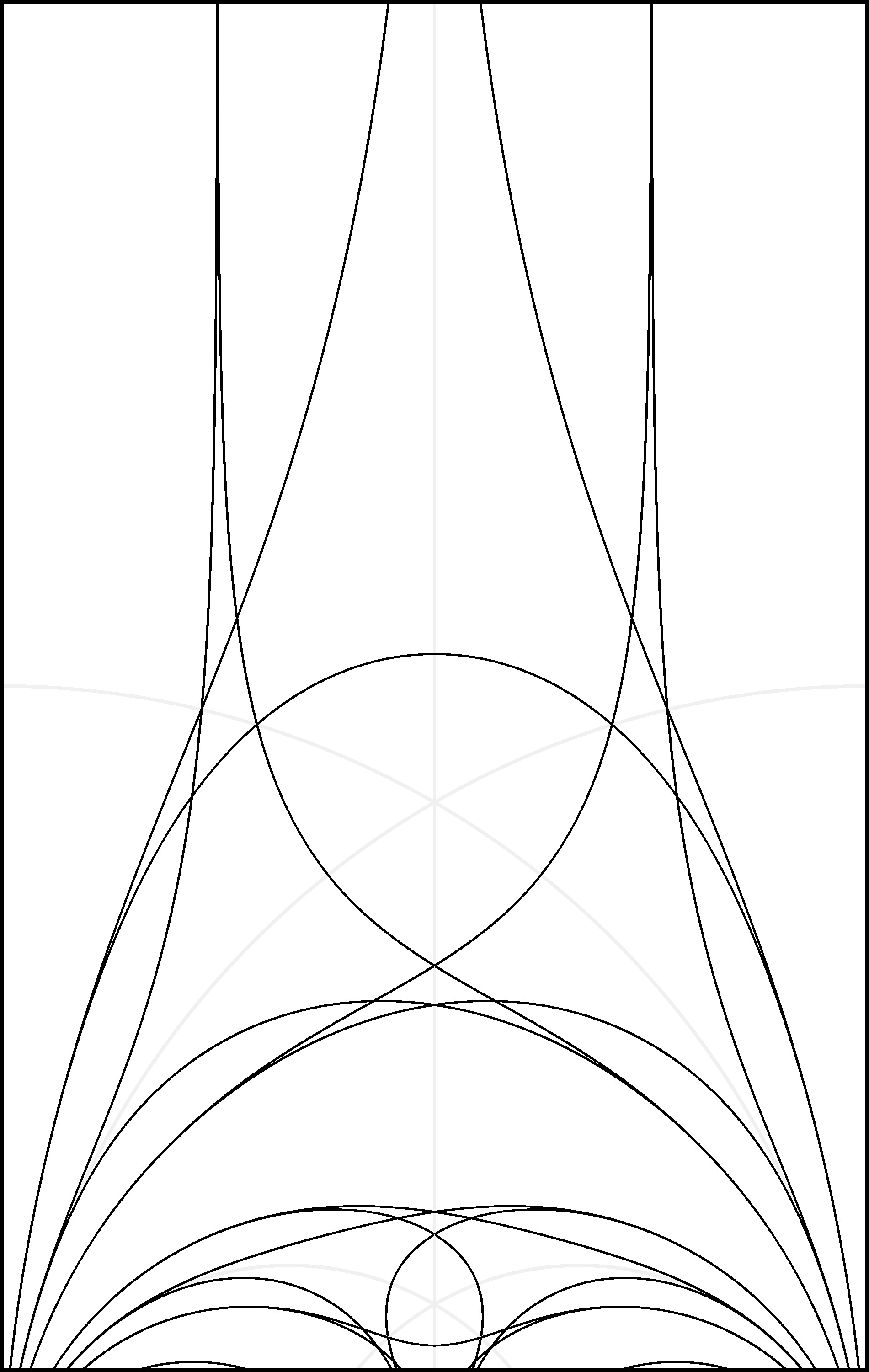}
\caption{Lin--Wang curves for $m=2$, $\tau$-half-plane.}
\end{center}
\end{figure}
\begin{figure}[htb]\label{fig15}
\centering
\frame{
\begin{overpic}[scale=0.6]{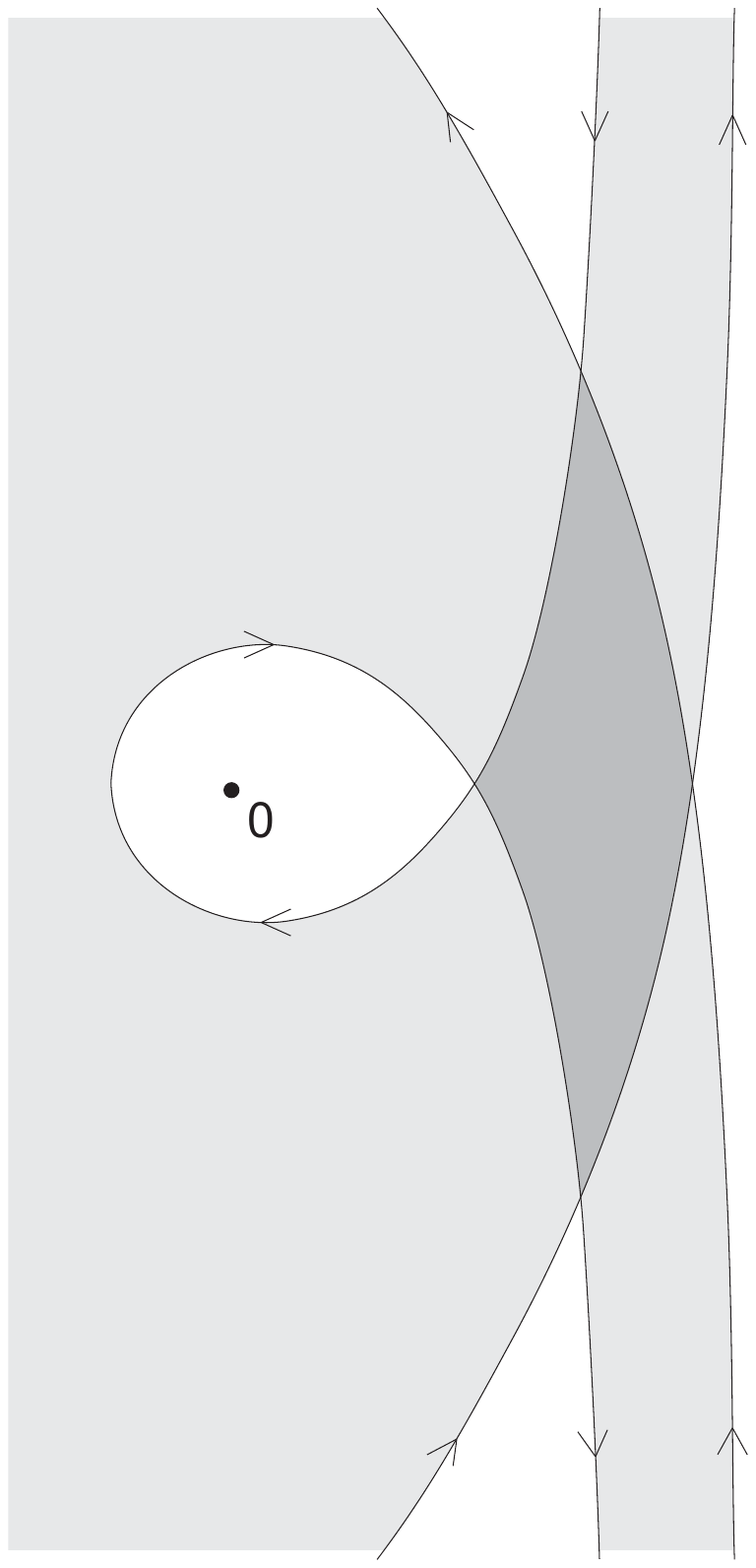}
\end{overpic}
}
\caption{$m=2$, $J$-plane. The curve with a loop is in $\L_2^I$,
the other two curves in $\L_2^{II}$.
Shaded area is the {\em hypothetical} projection
of $\Sph_{1,1}(5)$; it is not known whether the restriction of the forgetful
map on $\Sph_{1,1}(5)$ is open.}
\end{figure}
\begin{figure}\label{fig16}
\begin{center}
\includegraphics*[width=4.5in]{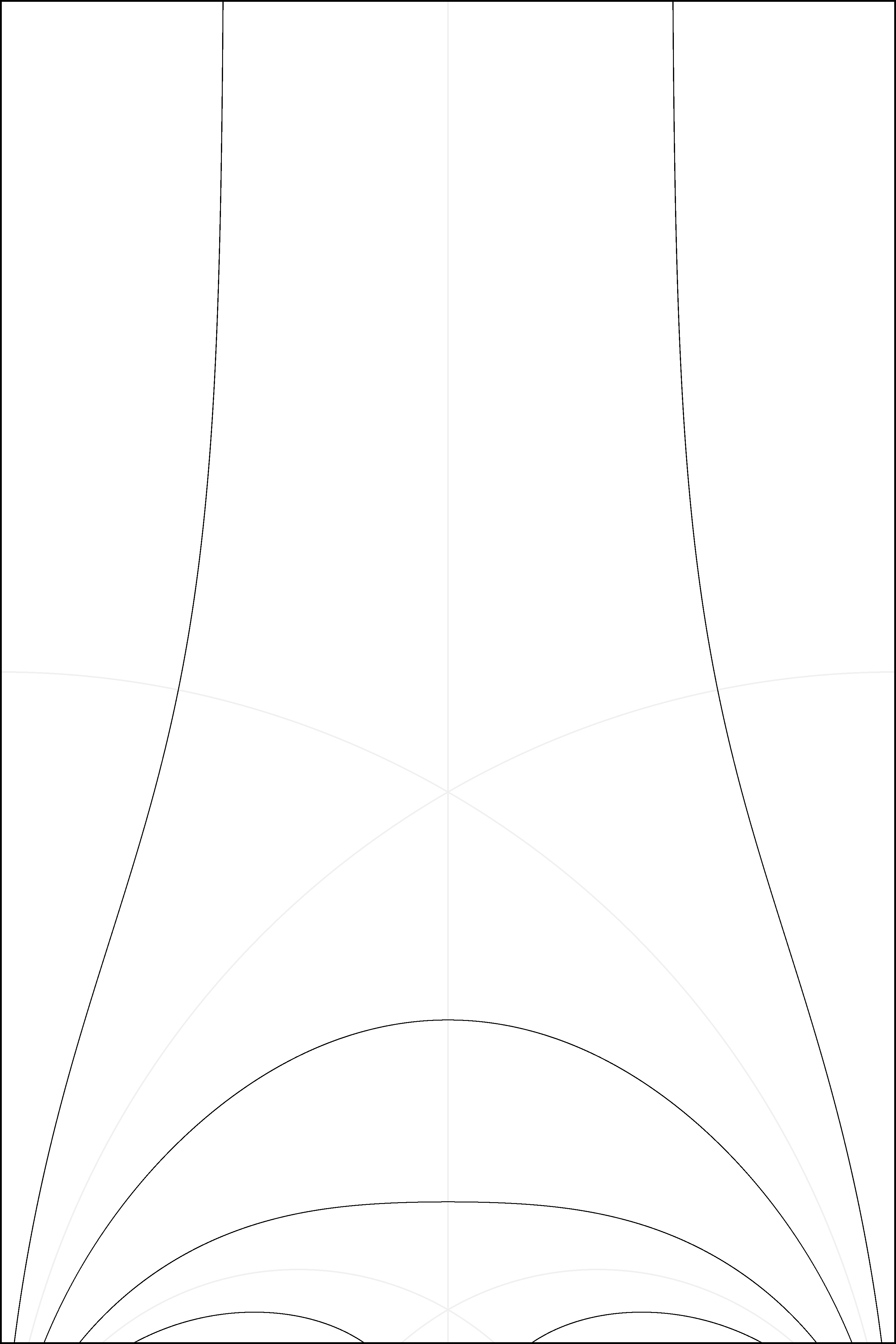}
\caption{$m=3$, $\tau$-half-plane, Lin--Wang curves from $\L_3^I$.}
\end{center}
\end{figure}
\begin{figure}\label{fig17}
\begin{center}
\includegraphics*[width=4.5in]{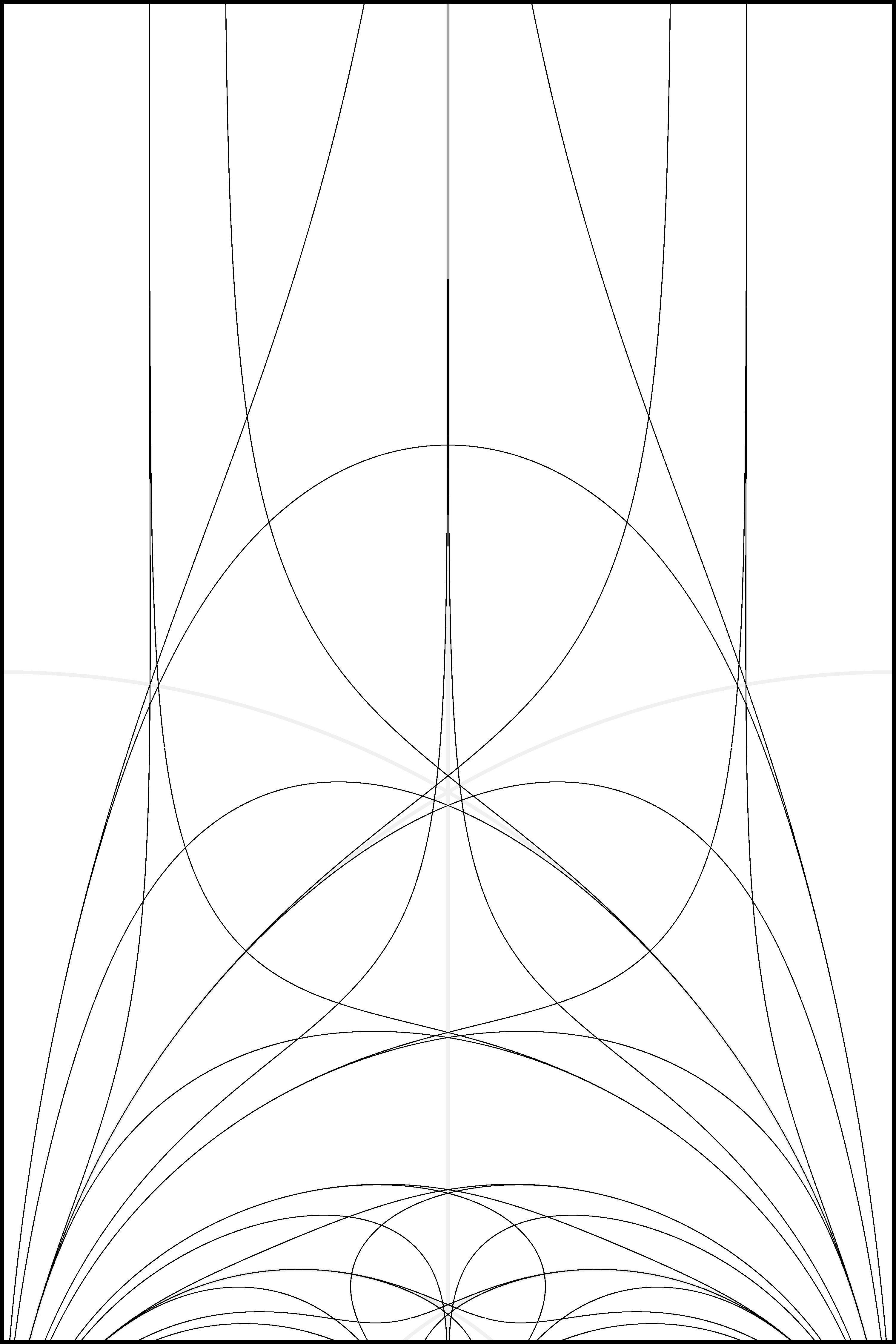}
\caption{$m=3$, $\tau$-half-plane, Lin--Wang curves from $\L_3^{II}$.}
\end{center}
\end{figure}
\begin{figure}[htb]\label{fig18}
\centering
\frame{
\begin{overpic}[scale=0.08]{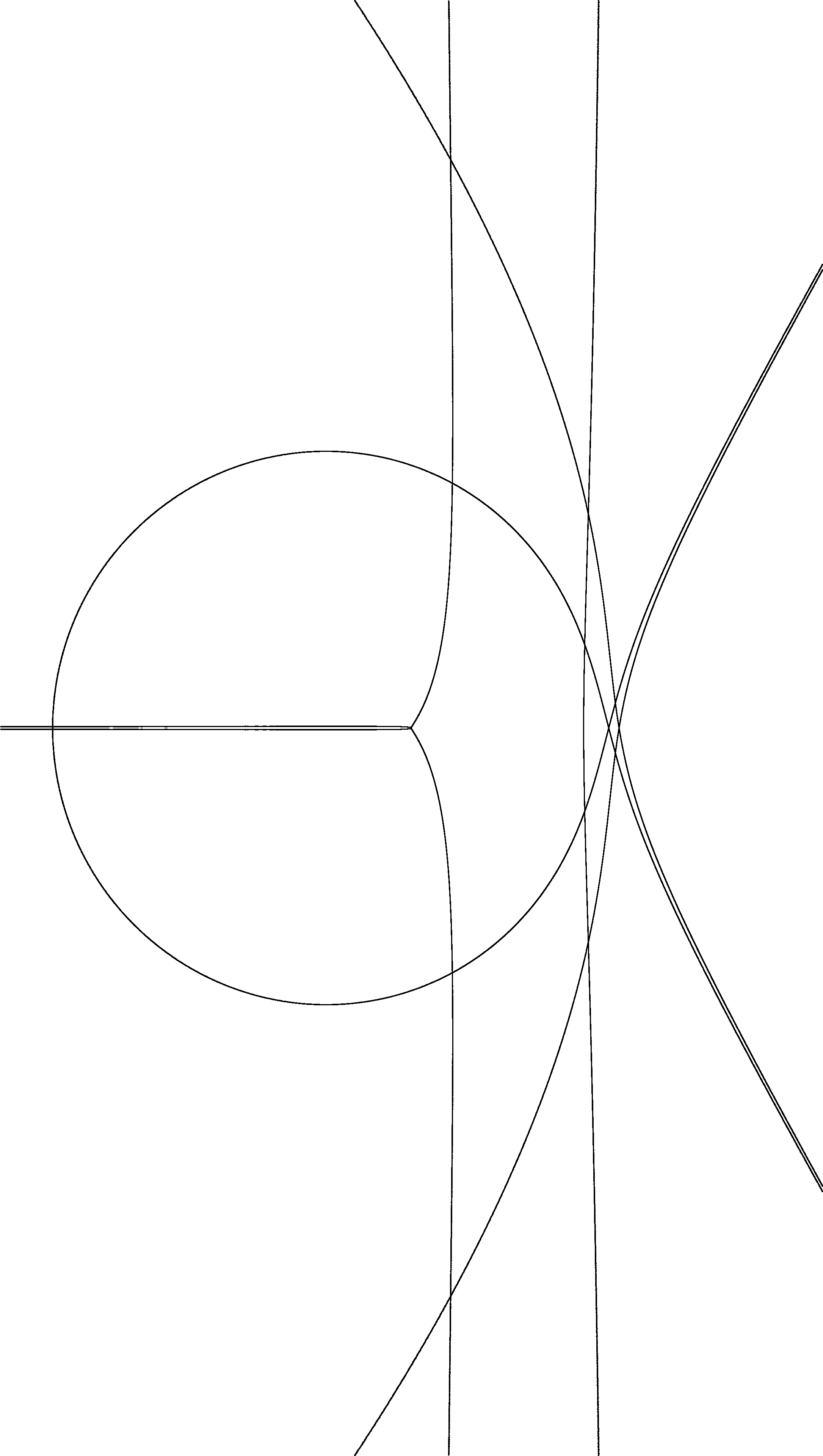}
\put(28.2,50){\circle*{1}}
\put(49.15,50){\circle*{1}}
\put(26.6,47){$0$}
\put(49.3,47){$1$}
\end{overpic}
}
\quad
\frame{
\begin{overpic}[scale=0.08]{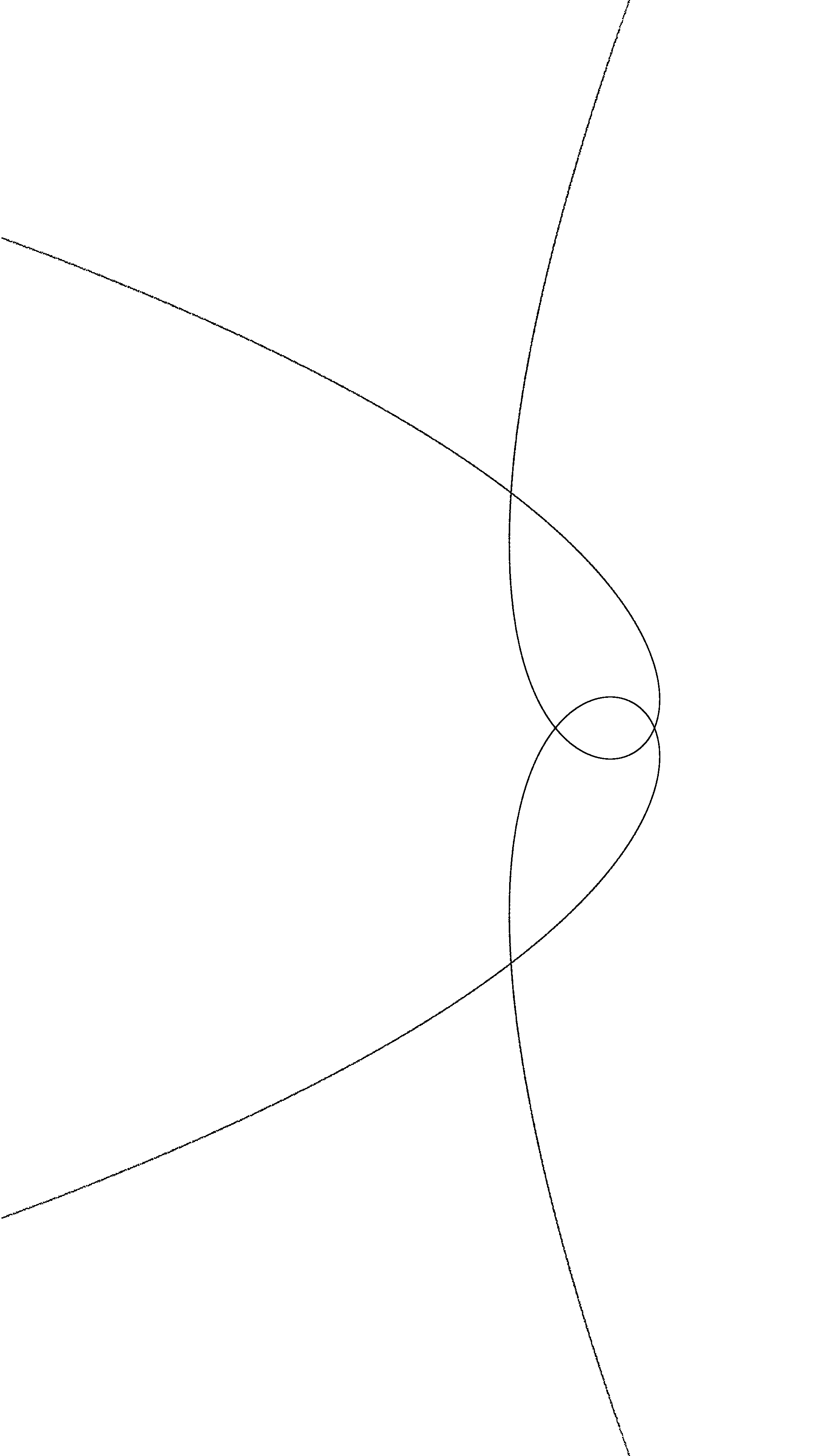}
\put(42.4,50){\circle*{1}}
\put(39.9,48.5){$0$}
\end{overpic}
}
\caption{$m=3$, $J$-plane, Lin--Wang curves from both components.
Magnification of detail on the right, this detail is on component $\L_3^{II}$.}
\end{figure}

\end{document}